\documentclass[11pt]{article}
\usepackage[verbose,tmargin=2.4cm,bmargin=2.4cm,lmargin=2.4cm,rmargin=2.4cm]{geometry}
\usepackage[english]{babel}
\usepackage[utf8]{inputenc}
\usepackage{graphicx}
\usepackage{setspace}
\usepackage{parskip}
\usepackage{float}
\usepackage{booktabs}
\usepackage{mathrsfs}
\usepackage[page]{appendix}

\usepackage{commath}
\usepackage{lmodern}
\usepackage{framed}

\usepackage{titling}       
\usepackage{ragged2e}
\usepackage{subcaption}
\usepackage{siunitx}
\usepackage{caption}
\usepackage[most]{tcolorbox}
\usepackage{placeins}

\usepackage{amsmath,amsfonts,amssymb,amsthm}
\usepackage{mathtools}
\usepackage{physics}

\usepackage{dsfont}
\usepackage[linesnumbered,ruled]{algorithm2e}
\usepackage{booktabs}
\usepackage{scalerel}
\usepackage{tikz}
\usepackage{bm}
\usepackage{amsbsy}
\usetikzlibrary{automata,arrows,positioning,calc}
\usepackage{xcolor}

\usepackage{url}
\usepackage[pdfstartview=FitV, urlcolor=blue]{hyperref}
\usepackage{authblk}

\newcommand{\vx}{\mathbf{x}}
\newcommand{\vy}{\mathbf{y}}
\newcommand{\x}{\mathbf{x}}
\newcommand{\y}{\mathbf{y}}
\newcommand{\yhat}{\hat{\y}}
\newcommand{\veca}{\mathbf{a}}
\newcommand{\ba}{\mathbf{a}}
\newcommand{\bahat}{\hat{\ba}}
\newcommand{\bq}{\mathbf{q}}
\newcommand{\bgamma}{\pmb{\gamma}}
\newcommand{\bv}[1]{\mathbf{#1}}

\newcommand{\TO}{T_{{\scaleto{O}{3pt}}}}
\newcommand{\TR}{T_{{\scaleto{R}{3pt}}}}
\newcommand{\A}{\mathcal{A}}
\newcommand{\M}{\mathcal{M}}
\newcommand{\E}{\mathbb{E}}
\newcommand{\R}{\mathbb{R}}
\newcommand{\Q}{\mathbb{Q}^{M-1}}
\newcommand{\Prob}{\mathbb{P}}
\newcommand{\sscap}[2]{#1^{}_{\scaleto{#2}{3pt}}}
\newcommand{\TG}{\sscap{T}{G}}

\newcommand{\lp}{\left(}
\newcommand{\rp}{\right)}

\newcommand{\update}{\textcolor{black}}

\newtheorem*{theorem*}{Theorem}

\newtheorem*{example*}{Example}

\newtheorem*{fact*}{Fact}
\newtheorem{prop}{Proposition}
\newtheorem*{prop*}{Proposition}
\newtheorem{rem}{Remark}
\newtheorem*{rem*}{Remark}
\newtheorem{assumption}{Assumption}
\newtheorem*{assumption*}{Assumption}

\title{\LARGE \bf
Occasionally Observed Piecewise-deterministic Markov Processes
}
\author[1]{Marissa Gee}
\author[2]{Alexander Vladimirsky}
\affil[1]{Department of Mathematics and Statistics, Kenyon College, Gambier, OH, {\tt gee1@kenyon.edu}}
\affil[2]{Department of Mathematics, Cornell University, Ithaca, NY, {\tt vladimirsky@cornell.edu}}

\date{}

\begin{document}
\maketitle

\begin{abstract}
\noindent
Piecewise-deterministic Markov processes (PDMPs) are often used to model abrupt changes in the global environment or capabilities of a controlled system.
This is typically done by considering a set of ``operating modes'' (each with its own system dynamics  and performance metrics) and assuming that the mode can switch stochastically while the system state evolves.  Such models have a broad range of applications in engineering, economics, manufacturing, robotics, and biological sciences.
Here, we introduce and analyze an ``occasionally observed'' version of mode-switching PDMPs.
We show how such systems can be controlled optimally if the planner is not alerted to mode-switches as they occur but may instead have access to infrequent mode observations.  We first develop a general framework for handling this through dynamic programming on a higher-dimensional mode-belief space.
While quite general, this method is rarely practical due to the curse of dimensionality.   
We then discuss assumptions that allow for solving the same problem much more efficiently, with the computational costs growing linearly (rather than exponentially) with the number of modes.  We use this approach to derive Hamilton-Jacobi-Bellman PDEs and quasi-variational inequalities encoding the optimal behavior for a variety of planning horizons (fixed, infinite, indefinite, random) and mode-observation schemes (at fixed times or on-demand).
We discuss the computational challenges associated with each version and illustrate the resulting methods on test problems from surveillance-evading path planning. 
We also include an example based on robotic navigation: a Mars rover that minimizes the expected time to target while accounting for the possibility of unobserved/incremental 
damages and dynamics-altering breakdowns.

\end{abstract}


\section{Introduction} \label{sec:oopdmps-introduction}

Piecewise-deterministic Markov processes (PDMPs) \cite{davis1984pdmps} 
provide an excellent framework for modeling non-diffusive stochastic processes, in which deterministic dynamics are punctuated by random jumps.
PDMPs are also very useful for modeling abrupt changes in a global environment with a set of operating {\em modes} $\mathcal{M} = \{1, 2, \hdots, M\},$ 
where each mode specifies its own system dynamics and performance metrics while the mode-to-mode switches are governed by a known continuous-in-time Markov chain.
Such processes (and the task of controlling them optimally) arise naturally in many applications including optimizing production rates in failure-prone manufacturing systems \cite{akella1986manufacturing}, determining optimal premium levels for insurance plans \cite{schal1998insurance}, 
coexistence modeling in population ecology \cite{benaim2016lotka, costa2016piecewise, HeningStrickler2019}, 
and fish harvesting in an environment fluctuating due to El Niño \cite{cartee2023UQ}. 
PDMPs are also very relevant in many path-planning applications; e.g., a planetary rover might often change its trajectory in case of a partial breakdown 
that diminishes its speed on rough terrain.  Notably, if the rate of breakdowns is known and possibly dependent on terrain properties, this risk 
will affect the choice of trajectories even before that switch (to a ``partially broken'' mode) happens \cite{gee2022breakdowns}.  
In the same way, an animal deciding on a path to take to its feeding grounds will take into account the (spatially dependent) risk of predation even before any predators actually spot it \cite{rands2017hesitate, LOF}. 
Similarly, an evader trying to minimize the risk of detection will try to hedge its path-choices against the future stochastic changes in surveillance patterns.

In all these cases, the classical assumption in controlling PDMPs is that any mode switch becomes immediately known to the controller, who can alter their plans accordingly.  But in many applications this assumption is quite far from reality.  An accumulation of hardware/software issues might switch a planetary rover to an intermediate ``breakdown-prone'' mode without affecting the immediately observed dynamics.  A controller will not know about this switch until they choose to run diagnostics, which would normally incur an additional time penalty.  Whether and when to run such diagnostics (as well as how to modify the path until then) becomes a crucial question.
Similarly, an animal may not immediately know that it is being stalked by a predator unless it takes a break from foraging to look around.
In the same way, a surveillance evader typically will not have immediate knowledge of a change in surveillance patterns.  In all these examples, 
the current mode is observed only occasionally and the main challenge is how to govern the system in between observations.  
The purpose of this paper is to derive the PDEs defining the optimal behavior for a subclass of such Occasionally Observed Piecewise-Deterministic Markov Processes (OOPDMPs)
and develop efficient numerical methods for them.

A general approach for controlling ``partially-observed'' processes involves maintaining a {\em belief} (i.e., a probability distribution over the set of possible system states),
updating that belief whenever partial observations become available, and using dynamic programming on a {\em space of beliefs} to find the optimal control policy \cite{kaelbling1998pomdps}.
While flexible and powerful, this approach is computationally expensive and often infeasible for many applications.  In the setting of OOPDMPs, it would result in algorithms whose computational and storage costs grow exponentially with \update{the number of modes} $M.$  We circumvent this difficulty by focusing on problems in which the evolution of beliefs can be fully computed using the last mode observation and the time since it was observed.  While restrictive, this assumption allows us to derive algorithms on a much lower dimensional state space.
\update{Beyond computational costs, we also note a clear difference in types of addressed problems:  the usual ``partially-observed'' framework is most useful when state observations are noisy and partial but frequent, while the dynamics and costs are subject to nearly continuous random perturbations.  
In contrast, our approach is primarily suitable for problems, in which the controlled continuous dynamics are fully observed, performance metrics 
remain predictable in each mode, mode switches are discrete and not detectable, and mode observations are exact though infrequent.} 

We start by reviewing the control of fully-observed PDMPS and existing belief-state programming approaches in Section \ref{sec:background}.
We then present the OOPDMP framework in detail in Section \ref{sec:models}.
The structure of the governing PDEs that we derive for these models is heavily dependent on two key features:  (a) when the process stops (e.g., at a fixed or random time, upon reaching a target, or never) and (b) when and how the mode observations become available (e.g., never, at predetermined
times, or upon request -- with a bound on the number of requested observations or with a fixed cost attached to each of them).  
This presentation is followed by an overview of the numerical solution techniques in Section \ref{sec:numerics}.
The method is then illustrated in Section \ref{sec:numexp} using computational examples of surveillance evasion and planetary rover path planning.
We conclude by discussing the limitations and directions for future work in Section \ref{sec:conclusion}.

\section{Background}
\label{sec:background}
\subsection{Fully Observed Piecewise-deterministic Markov Processes}\label{sec:pdmps}
In deterministic optimal control problems, we are generally interested in a process with controlled dynamics
\begin{align}
\label{eq:generic_dynamics}
\vy'(s) &= \vb f\left(\vy(s), s, \veca \right)\\ 
\vy(t) &= \vx
\end{align}
where $\vy(s) \in \Omega \subset \R^d$ is the state of the system at the time $s \geq t$ 
and a control parameter $\veca$ takes values in some compact set $A$.
The goal is usually to minimize the total cost up to the termination, comprised of the running cost $K$ integrated along the trajectory and the terminal cost $\psi.$

In general piecewise-deterministic Markov processes (PDMPs) \cite{davis1984pdmps}, the continuous evolution of the state described by \eqref{eq:generic_dynamics} is punctuated by random discrete jumps, taking $\y$ 
to new points in $\Omega$.  The process description is augmented by a (state-dependent) jump intensity and transition kernel (a probability measure over possible post-jump states).
However, we are interested in applications where ``jumps'' are only used to represent abrupt changes in a global operating environment, controlled system's capabilities, or performance measures.  This allows for a simplified description of the relevant PDMPs, with a special structure that will be exploited throughout this paper. 
Instead of jumps in $\y(s),$ we assume that switches occur in
an additional discrete state variable $\mu(s)$ which takes values in some finite set $\mathcal{M} = \{1, 2, \hdots, M\}$ and which we refer to as the \emph{operating mode} or simply the \emph{mode} of the problem.
Transitions between modes are then stochastic with known, possibly state- and time-dependent rates $\lambda_{ij}(\vx,t)$.
Thus, at a some state and mode $(\vy(s), \mu(s))$, the likelihood of a mode transition over some time $\tau$ is given by
\begin{equation}
\Prob\left(\mu(s + \tau) = j \,|\, \mu(s) = i\right) = \lambda_{ij}\lp \vy(s), s\rp \tau + o(\tau).
\end{equation}
Within a single mode $i \in \M$, the system now evolves according to the deterministic dynamics
\begin{align}
\vy'(s) &= \vb f_i\left(\vy(s), s, \veca \right).
\end{align}

In controlled PDMPs of this type,
the objective is to minimize the {\em expected} cumulative cost 
computed using mode-dependent running costs $K_i(\vx,t,\veca)$ and terminal costs $\psi_{i}(\vx).$ 
Throughout the paper, we make the following assumptions 
about
$\Omega,$ $\vb f_i,$ $K_i,$ $\psi_i,$ and $\lambda_{ij}.$
\begin{assumption}\label{assump:regularity}
For all $i,j \in \M$:
\begin{enumerate}
\item The domain $\Omega \subset \R^d$ is open, bounded, connected, and has a piecewise-smooth boundary.
\item Whenever the problem formulation includes a target set $\Gamma \subset \Omega,$\\ that $\Gamma$ is closed and has a non-empty interior.
\item Each $\vb f_i$ is bounded and Lipschitz-continuous.
\item All $K_i$, $\psi_i$, and $\lambda_{ij}$ are bounded, nonnegative, and uniformly continuous.
\item For every $\x \in \bar{\Omega}$ and $t \in \R,$ all sets $\left\{  \left( \vb f_i(\x, t, \ba), K_i(\x,t,\ba) \right) \, \mid \, \ba \in A \right\}$ are convex.  
\end{enumerate}
\end{assumption}

We now give a brief overview of optimal control of \emph{fully observed} PDMPs for a variety of planning horizons before examining the cases with limited information.
 
$\bullet$
In \emph{finite horizon} (also referred to as \emph{fixed horizon}) problems, the process continues until some known and pre-specified finite time $T$.
Given a policy $\veca(\vx, t, i)$ that specifies a control  
as a function of state, time, and mode, 
if we start at the time $t$ from $(\vx = \vy(t), i = \mu(t)),$ 
the \emph{expected} cost-to-go can be defined as
\begin{align}
J_i(\vx,t,\veca(\cdot)) &= \E\left[ \int_{t}^{ T } K_{\mu(s)}\left(\vy(s),s,\veca(\vy(s), s, \mu(s))\right) \, ds + \psi_{\mu(T)}\left(\vy(T)\right) \right].
\end{align}
Given the initial mode $i$,  the current mode $\mu(s)$ by the time $s \in (t,T]$ is a random variable and the expectation is taken with respect to
a random sequence of mode switches.  However, the controller learns about each mode transition immediately and adjusts the control accordingly
(since $\veca(\cdot)$ is a control policy in feedback form and $\mu(s)$ is among its arguments).

We can now define the value function, which encodes the \emph{optimal} expected cost-to-go for every possible state, time, and mode, as
\begin{equation}
       u_i(\vx,t) = \inf_{\veca(\cdot) \in \A} J_i(\vx,t,\veca(\cdot)) \label{eq:v-finite-known-modes}
\end{equation}
where $\A$ is the set of measurable functions from $\Omega\times [0,T] \times \mathcal{M} $ to $A$.
We note that Assumption \ref{assump:regularity}.4 can be used to prove the existence of optimal controls  
and the above infimum is actually attained.
A standard argument \cite{davis1999viscosity}
shows that these value functions $u_i(\vx,t)$ are the unique viscosity solutions\footnote{
In general, whenever we refer to value functions as ``satisfying'' PDEs we mean that they are $\Omega$-constrained viscosity solutions \cite{bardi2008, davis1999viscosity}.
The same also applies to boundary conditions, which are meant to be satisfied ``in viscosity sense''; see  \cite[Chapter 5]{bardi2008}.
Throughout the paper,  
this interpretation yields value functions defined by minimizing over only those controls that guarantee that the process does not leave 
$\bar{\Omega}$ 
though might be moving for some time along $\partial \Omega$.
}
to a weakly-coupled system of Hamilton-Jacobi-Bellman (HJB) equations of the form 
\begin{align}
         - \frac{\partial u_i}{\partial t} &= H_i \lp \vx, t, \nabla u_i(\vx, t) \rp + \sum_{j \neq i} \lambda_{ij}(\vx) \, \left(u_j(\vx, t) - u_i(\vx, t)\right), \quad t \in [0,T), \label{eq:fin-pde}\\
        u_i(\vx,T) &= \psi_i(\vx), \label{eq:fin-tc}
\end{align}
where $i=1,...,M,$ and we have a mode-dependent Hamiltonian $H_i$ defined by
\begin{equation}
H_i\lp\vx, t, \nabla u_i(\vx, t)\rp = \min_{\veca \in A}\left\{K_i(\vx,t,\veca) + \mathbf{f}_i(\vx, t, \veca) \cdot \nabla u_i(\vx, t) \right\}. \label{eq:fin-hamil}
\end{equation}
This system can be solved numerically backwards in time using standard discretization schemes for HJB equations.

In all other versions of the problem considered below, we will assume that 
the dynamics $\mathbf{f}_i,$ costs $(K_i, \psi_i),$ and switching rates $\lambda_{ij}$ are autonomous (i.e., not dependent directly on time).
We note that, for the finite horizon version above, the value functions $u_i$ would still remain $t$-dependent even with this assumption, 
since it is important to know the time $(T-t)$ remaining until the end of the process before the optimal control value can be selected.
In contrast, in all of the following versions of the problem, the value functions $u_i(\vx)$ will be stationary.
We will slightly abuse the notation by using 
$H_i\lp \vx, \nabla u_i(\vx)\rp$  
to denote the corresponding time-independent Hamiltonians.

$\bullet$ In \emph{infinite-horizon} problems, the process never terminates, but the running cost is discounted with time.
Since the starting time is now irrelevant, we can define
the expected overall cost as
\begin{equation}
J_i(\vx, \veca(\cdot)) = \E \left[ \int_0^\infty e^{-\beta s} K_{\mu(s)} \left(\vy(s), \veca \left( \vy(s), \mu(s) \right) \right) ds \right] \label{eq:inf-full-cost}
\end{equation}
where $\beta > 0$ is a discount factor. The value functions
satisfy a system of {\em stationary} PDEs
\begin{align}
0 &= H_i\lp \vx, \nabla u_i(\vx)\rp - \beta u_i(\vx) + \sum_{j \neq i} \lambda_{ij}(\vx) \, \left(u_j(\vx) - u_i(\vx)\right). \label{eq:inf-full-pde}
\end{align}

$\bullet$ In \emph{indefinite horizon} (also known as \emph{exit-time}) problems,
there is no time-discounting and
the process terminates once it reaches (enters the interior of) a target set $\Gamma \subset \Omega$.  As a result, the planning horizon 
$\TG = \inf \left\{t \mid \vy(t) \in \Gamma \backslash \partial \Gamma \right\}$
becomes a random variable  dependent on the initial state $(\vx,i),$ the selected control policy $\veca,$ and the sequence of random mode-switches.  
The corresponding value functions $u_i(\vx)$
satisfy a system of PDEs
\begin{align}
          0 &= H_i(\vx, \nabla u_i(\vx)) + \sum_{j \neq i} \lambda_{ij}(\vx) \, \left(u_j(\vx) - u_i(\vx)\right), & \vx &\in \Omega \backslash \Gamma; \label{eq:indef-full-pde}\\
        u_i(\vx) &= \psi_i(\vx), \quad &\vx &\in \partial \Gamma\label{eq:indef-full-bc}.
\end{align}

$\bullet$ Finally, in \emph{randomly-terminated} problems, the termination is viewed as a result of a non-homogeneous Poisson process,
where the rate of termination $\gamma$ might be also mode and/or location dependent.
The value functions can be found as viscosity solutions of
the stationary system
\begin{align}
0 &= H_i\lp \vx, \nabla u_i(\vx) \rp + \gamma_i(\vx) \lp \psi_i(\vx) - u_i(\vx) \rp + \sum_{j \neq i} \lambda_{ij}(\vx) \, \left(u_j(\vx) - u_i(\vx)\right). \label{eq:rt-full-pde}
\end{align}

These stationary PDE systems can be computationally challenging: due to the coupling between different modes, 
fast (non-iterative) numerical methods are not usable even if they were applicable for the mode-decoupled version (with all $\lambda_{ij} = 0$).
The systems \eqref{eq:inf-full-pde}, \eqref{eq:indef-full-pde}, and \eqref{eq:rt-full-pde} are generally solved using iterative methods; e.g., by looping through all modes and finding a numerical approximation for the current $u_i$  while temporarily holding all other $u_j$'s fixed, repeating this process until the value changes fall below a specified threshold \cite{gee2022breakdowns}.

\renewcommand{\arraystretch}{1.3}
\begin{table}
\update{
    \centering
    \begin{tabular}{| c | p{11.5cm} | m{1.3cm}|}
    \hline
     Notation & Description  & First used in\\
     \hline
     $\vb f_i, \, K_i, \, \psi_i$ & Dynamics, running cost, and terminal cost in the $i$-th mode  & \S\ref{sec:pdmps}\\
     \hline
     $\lambda_{ij} \geq 0$ & Transition rate from mode $i$ to mode $j$ & \S\ref{sec:pdmps}\\
     \hline
     $\Lambda$ & Matrix of transition rates & \S\ref{sec:pdmps}\\
     \hline
     $\x \in \Omega$ & A generic continuous state (also used as a generic initial condition) & \S\ref{sec:pdmps}\\ 
     \hline
     $\y_i(s) \in \Omega$ & Current continuous state at the time $s$ & \S\ref{sec:pdmps}\\
     \hline
     $\mu(s) \in \M$ & Current mode at the time $s$ & \S\ref{sec:pdmps}\\
     \hline
     $\beta > 0$ & Time-discounting rate in infinite-horizon problems & \S\ref{sec:pdmps} \\
     \hline
     $\gamma_i > 0$ & Mode-dependent termination rate in randomly-terminated problems & \S\ref{sec:pdmps}\\
     \hline
     $u_i$ & Value function for the $i$-th mode if all mode transitions are observed & \S\ref{sec:pdmps}\\
     \hline
     $\bv q \in \mathbb{Q}^{M-1}$  & Initial probability distribution over $M$ modes & \S\ref{sec:belief-prog}\\    
     \hline
     $\bv b(s) \in \mathbb{Q}^{M-1}$  & Current probability distribution over $M$ modes at time $s$ & \S\ref{sec:belief-prog}\\ 
     \hline
     $ \overline{K}_{\vb q}, \,  \overline{\psi}_{\vb q}$ & $\vb q$-dependent expected running cost and expected terminal cost & \S\ref{sec:reduced-b-prog}\\[2pt]
     \hline
     $G_{\vb q}(\vx, t, \nabla v)$ & $\vb q$-dependent Hamiltonian for finite-horizon problems & \S\ref{sec:reduced-b-prog}\\
     \hline
     $v_{\vb q}$ & Value function when mode transitions are not observed, starting from the mode distribution $\vb q$ & \S\ref{sec:reduced-b-prog}\\
     \hline
     $m \in \M$ & The last observed mode & \S\ref{sec:reduced-b-prog}\\
     \hline
     $v_{m}, \, G_m$ & Value function and the Hamiltonian when mode transitions are not observed, starting deterministically from mode $m \in \M$ & \S\ref{sec:reduced-b-prog}\\
     \hline
     $T_l$ & Time of the $l$-th mode observation (if known) & \S\ref{sec:finite}\\
     \hline
     $\tilde{T}_l = T_{l+1} - T_l$ & Time interval between two subsequent mode observations & \S\ref{sec:finite}\\
     \hline
     $v_m^l, \,  \overline{K}_m^l, \,  \Theta_m^{l+1}$ & Value function, expected running cost, and expected terminal cost on that time interval 
     if $m$ is the last observed mode at the time $T_l$. (Note: superscript $l$ is omitted if dealing with periodic observations in an infinite-horizon setting, as in \S \ref{sec:infinite}.) & \S\ref{sec:finite}\\
     \hline
     $F_{m} (\vx, t, \nabla v)$ & Initial-mode-dependent Hamiltonian for the infinite-horizon problem with periodic observations & \S\ref{sec:infinite} \\
     \hline
     $\TG$ & Policy-dependent time-to-target in indefinite-horizon problems & \S\ref{sec:indef}\\
     \hline
     $\bv q_s \in \mathbb{Q}^{M-1}$  & Stationary mode  distribution corresponding to a constant $\Lambda$& \S\ref{sec:indef}\\    
     \hline
     $\TO$ & Time of the next on-demand mode observation & \S\ref{sec:indef}\\
     \hline
     $C$ & Cost of possible on-demand mode observations & \S\ref{sec:indef-paid-for}\\
     \hline
     $\TR$ & Time of random premature termination & \S\ref{sec:rt}\\
     \hline
     $\phi_i$ & Cost of random premature termination when in mode $i$ & \S\ref{sec:rt}\\
     \hline
    $\overline{\phi}_m$ & Expected cost of random premature termination if starting from mode $m \in \M$ & \S\ref{sec:rt}\\ 
     \hline
     $R_m(\vx, t, v, \nabla v)$ & Operator combining Hamiltonian and expected cost of random premature termination for the randomly-terminated problem & \S\ref{sec:rt} \\
     \hline
    \end{tabular}
    \caption{\textbf{Summary of notation and variable names.}
    }
    \label{tab:params}
}
\end{table}


\subsection{Belief-state Programming}\label{sec:belief-prog}
Many optimal control problems with partial or incomplete information rely on belief-state programming to compute optimal policies \cite{astrom1965optimal,kaelbling1998pomdps}.
The belief encodes the available information about the current state, usually as a probability distribution over the state space.
With an appropriate method for updating the belief based on observations, it is then possible to convert a partially observed problem into a fully observed problem over the belief space.
Unfortunately, exactly computing optimal policies as a function of the belief is often computationally intractable due to the curse of dimensionality.
Given an optimal control problem over a finite state space with $N$ elements, the belief takes the form of a vector in the standard $(N-1)$-dimensional probability simplex $\mathbb{Q}^{N-1}$, and the cost of performing belief-state programming thus scales exponentially with the size of the state space.
For problems with continuous state space, belief-state programming instead involves optimizing over an infinite-dimensional space, which is generally only tractable if the belief can be parameterized using a finite number of parameters, e.g., as in the classic Kalman filter \cite{kalman1960filtering}.

There is existing work developing schemes for optimal control of switching systems with noisy state observations.
One setting common in {\em general} PDMPs is to assume that the planner receives noisy observations of the state 
immediately after each random jump \cite{brandejsky2013stopping, bauerle2018popdmps}.
As a result, the random jump times are always perfectly observed, and to be numerically tractable the number of possible post-jump states must be finite.
A different formulation is sometimes used in stochastic switching systems,  
where it is assumed that the system can be in one of a finite number of modes, each with its own version of stochastic dynamics.
The planner receives continuous observations of the current mode corrupted by Gaussian noise, and the belief is computed using the classic Wonham filter, see, e.g., \cite{yu2014asset}, \cite{tran2014lotka}.
However, since the control problem must still be solved over the entire belief space, this approach is only feasible when the number of modes is small.

We now \update{develop} a framework for belief-state programming in PDMPs of the type described in Section \ref{sec:pdmps} (taking the finite horizon case as a representative example).
We assume that the continuous portion of the state $\vy(s)$ is always known, but the mode switches remain unobserved.
We present a simplified case where the planner begins with an initial mode probability distribution $\vb q \in \mathbb{Q}^{M-1}$
and does not receive any future mode observations.\footnote{The same approach can be used when the planner receives (exact or noisy) infrequent observations of the current mode. A version with observation times known in advance 
is presented in the Appendix, Section \ref{app:noisy-obs}.}
\update{The following fundamental assumptions will be made throughout, to ensure} that the current mode cannot be deduced from the observed dynamics and increases in cumulative cost within each mode.
\begin{assumption}
\label{a:no_mode_info_from_traj}
To ensure the planner is not exactly aware of the current mode, we assume
\begin{enumerate}
\item The dynamics are mode-independent, so $\bv f_i(\vx, t, \veca) = \bv f(\vx, t, \veca)$.
\item \update{No differences in accumulating cost are observed until the process terminates.}
\end{enumerate}
\end{assumption}
Many of the examples discussed in Section \ref{sec:oopdmps-introduction} satisfy these assumptions\update{; e.g.,
an evader's range of attainable velocities is not directly influenced by the current mode/surveillance pattern, and
a rover may enter a ``breakdown-prone'' mode without changing its immediate dynamics.
For an evader primarily concerned with surveillance, the accumulating cost will be mode-dependent, but will not be observed directly: 
it is hard to asses your increasing surveillance exposure if you do not know how surveillance patterns have been changing up till now.  
For a rover, both its normal and breakdown-prone modes might have the same fully observed accumulating cost (e.g., the time or energy used so far along the way to the target), but these modes will have different rates of switching to a new (slower, partially broken) operating regime.  We will provide a detailed treatment of the latter case in \S \ref{sec:rt} and \S\ref{sec:rover}.
 }

We assume that the planner only has access to initial \update{mode} information in the form of an initial belief $\vb b(t) = \vb q$.
As previously noted, the current mode of the PDMP problem can be thought of as the state of a 
\update{continuous-in-time Markov chain (CTMC)}.
Specifically, given transition rates $\lambda_{ij}(\vx,t) \geq 0$, the CTMC has rate matrix $\Lambda(\vx,t)$, where
\begin{align*}
   [\Lambda(\vx,t)]_{ij} =
   \begin{cases}
         \lambda_{ij}(\vx,t), & i \neq j,\\
    		 -\sum_{k \neq i} \lambda_{ik}(\vx, t), & i = j.
	\end{cases}
\end{align*}
Since the belief is simply a distribution over the possible modes, along a trajectory $\vy(s)$ it evolves according to the deterministic dynamics
\begin{align}
\vb b'(s) &= \vb b(s) \, \Lambda\lp\vy(s),s\rp, 
\qquad \update{s \geq t;}
\label{eq:belief-statedep-ode}\\
\vb b(t) &= \vb q.
\end{align}
Due to the dependence of $\Lambda(\vx,t)$ on the continuous state, the belief will be policy-dependent.

We can now define the \emph{expected} running and terminal costs given any current belief $\vb q$ as
\begin{align}
\overline{K}(\vx, \vb q, t, \veca) = \sum_{i=1}^M q_i K_i(\vx, t, \veca), \quad \text{and} \quad \overline{\psi}(\vx, \vb q) = \sum_{i=1}^M q_i \psi_i(\vx).
\end{align}
The value function $w(\vx, \vb q, t)$ can then be defined to encode the minimal expected cost if we start at a time $t \in [0,T]$
with the continuous state $\vy(t) =\vx$ and the belief $\vb b(t) = \vb q.$  I.e.,
\begin{align}
	w(\vx, \vb q, t) = \inf_{\veca(\cdot)} \left\{ \int_t^T \overline{K}\lp\vy(s), \vb b(s), s, \veca\lp\vy(s), \vb b(s), s\rp\rp ds + \overline{\psi}\lp\vy(T), \vb b(T)\rp \right\}.
	\label{eq:belief_value_func}
\end{align}
This $w$ could be found by solving a PDE over $\Omega \times \mathbb{Q}^{M-1} \times [0,T].$ Specifically, 
\begin{align}
         - \frac{\partial w}{\partial t} &= 
         \min_{\veca \in A}\left\{\overline{K}(\vx, \vb q, t, \veca) + \mathbf{f}(\vx, t, \veca) \cdot \nabla_{\vx} w 
         + Q\lp \vx, \vb q, t \rp \cdot \nabla_{\vb q} w \right\},\label{eq:belief_pde}\\
        w(\vx, \vb q, T) &= \overline{\psi}(\vx, \vb q);
\end{align}
where $Q \lp \vx, \vb q, t \rp = \vb q  \Lambda \lp \x,t \rp$ in accordance with \eqref{eq:belief-statedep-ode}.

This problem could then be solved numerically over $\Omega \times \mathbb{Q}^{M-1} \times [0,T]$, with the computational cost growing exponentially with the number of modes.  In what follows, we make additional assumptions that allow us to avoid this curse of dimensionality by instead solving PDEs over $\Omega \times [0,T].$

\begin{rem}
\label{rem:pomdp_methods}
We note that the prohibitive computational cost of belief-space dynamic programming is well-known even in a fully discrete setting.
This is why approximate dynamic programming and heuristic-based algorithms are frequently used in Partially Observed Markov Decision Processes (POMDPs)
\cite{littman1995learning, kurniawati2022partially, kaelbling1998pomdps}.  A more efficient algorithm for a special class of POMDPs with {\em intermittent} state observations was recently introduced in \cite{basich2022planning} by reconstructing the belief from the last (exact) observation and a sequence of actions taken since then.  
Unfortunately, several features make the same ideas inapplicable in our continuous setting.  In \cite{basich2022planning}, no information about the state is obtained 
in between observations and each state transition might yield a new observation with a known probability.  
In our setting, it is only the {\em discrete mode} $\mu(s)$ that is infrequently observed (the continuous component of the state $\y(s)$ is always fully known), and 
mode-observations are not obtained at random. 
More importantly, the computational cost of the method developed in \cite{basich2022planning} grows exponentially with the number of available actions, 
while in our applications the set of control values $A$ is typically infinite.
\update{The same issue also strongly impacts the efficiency of sampling-based methods \cite{brown2012survey, silver2010planning}.
}
\end{rem}

\section{OOPDMP models}
\label{sec:models}
\subsection{Reduced Belief-state Programming: finite horizon, no observations}\label{sec:reduced-b-prog}
Turning to our approach, we now restrict the class of considered problems in order to make them computationally tractable.
We will consider only those problems where the evolving belief is not influenced by the controller's choices, and is 
thus 
only dependent on time and our prior distribution over $\M.$ 
To guarantee this, \update{from here on we will} make the following additional assumption.
\begin{assumption}
\label{a:const_switching_rates}
The transition rates between modes are constants: 
$\lambda_{ij}(\vx,t) = \lambda_{ij}$.
\end{assumption}
If the planner has initial belief $\vb q$ at time $t = 0$ 
and receives no further mode observations, the belief at some future time $s \geq 0$ can be computed using the matrix exponential\footnote{
For the purposes of finite-horizon problems, we can similarly handle $\lambda_{ij}(\vx,t) = \lambda_{ij}(t),$ though $\vb b'(s) = \vb b(s) \Lambda(s)$ would then
generally be solved numerically. 
The assumption of constant $\lambda_{ij}$ is made in view of the other (stationary) problem classes.}
 $\vb b(s) = \vb q \exp(s \Lambda)$.
In some applications, this is quite a reasonable assumption to make, as in surveillance-evasion scenarios in which the switching surveillance patterns are not impacted by the evader's chosen trajectory.
In others it is less realistic; e.g., a planetary rover's rate of switching to a breakdown-prone mode would in reality be affected by the terrain through which that rover is traveling.

Let us again focus on a finite horizon control problem with no observations.
This could, for example, represent an
evader remaining in the domain until some time $T$, minimizing their exposure to surveillance according to mode-dependent surveillance intensity functions $K_i(\vx)$, with access to only an initial distribution $\vb q$ over the possible surveillance modes/patterns. (An example of this type is included in Section \ref{sec:fin-results}.)
For a \emph{fixed} initial belief $\vb q$, we can now compute the expected running and terminal costs with respect to the current belief as $\overline{K}_{\vb q}(\vx, t, \veca) = \overline{K}\lp\vx, \vb b(t), t, \veca\rp$ and $\overline{\psi}_{\vb q}(\vx, t) = \overline{\psi}\lp\vx, \vb b(t)\rp$ where the subscript encodes the dependence of $\vb b(t)$ on that initial belief.
Using open-loop control policies $\ba: \R \mapsto A,$
we can similarly define a value function for a fixed $\vb q$ as 
\begin{align}
	v_{\vb q}(\vx, t) = \inf_{\veca(\cdot)} \left\{ \int_t^T \overline{K}_{\vb q} \left( \vy(s), s, \ba(s) \right) ds + \overline{\psi}_{\vb q}\lp\vy(T),T\rp \right\},
\end{align}
which can be computed by numerically solving 
\begin{align}
	-\frac{\partial v_{\vb q}}{\partial t}(\vx,t) &= G_{\vb q}(\vx, t, \nabla v_{\vb q}),  &\vx& \in \Omega, \, t \in [0, T]; \label{eq:fin-none-pde}\\
	v_{\vb q}(\vx,T) &=\overline{\psi}_{\vb q}(\vx, T), &\vx& \in \Omega \label{eq:fin-none-tc},
\end{align}
where $G_{\vb q}(\vx, t, \nabla v_{\vb q}) = \min_{\veca\in A}\left\{\overline{K}_{\vb q}(\vx, t, \veca) + \vb f(\vx, t,\veca) \cdot\nabla v_{\vb q}(\vx,t) \right\}.$
The argmin in $G_{\vb q}$ can be then used to recover the optimal policy in feedback form $\ba_* = \ba_*(\x,t).$
This is equivalent to solving a standard, fully observed finite horizon control problem over the domain $\Omega \times [0,T]$, which is significantly cheaper than solving equation \eqref{eq:belief_pde}.
It should be noted, however, that computing $v_{\vb q}$ recovers optimal feedback policies for \emph{any} starting state $\vx$ but \emph{only one} starting belief $\vb q$.
(Unlike in \eqref{eq:belief_value_func}, here the feedback policy no longer includes the current belief $\vb b$ among its arguments,
but the initial belief $\vb q$ still very much influences the cost of every such policy.)

So far, we have focused solely on the case that the planner has initial information about the mode and receives no further mode observations.
When we examine observation schemes, we will be particularly interested in the special case that the planner perfectly knows the initial mode $m = \mu(0)$ (encoded as a belief by the $m$th standard basis vector $\vb e_m$).
In this case, we will slightly abuse notation and write $\overline{K}_m$ and $\overline{\psi}_m$ to denote the expected running and terminal costs associated with starting deterministically in mode $m$ (so $\vb b(0) = \vb e_m$).
If $v_m(\vx, t)$ is similarly defined to be the value function associated with $\vb b(0) = \vb e_m$, then it will satisfy equation \eqref{eq:fin-none-pde} with an initial-mode-dependent Hamiltonian $G_m(\vx, t, \nabla v_m) = G_{\vb e_m}(\vx, t, \nabla v_{\vb q})$.
\update{For readers' convenience, we summarize these and other notational conventions in Table \ref{tab:params}.}

In the rest of Section \ref{sec:models}, we analyze various scenarios in which the planner has access to infrequent mode observations after the initial time $t = 0$. 
Our goal is not to enumerate all combinations of planning horizons and observation schemes.
Instead, we focus on several representative examples and highlight the differences in  
theoretical and computational aspects of the resulting control problems. 
In some cases (Section \ref{sec:finite}), observations can be incorporated with minimal changes to the structure of the original problem. 
In \ref{sec:infinite}, we show how the introduction of regular observations converts a stationary problem over an infinite time horizon into a time-dependent problem over a finite interval $[0,T]$ with nonlocal coupling.
In \ref{sec:indef}, we examine how the introduction of time-dependence into previously stationary problems can be exploited to design efficient numerical methods.
We close with a framework in \ref{sec:rt} motivated by applications in which unobserved mode transitions affect the likelihood of a (premature) random termination. 

\subsection{Finite Horizon, Observations at Known Times}\label{sec:finite}
We continue our focus on finite horizon problems, but now assume that the planner knows the starting mode $\mu(0)$ exactly and receives exact observations of the current mode at $L$ predetermined times $0 < T_1, \hdots, T_L < T$.
These observations create a natural partition of the time horizon $[T_0=0,\, T_{L+1} = T]$ into $L+1$ intervals of the form $[T_l, T_{l+1}]$, for $l \in \{1, 2, \hdots, L\}$.
When restricted to just one interval $[T_l, T_{l+1}]$, the problem has the same structure as the case with only initial information outlined above, with the mode known at $T_l$ and a terminal condition at $T_{l+1}$ (now possibly based upon the expected result of the observation). 
If we let $\tilde{T}_l = T_{l+1} - T_l$ be the duration of each such interval, we can thus solve over $[0, \tilde{T}_l]$ between each pair of observations.
Specifically, suppose at time $T_l$ the planner learns that $\mu(T_l) = m$, and let $\vb b(s)$ be the resulting evolving belief with $\vb b(0) = \vb e_m$ for $s \in [0, \tilde{T}_l]$.
We denote the expected running cost on the $l$-th interval as $\overline{K}^l_m(\vx, s, \veca) = \sum_{n=1}^M K_m\lp\vx, \vb b\lp s\rp,s+ T_l, \veca\rp$ where we compute each $K_m$ at $t = s + T_l$, the total time passed.
Similarly, let $\vb f^l(\vx, s, \veca) = \vb f^l(\vx, s+T_l, \veca)$ denote the shifted dynamics on the $l$-th interval.

We now let $v_m^l(\vx, s)$ denote the value function over the $l$-th interval with respect to the initial belief $\vb q = \vb e_m$. 
Since the terminal cost for the final interval is given by $\overline{\psi}_m(\vx, s)$, we can write 
\begin{equation}
v_m^L(\vx, s) = \inf_{\veca(\cdot)} \left\{\int_s^{\sscap{\tilde{T}}{L}} \overline{K}_m^L \lp\vy(r), r, \veca(r) \rp \,dr + \overline{\psi}_m\lp\vy(\sscap{\tilde{T}}{L}), \sscap{\tilde{T}}{L}\rp \right\}.
\end{equation}
where $\vy(r)$ is computed with respect to the dynamics $\vy'(r) = \vb f^L\left( \vx, r, \veca(r) \right)$ for $r \in [0, \sscap{\tilde{T}}{L}].$
For earlier intervals, the terminal cost is given by $\sum_{n=1}^M b_n(\tilde{T}_l) v^{l+1}_n(\vx,0)$, which encodes the expected optimal cost-to-go after the next observation is received at time $T_{l+1}$. 
Note that the previous observation is assumed to be $\mu(T_l) = m$; so, $\vb b(\tilde{T}_l) = \vb e_m \exp(\tilde{T}_l\Lambda).$ 
Linear combinations of this form will be useful throughout the problems we consider; so, we define
\begin{equation}
\Theta_m^l(\vx, s) = \sum_{n=1}^M b_n(s) v^l_n(\vx, 0) \label{eq:exp-obs-value}
\end{equation}
to represent this more compactly. (We will also use $\Theta_m(\vx, t)$ defined in the same way when there is no dependence on the observation number $l.$)
For $l = 0, \hdots, L-1$, we can now write
\begin{equation}
v_m^l(\vx, s) = \inf_{\veca(\cdot)} \left\{\int_s^{\tilde{T}_{l+1}} \overline{K}_m^l \lp\vy(r), r, \veca(r) \rp \,dr + \Theta_m^{l+1}(\vy(\tilde{T}_l), \tilde{T}_l)\right\}.
\end{equation}
Defined in this way, the value functions solve the following system of PDEs (with coupling through the terminal condition)
\begin{align}
	-\frac{\partial v_m^l}{\partial s} &= G^l_m\lp\vx, s, \nabla v_m^l\rp,  &\vx& \in \Omega, s \in [0, \tilde{T}_{l}] \\
	v_m^l(\vx,\tilde{T}_l) &= \begin{cases}
		\Theta^{l+1}_m(\vx, \tilde{T}_l), & \text{ if } l < L,\\[2pt]
		\overline{\psi}_m(\vx, \tilde{T}_l) & \text{ if } l = L.
	\end{cases} &\vx& \in \Omega,
\end{align}
with $G^l_m(\vx, s, \nabla v_m^l) = \min_{\veca\in A}\left\{\overline{K}^l_m(\vx,s,\veca) + \vb f^l(\vx,s,\veca) \cdot\nabla v_m^l(\vx,s) \right\}$.
This system can be solved numerically backwards in time using standard discretization schemes.
The cost of doing so is equivalent to that of solving $M$ copies of equation \eqref{eq:fin-none-pde}, since each PDE is solved over $[0,\tilde{T}_l]$ and $\sum_{l=0}^L \tilde{T}_l= T$.

\begin{rem}
Once the $v_m^l$ are computed for all $m \in \M$ and $l = 1, 2, \hdots, L$, we can compute the value function $v^0_{\vb q}$ corresponding to an arbitrary initial distribution $\vb q$ by solving
\begin{align}
        -\frac{\partial v^0_{\vb q}}{\partial s}(\vx,s) &= G_{\vb q}^0(\vx, s, \nabla v_{\vb q}^0) &\vx& \in \Omega, s \in [0, \update{\tilde{T}_0}) \\
        v^0_{\vb q}(\vx, \update{\tilde{T}_0}) &= \Theta_{\vb q}^1(\vx, \update{\tilde{T}_0}) &\vx& \in \Omega.
\end{align}
\end{rem}

\subsection{Infinite Horizon, Periodic Observations}\label{sec:infinite}
We now turn our attention to an example of infinite horizon problems with periodic mode observations.
Here and for the rest of Section \ref{sec:models}, we will make the following standard assumption
\begin{assumption}
\label{a:autonomous_cost_and_dynamics}
The mode-dependent dynamics and running costs are autonomous: $\vb f(\vx, t, \veca) = \vb f(\vx,\ba)$ and $K_i(\vx, t, \veca) = K_i(\vx, \veca).$ 
\end{assumption}
Our goal is to minimize the expected discounted cost over all time, as in \eqref{eq:inf-full-cost}, but with the current mode $\mu(s)$ observed only occasionally/periodically,
at times $s=T, 2T, 3T, ...,$  
\update{where the time between observations $T$ is known in advance.}
As in the previous subsection, we are interested in controlling the system in between these observations.
Since the structure of this problem is identical on each interval $[lT, (l+1)T],$  it is sufficient to find optimal controls over the first interval $[0,T].$
Indeed,  at $t=T$ the planner exactly observes the current mode $\mu(T)$ and faces an infinite horizon optimization problem identical to the one posed at $t = 0.$

Even though $K_i$ and $\vb f$ are autonomous, the value functions here depend on $t \in [0,T],$ interpreted as the time since the last observation.
This is a significant difference from the fully-observed case \eqref{eq:inf-full-cost}.
If $v_m(\vx,t)$ is defined as the optimal expected cost starting from $\y(t) = \x$ under the assumption that $\mu(0) = m,$ 
that value function must satisfy
\begin{equation}
v_m(\vx,t) = \int_t^T e^{-\beta s} \overline{K}_m\lp\vy(s), s, \veca(s)\rp \,ds + e^{-\beta T} \Theta_m\left(\vy(T), T\right).
\end{equation}
where $\Theta_m(\vx, t) = \sum_{n=1}^M b_n(t) v_n(\vx, 0)$ and $\vb b(t) = \vb e_m \exp(t \Lambda).$
Due to $\Theta_m,$  these value functions are coupled; moreover, this coupling is \emph{nonlocal in time}:
each $t$-slice of $v_m$ depends on the zeroth time slice of all $v_n.$
This unusual feature will also have consequences in our approximating these value functions numerically. 

For any $\tau < T-t,$ this value function satisfies the optimality  principle
\begin{align}
\label{eq:inf_hor_optim_principle}
v_m(\vx, t) 
&= \inf_{\veca(\cdot) \in \mathcal{A}} \left \{ \int_t^{t+\tau} e^{-\beta s} \overline{K}_m\left( \vy(s), s, \veca_m(s) \right) \, ds + v_m \left(\vy(t+\tau), t+\tau \right) \right \},
\end{align}
which is analogous to that of a finite horizon control problem -- though with the running cost $e^{-\beta t} \overline{K}_m(\vx, t, \veca)$ and a nonstandard terminal condition.
The usual formal argument based on Taylor-expanding \eqref{eq:inf_hor_optim_principle} yields
\begin{align}
- \frac{\partial v_m}{\partial t} &= F_m(\vx, t, \nabla v_m),  &\vx& \in \Omega, t \in [0, T), \label{eq:inf-per-pde}\\
v_m(\vx, T) &= e^{-\beta T} \Theta_m(\vx, T), &\vx& \in \Omega,\label{eq:inf-per-tc}
\end{align}
where now $F_m(\vx, t, \nabla v_m) = \min_{\veca \in A} \left\{ e^{-\beta t} \overline{K}_m(\vx, t, \veca) + \vb f(\vx, \veca) \cdot \nabla v_m(\vx, t) \right\}$. 

\begin{rem}
While equations \eqref{eq:inf-per-pde} - \eqref{eq:inf-per-tc} do not look like traditional infinite horizon HJB equations, they are equivalent when the problem has no mode-dependence, i.e., $K_m(\vx, \veca) = K(\vx, \veca)$. 
It is then sufficient to compute a single value function which satisfies
\begin{align}
- \frac{\partial v}{\partial t} &=\min_{\veca \in A} \left\{ e^{-\beta t} K(\vx,\veca) + \bv f(\vx, \veca) \cdot \nabla v(\vx, t) \right\},  \qquad t \in [0, T) \label{eq:inf-per-pde-no-modes}\\
v(\vx, T) &= e^{-\beta T} v(\vx, 0).
\end{align}
Setting $v(\vx,t) = e^{-\beta t} u(\vx)$ 
and substituting this into equation \eqref{eq:inf-per-pde-no-modes} recovers the standard HJB PDE for discounted infinite horizon problems:
$-\beta u + \min_{\veca \in A} \left\{ K(\vx,\veca) + \bv f(\vx, \veca) \cdot \nabla u(\vx) \right\} = 0.$
\end{rem}

\begin{rem}
Just as in the previous subsection, once we have solved for $v_1, \hdots, \sscap{v}{M},$ it is relatively inexpensive to compute $v_{\vb q}$ for an arbitrary initial distribution $\vb q$.
In this case, we have
\begin{align}
- \frac{\partial v_{\vb q}}{\partial t}  &= F_{\vb q}(\vx, t, \nabla v_{\vb q}),  \qquad t \in [0, T), \label{eq:inf-per-pde-q}\\
v_{\vb q}(\vx, T) &= e^{-\beta T} \Theta_{\vb q}(\vx, T), \label{inf-per-tc-q}
\end{align}
where $\Theta_{\vb q}(\vx, T)$ can be computed a priori for all $\vx \in \Omega$ once all $v_m$ are known.

A similar approach (for planning starting from an arbitrary $\bq \in \Q$ and up to the first observation) can be also used with
all models considered in the following subsections.  To limit the length of the paper, we do not include the actual PDEs and quasi-variational inequalities 
for $\bq \neq \vb e_m$ in each horizon/observation model.
However, some of our examples in Section \ref{sec:numexp} illustrate this more general case.
\end{rem}

\subsection{Indefinite Horizon Control Problems}\label{sec:indef}
Turning to the setting where the process terminates upon 
entering the interior of
 target $\Gamma,$
we make the same Assumption \ref{a:autonomous_cost_and_dynamics} and 
start with 
the version of this problem in which the planner only has access to an initial mode distribution $\vb q$  
without
further mode observations.
Just as in the fully-observed case, the planning horizon $\TG = \inf\{s \geq t \mid \y(t) = \x \in \Omega \backslash \Gamma, \, \y(s) \in \Gamma \backslash \partial \Gamma\}$ 
depends on the initial state and 
the chosen policy $\ba:\R \mapsto A.$  
But due to Assumption \ref{a:no_mode_info_from_traj} and the lack of observed mode-transitions, this  $\TG$ will have the same deterministic value regardless of the starting mode.
We define the value function for a fixed $\vb q$ as
\begin{equation}
\label{eq:indef_hor_val_func}
v_{\vb q}(\vx, t) = 
\inf_{\veca(\cdot)} \left\{\int_{t}^{\TG\left(\vx, \veca(\cdot)\right)} 
\overline{K}_{\vb q}\left(\vy(s), s, \veca(s)\right) \, ds + \overline{\psi}_{\vb q}\left(\vy\left(
\TG\left(\vx, \veca(\cdot)\right)
\right)\right) \right\}.
\end{equation}
As in Section \ref{sec:infinite}, the value function is now \emph{time-dependent}, due to the dependence of the optimal policy on the evolving belief (through $\overline{K}_{\vb q}$ and $\overline{\psi}_{\vb q}$).
A classical argument yields the time-dependent governing HJB PDE
\begin{align}
         -\frac{\partial v_{\vb q}}{\partial t}(\vx,t) &= G_{\vb q}(\vx, t, \nabla v_{\vb q}), &\vx& \in \Omega\setminus\Gamma, \, t \geq 0 \label{eq:indef-none-pde}\\
        v_{\vb q}(\vx,t) &= \overline{\psi}_{\vb q}(\vx, t), &\vx& \in \partial \Gamma, \, t \geq 0 \label{eq:indef-none-bc}.
\end{align}
where we recall that $G_{\vb q}(\vx, t, \nabla v_{\vb q}) = \min_{\veca\in A}\left\{\overline{K}_{\vb q}(\vx,t,\veca) + f(\vx,\veca) \cdot \nabla v_{\vb q} \right\}$.
We note that equation \eqref{eq:indef-none-pde} is identical to the governing PDE in \eqref{eq:fin-none-pde} for a finite horizon control problem with no observations,
but specifies only a boundary condition in \eqref{eq:indef-none-bc} instead of a terminal condition as in \eqref{eq:fin-none-tc}.
Solving this PDE for $t \in [0, +\infty)$ is not practical, and we need to obtain an equivalent formulation on a finite time interval.

\begin{prop}\label{prop:bounded-TG}
Suppose Assumptions \ref{assump:regularity}-\ref{a:autonomous_cost_and_dynamics} hold and in addition
the controlled dynamics is ``geometric''; i.e., $\vb f(\vx, \veca) = f(\vx, \veca) \veca,$
where $\veca \in A = \mathbb{S}^1$ is a unit vector specifying the chosen direction of motion while $f$ is the corresponding speed of motion.
Suppose also that the explicit bounds
\begin{itemize}
\item $0 < K_{\min} \leq K_{\update{i}}(\vx, \veca) \leq K_{\max}$ 
\item $0 < f_{\min} \leq f(\vx, \veca) \leq f_{\max}$
\item $0 \leq \psi_m(\vx) \leq \psi_{\max}$
\end{itemize}
hold for all $\vx \in \Omega$, $\veca \in \mathbb{S}^1$, $m\in\mathcal{M}.$\\
Let $z(\vx)$ be the minimal time needed to reach the target $\Gamma$ starting from $\x \in \Omega \backslash \Gamma.$
Suppose that $\ba_*(\cdot)$ is any optimal control and $\y_*(\cdot)$ is the corresponding optimal trajectory
for a starting point $\x=\y_*(0)$ and some initial mode distribution $\bq \in \Q.$ 
Then
\begin{equation}
\TG\left(\x, \ba_*(\cdot)\right)
 \; \leq \; \frac{z(\vx)K_{\max}}{K_{\min}}  + \frac{\psi_{\max}}{K_{\min}}. \label{eq:TG-bound}
\end{equation}
\end{prop}
\begin{proof}
Under these assumptions,  the minimum time to target $z(\x)$ is bounded and locally Lipschitz-continuous. 
Moreover, it can be recovered as a viscosity solution of
the stationary HJB PDE
\begin{align}
\min_{\veca \in \mathbb{S}^1} \left\{f(\vx, \veca) \veca \cdot \nabla z(\vx) \right \} \, + \, 1 &=0, &\vx&\in \Omega \setminus \Gamma; \label{eq:min-time-pde}\\
z(\vx) &= 0, & \vx &\in \partial \Gamma.
\end{align}
Based on Assumption \ref{assump:regularity}, there exists a time-optimal open loop control $\bahat: \R \mapsto A$ and the corresponding trajectory $\yhat(s)$
such that $\TG\left(\x, \bahat(\cdot)\right) = z(\x).$  We note that $\y_*(s)$ 
cannot take too long on the way to
$\Gamma$ since otherwise the expected cost incurred along $\yhat(s)$ might be lower. 
Indeed,
\begin{align}
v_{\bq}(\vx, 0) \; \leq \;
J_{\vb q}(\x,0, \bahat(\cdot)) &= 
\int_0^{z(\x)} \overline{K}_{\vb q}\left( \yhat(s), s, \bahat(s) \right)  \, ds 
\, + \, 
\overline{\psi}_{\vb q} \left(\yhat \left(  z(\x)  \right)\right)
\; \leq \; z(\vx) K_{\max} + \psi_{\max}.
\label{eq:time_subopt_bound}
\end{align}
At the same time, by the optimality of $\ba_*(\cdot),$
\begin{align}
\nonumber
v_{\bq}(\vx, 0) \; = \;
J_{\bq}\left(\x, 0, \ba_*(\cdot)\right) &=
\int_0^{\TG\left(\x, \ba_*(\cdot)\right)} 
\overline{K}_{\bq}\left( \y_*(s), s, \ba_*(s) \right)  \, ds 
\, + \, 
\overline{\psi}_{\bq} \left(\y_* \left(  \TG\left(\x, \ba_*(\cdot) \right)  \right) \right)\\
& \geq  \; K_{\min} \TG\left(\x, \ba_*(\cdot)\right).
\end{align}
Combining these and dividing both sides by $K_{\min},$ we obtain \eqref{eq:TG-bound}.
\end{proof}

Since $z(\x)$ is bounded, we can choose any 
$T > \frac{ K_{\max}}{K_{\min}} \max_{\vx \in \Omega} \left\{z(\vx)\right\}  + \frac{\psi_{\max}}{K_{\min}}$
as the planning horizon, solving \eqref{eq:indef-none-pde}-\eqref{eq:indef-none-bc} on a time interval $[0,T]$
with a terminal condition
\begin{equation}
v_{\vb q}(\vx, T) = \begin{cases}
	\overline{\psi}_{\vb q}(\vx, T), & \text{if } \vx \in \partial \Gamma; \\
	+\infty. & \text{if } \vx \notin \Gamma.
	\end{cases} \label{eq:indef-none-tc}
\end{equation}
Importantly, this terminal condition will not impact
the value function in the time slice $t=0$ and on all relevant characteristics $(\y_*(s), s)$ from that time slice to $\Gamma \times [0, +\infty).$
(All optimal trajectories starting from every $\x \in \Omega$ at $t=0$ will reach $\Gamma$ by the time $T$.)

\begin{rem}
\label{rem:stationary_q}
One case in which this machinery is not needed is when $\vb b(0) = \vb q_s$, where $\vb q_s$ is a stationary distribution of the CTMC associated with the mode switching process.
Since this initial belief implies $\vb b(t) = \bq_s$ for all $\update{t \geq 0},$
the value function $v_{\vb q_s}$ would be no longer time-dependent and could be recovered by solving to a stationary PDE
$\min_{\veca\in A}\left\{\overline{K}_{\vb q_s}(\vx, \veca) + f(\vx, \veca) \cdot \nabla v_{\vb q_s} \right\} = 0$ on $\Omega \backslash \Gamma$ with 
$v_{\vb q_s} = \overline{\psi}_{\vb q_s}$ on $\partial \Gamma.$
But for $\Lambda \neq 0$ and any other initial belief $\vb q,$ the time-dependent PDE \eqref{eq:indef-none-pde} is unavoidable.
\end{rem}

We now analyze two possible mode-observation schemes for indefinite horizon problems:
(a) an a priori limited number of on-demand mode-observations that are available ``for free'' and
(b) unlimited on-demand mode-observations that are available at a fixed cost per request.
In both cases, the value functions will satisfy quasi-variational inequalities, but their structure 
(and the associated computational challenges) will be quite different.
The following operator will be a useful tool in discussing both:
given a time $\TO$ (which we will generally interpret as an observation time) and an expected cost-to-go-after-observation $B$,
we define
\begin{align}
\mathcal{J}_m\left(\vx, t, \veca(\cdot), \TO, B\right) &= \int_t^{\min\left\{\TG, \TO\right\}} \overline{K}_m\left(\vy(s),s,\veca(s)\right) ds \nonumber\\
&+ \mathds{1}_{\left[\TO < \TG\right]}B \; + \; \mathds{1}_{\left[\TG \leq \TO\right]} \overline{\psi}_m\left(\vy\left(\TG\right), \TG\right). 
\label{eq:indef-cost-operator}
\end{align}

\subsubsection{Bounded Number of On-demand Observations}\label{sec:indef-bdd}

Suppose a planner has access to at most $L$ mode observations and can choose when to use them while navigating the domain.
The value functions $v_m^l(\x,t)$ describing the minimal expected cost to go will now depend on the time since the last mode observation $t,$
the last observed mode $m$, the current continuous state $\x$, and the number of observations already made $l \in \{0,...,L\}.$
This is the same notation already introduced in Section \ref{sec:finite}, where observation times were prescribed, and we can use 
the  same $\Theta_m^{l+1}(\x, t)$ defined in \eqref{eq:exp-obs-value}
to represent the expected remaining cost after an immediate mode observation.  
For $l < L,$ the challenge is to choose not only a policy $\ba(\cdot)$ but also the next observation time $\TO.$
Using the operator introduced in \eqref{eq:indef-cost-operator}, the value functions must then satisfy
\begin{align}
v^l_m(\vx, t) \; = \; \inf_{\substack{\ba(\cdot)\\ \TO \geq t}} 
&\left\{
\int_t^{\min\left\{\TG, \TO\right\}} \overline{K}_m
\left(
\vy(s),s,\veca(s)
\right) ds 
\right. \nonumber\\
&+ \left.
\mathds{1}_{
\left[
\TO < \TG
\right]
}\Theta_m^{l+1} \!\!
\lp 
\y 
\left(
\TO 
\right), \TO 
\rp  \; + \; \mathds{1}_{\left[\TG \leq \TO\right]} \overline{\psi}_m\left(\vy\left(\TG\right), \TG\right)
\right\}
\end{align}
Assuming the next observation is not taken for a small time $\tau,$  we can write the following optimality principle
\begin{align}
v_m^l(\x, t) = \inf_{\ba(\cdot)} &\left \{ \int_t^{t+\tau} \overline{K}_m\left(\vy(s), s, \ba(s) \right) \, ds \right. \nonumber \\
&\quad \left.+ \min \left\{ v^l_m(\vy(t+\tau), t+\tau), \, \Theta_m^{l+1} \left(\vy(t + \tau), t+\tau \right) \right\} \right \} + o(\tau),
\end{align}
where the inner minimization encodes the choice of whether or not to make an observation at the time $t+\tau$.
A Taylor-series expansion yields a governing system of $M \times L$ quasi-variational inequalities,
\begin{align}
0 &= \max\left\{
- \frac{\partial v_m^l}{\partial t} - G_m(\vx, t, \nabla v^l_m), \; v_m^l(\vx, t) - \Theta_m^{l+1}(\vx, t) \right\}, &\vx& \in \Omega\setminus\Gamma, \, t \geq 0; \label{eq:indef-bdd-pde}\\
v_m^l(\vx, t) &= \overline{\psi}_m(\vx,t), &\vx &\in \partial \Gamma, \, t \geq 0, \label{eq:indef-bdd-bc}
\end{align}
solved for $l = 0,...,(L-1).$
Each $v_m^L(\vx,t)$ (for the no observations remaining case) can be found from the HJB PDE \eqref{eq:indef-none-pde}-\eqref{eq:indef-none-bc} 
with $\vb q = \vb e_m$.

If the conditions in Proposition \ref{prop:bounded-TG} are satisfied, then this problem is similarly easy to restrict to a finite time interval.
Taking too long to wait for a mode observation or to reach $\Gamma$ would again be dominated by the expected cost incurred while following a time-optimal trajectory.

\subsubsection{On-demand Observations with Positive Cost}
\label{sec:indef-paid-for}
Suppose a planner has to pay the cost $C(\x)>0$ for a mode observation requested at that location.
(For a planetary rover, this could be a time-penalty or energy-penalty for running the diagnostics necessary to determine 
whether it has transitioned into a breakdown-prone mode.  For a prey animal, this could be energy losses due to taking 
some time away from foraging to check for possible predators nearby.)
If there is no enforced maximum $L,$ 
it is no longer necessary to keep track of the number of past observations.
Thus, the value functions $v_m(\vx, t)$ can be defined  through
\begin{align}
v_m(\vx, t) \; = \;
\inf_{\substack{\ba(\cdot)\\ \TO \geq t}}  
\left\{ \mathcal{J}_m
\left( \x, \, t, \, \ba(\cdot), \, \TO, \, \,
C\left(\y\left(\TO\right)\right) 
+ \Theta_m 
\left( \y\left(\TO\right), \TO\right)\right) 
\right\}.
\end{align}

A standard argument based on Taylor-expanding the optimality principle yields a system of $M$
quasi-variational inequalities 
\begin{align}
0 &= \max\left\{ 
-\frac{\partial v_m}{\partial t} - G_m\left(\vx, t, \nabla v_m\right), \; v_m(\vx, t) - C(\vx) - \Theta_m(\vx, t)\right\}, 
&\vx& \in \Omega\setminus\Gamma, \, t \geq 0; \label{eq:indef-paid-pde}\\
v_m(\vx, t) &= \overline{\psi}_m(\vx, t), &\vx& \in \partial \Gamma, \, t \geq 0. \label{eq:indef-paid-bc}
\end{align}
The key distinction between \eqref{eq:indef-bdd-pde} and \eqref{eq:indef-paid-pde} is that the latter system has a (nonlocal in time) coupling through $\Theta_m$.
In \eqref{eq:indef-bdd-pde}, each $v_m^l$ only depends on $v_1^{l+1}, \hdots, v_M^{l+1}$, the value functions when there is one fewer observation remaining.
As a result, the system \eqref{eq:indef-bdd-pde} can be solved in a single pass (from $l = L$ to $l=0$) while the numerical solution of \eqref{eq:indef-paid-pde} is unavoidably iterative.

\subsection{Indefinite Horizon with (Random) Premature Terminations}\label{sec:rt}
In many applications, a process that deterministically terminates upon reaching the target $\Gamma$ might also terminate earlier as a result of some random event experienced en route.
This ``premature'' random termination typically incurs a much higher cost $\phi_{\update{i}}(\vx).$ 
Such termination could be either a true completion of this controlled process (e.g., a death of a foraging animal due to predation) or an observed and permanent switch 
to an entirely different operating environment (e.g., a planetary rover suffering a breakdown and reducing the speed from there on \cite{gee2022breakdowns}
or an illegal \update{forest} logger switching to new path-planning goals after being apprehended by a ground patrol \cite{CarteeVlad_Poaching}).
In the latter cases, $\phi_{\update{i}}(\vx)$ might be used to encode the remaining costs to go in the new operating environment.  We will call all such models {\em randomly terminated} and
will assume that the  premature termination result from a nonhomogeneous Poisson process, whose rate $\gamma$ might generally depend on the current position in $\y(t) \in \Omega$ or 
the current (unobserved) mode $\mu(t)$.

We first focus on the former case (i.e., the termination rate $\gamma(\x))$ and assume no additional mode observations beyond $\mu(0) = m.$
Since a premature termination is always observed, we can again define the value functions using the operator \eqref{eq:indef-cost-operator} as
\begin{equation}
v_m(\vx, t) = \inf_{\veca(\cdot)} \left\{\E_{\TR}^{} \left[ \mathcal{J}_m\left(\x, t, \, \ba(\cdot), \, \TR, \, 
\overline{\phi}_m\left(  \vy(\TR), \, \TR\right) \right) \right]
\right\},
\end{equation}
where $\overline{\phi}_m(\vx,t) = \sum_{n=1}^M b_n(t) \phi_n(\vx)$ with $\vb b(0) = \vb e_m$. 
These value functions can be thus recovered by solving a system of PDEs 
\begin{align}
- \frac{\partial v_m}{\partial t} &= G_m(\vx, t, \nabla v_m)
+ \gamma \left(\overline{\phi}_m(\vx,t) - v_m(\vx,t)\right),  &\vx& \in \Omega\setminus\Gamma, \, t \geq 0; \label{eq:rt-none-pde}\\
v_m(\vx, t) &= \overline{\psi}_m(\vx, t),  &\vx& \in \partial \Gamma, \, t \geq 0. 
\end{align}
If all $\phi_{\update{i}}(\x)$ are sufficiently high, it is never optimal to wait for a premature termination instead of trying to reach $\Gamma,$
and the approach described in Section \ref{sec:indef} can be also used here to solve \eqref{eq:rt-none-pde} on a bounded time interval.

We generally interpret a premature termination as  an undesirable event; e.g., a serious breakdown that impacts the dynamics of a traveling rover.
But to model the accumulation of unobservable damage that increases the likelihood of such breakdowns, it is natural to consider mode-dependent termination rates $\gamma_i$ (which we now take to be state-\emph{independent}, see Remark \ref{rem:mode-and-state-gamma}).
This introduces an additional source of information about the current mode, namely, if the process has not yet terminated, then the planner is less likely to be in modes with high termination rates (relative to a prediction based solely on the initial distribution and $\Lambda$).
To account for this additional information, the belief must now be defined as
\begin{equation}
b_i(t) = \Prob(\mu(t) = i \mid \vb b(0) = \vb q, \lnot \Xi(t)),
\end{equation}
where $\Xi(t)$ is the event that the process terminates prematurely in the interval $[0, t]$.

We can find an expression for $\vb b(t)$ by solving an auxiliary problem: a pure death process with inhomogenous death rates and mixing.
Suppose we have $M$ sub-populations that mix with rates $\lambda_{ij}$ and each has a death rate $\gamma_i.$
If $r_i(t)$ is the proportion of the initial population that is in sub-population $i$ at the time $t$, we have
\begin{align}
\bv r'(t) &= \vb r(t)(\Lambda - \text{diag}(\bgamma)) \label{eq:r-ode}
\end{align}
where $\bv r(t)$ is a row vector with $i$-th entry $r_i(t)$ and $\bgamma$ is a vector with $i$-th entry $\gamma_i$.
Specifying an initial condition $\vb r(0) = \vb q$, we can write
\begin{equation}
b_i(t) = \frac{
\Prob
\left( \left( \mu(t) = i \right) \cap \lnot \Xi(t) \mid \vb q \right)
}{\Prob(\lnot \Xi(t)\mid \vb q)} = \frac{r_i(t)}{\sum_{n=1}^M r_n(t)}. \label{eq:b-rt}
\end{equation}
Thus, $\vb b(t)$ can still be computed as an explicit function of time and the framework presented in Section \ref{sec:reduced-b-prog} is still applicable.\footnote{
It is also possible to directly derive a system of governing ODEs for $\vb b(t)$ using Bayes Theorem; see Appendix \ref{app:b-ode-bayes}.}
We can similarly modify \eqref{eq:rt-none-pde} to account for mode-dependent $\gamma_i$
\begin{align}
- \frac{\partial v_m}{\partial t} &= G_m(\vx, t, \nabla v_m)
+ \sum_{n=1}^M b_n(t) \gamma_n \left(\phi_n(\vx) - v_m(\vx,t)\right),  &\vx& \in \Omega\setminus\Gamma, \, t \geq 0 \label{eq:rt-none-gamma-pde}\\
v_m(\vx, t) &= \overline{\psi}_m(\vx, t),  &\vx& \in \partial \Gamma, \, t \geq 0.
\end{align}

\begin{rem}
To model different observation schemes in addition to random terminations,
it is straightforward to define an operator
\begin{align}
R_m(\vx, t, v_m, \nabla v_m) &= \min_{\veca \in A} \left\{\overline{K}_m(\vx, t, \veca) + f(\vx, \veca) \cdot \nabla v_m(\vx, t)\right\} + \sum_{n=1}^M b_n(t) \gamma_n \left(\phi_n(\vx) - v_m(\vx,t)\right)
\end{align}
and replace $G_m$ with $R_m$ in any of the systems considered in subsections \ref{sec:indef-bdd} and \ref{sec:indef-paid-for}.
\end{rem}

\begin{rem}\label{rem:mode-and-state-gamma}
One significant limitation is our inability to handle the dependence of the termination rate on both the current state and (unobserved) mode simultaneously.
With  termination rates $\gamma_i(\x),$ the planner's past actions influence the current belief, and thus the full belief-space dynamic programming of Section \ref{sec:belief-prog} becomes necessary.
\end{rem}

\section{Numerical Methods}
\label{sec:numerics}
To simplify the exposition, we will discuss numerical methods for 2D occasionally observed PDMPs on a unit square (with possible obstacles).
We will focus on Eulerian discretizations for the following subclass of problems.
\begin{assumption}
The dynamics and running costs
are \emph{isotropic}: 
\begin{itemize}
\item
$\vb f(\vx, t, \veca) = \veca f(\vx, t)$, with a scalar speed function $f(\vx, t)$\\ and $\veca \in \mathbb{S}^1$ specifying the chosen direction of motion;
\item
$K_i(\x, t, \ba) = K_i(\x, t)$ for all  $i \in \M.$
\end{itemize}
Similarly, we assume $\vb f(\vx, \veca) = \veca f(\vx)$ and $ K_i(\x, \ba) = K_i(\x)$
in the autonomous case
\end{assumption}
Under these assumptions, the optimal policy associated with a value function $v$ can be recovered by setting $\veca = -\nabla v/|\nabla v|$. 
This allows us to rewrite the Hamiltonians for the control problems we consider, e.g., $G_m(\vx, t, \nabla v) = \overline{K}_m(\vx, t) - f(\vx,t)|\nabla v|$.

We will use a uniform $(J+1) \times (J+1)$ Cartesian grid on $[0,1] \times [0,1]$ with grid spacing $h =1/J$
and gridpoints $(x^i, y^j) = (ih, jh)$.
When discretizing time, the choice of grid spacing $\Delta t$ is governed by a standard Courant-Friedrichs-Lewy (CFL) condition: 
$\Delta t \leq \frac{h}{ \sqrt{2} \max_{\x,t}f(\vx, t)}.$ 
The resulting time slices are $t^k =  k \Delta t$ with $k = 0,..., N$ and $\sscap{t}{N} = T.$
A grid function $V$ will denote the numerical solution with  $V^{i,j,k} = V(x^i, y^j, t^k) \approx v(x^i, y^j, t^k)$. 
We will similarly use $f^{ijk}$, $K_n^{ijk}$, and $\psi_n^{ij}$ to denote the discretized dynamics, running costs, and terminal costs for all modes $n \in \M.$
(For the sake of readability, the commas in superscripts will be omitted whenever this does not create an ambiguity.)

Throughout this section, we make use of standard first-order upwind finite difference discretizations of Hamilton-Jacobi-Bellman PDEs,
which are convergent to viscosity solutions via a usual (monotonicity + consistency + stability) argument \cite{barles1991convergence, oberman2006convergent}. 
For all problem types, our discretizations rely on the following  
one-sided difference operators:
\begin{gather}
D^{ijk}_{\pm x} V= \frac{V^{i\pm1,j,k} - V^{i,j,k}}{\pm h}, \quad  D^{ijk}_{\pm y} V = \frac{V^{i,j\pm1,k} - V^{i,j,k}}{\pm h},
\quad \text{and} \quad
\mathcal{D}^{ijk}_{-t} V= \frac{V^{i,j,k} - V^{i,j,k-1}}{\Delta t}.
\end{gather}
If obstacles are present, their boundaries will be assumed to be grid aligned, and
\update{the gridpoints inside the obstacle boundaries will be entirely omitted (or, equivalently, their values will be assumed to be infinite).}
To apply the above operators on the boundary of the square, we will similarly assume that  $V^{ijk} = +\infty$ whenever $i$ or $j$ fall outside of the $0,...,J$ range. 
This will have the effect of approximating the domain-constrained viscosity solutions, and preventing optimal trajectories from leaving $\bar{\Omega}.$

We will also define the \emph{upwind} difference operators
\begin{align}
\mathcal{D}^{ijk}_x V = \min\{D^{ijk}_{+x} V, -D^{ijk}_{-x} V, 0 \} \quad \text{and} \quad \mathcal{D}^{ijk}_y V = \min\{D^{ijk}_{+y} V, -D^{ijk}_{-y} V, 0 \}, 
\end{align}
which (together with the CFL condition) ensure that the computational stencil used to approximate the gradient of $V$ straddles the characteristic corresponding to the optimal trajectory of the original control problem. 

\update{Since all of our PDE discretizations can be re-interpreted as Markov Decision Processes (MDPs) with infinitely many available actions on the corresponding grid-like graphs,
the convergence of all iterative algorithms follows from the standard value iterations convergence results for MDPs \cite{kushnerdupuis1992}.}  

\subsection{Finite Horizon} \label{sec:num-fin}
If the planner does not receive any observations after $t = 0$ (as in Section \ref{sec:reduced-b-prog}), we are interested in solving the equation \eqref{eq:fin-none-pde} with the Hamiltonian $G_m(\vx, t, \nabla v_m) = \overline{K}_m(\vx,t) - f(\vx, t) \, |\nabla v_m|,$ corresponding to the isotropic problem and $\bq = \vb e_m.$ 
A straightforward application of the upwind difference operators  leads to the discretized equations
\begin{align}
V_m^{i,j,k-1} &= V_m^{ijk} + \Delta t\left(\overline{K}_m^{ijk} - f^{ijk} \sqrt{\left(\mathcal{D}_x^{ijk} V_m\right)^2 + \left(\mathcal{D}_y^{ijk} V_m \right)^2} \right), \label{eq:fin-dis}\\
  V_m^{ijN} &= \overline{\psi}_m^{ijN}.
\end{align}
for all $m, i, j,$ and $k$. 
This results in an explicit, backwards in time scheme for computing $V_m$ for all possible starting modes $m.$ 
The same approach can be also used with any initial distribution $\vb q \in \Q.$
Going forward, we will use
\begin{equation}
G_m^{ijk} = \overline{K}_m^{ijk} - f^{ijk} \sqrt{\left(\mathcal{D}_x^{ijk} V_m\right)^2 + \left(\mathcal{D}_y^{ijk} V_m \right)^2}
\end{equation}
to denote our discretized Hamiltonian.

When the planner receives observations of the mode at known times $T_l$ (as in Section \ref{sec:finite}), we can make use of the same discretization.
Algorithm \ref{alg:fin-sched-solver} describes the process\footnote{
\update{Even though our current implementation and the pseudocodes included in this section are serial, we note that many of these computations are 
``perfectly parallelizable'' and could be sped up significantly, e.g., through computing $v_m(\vx, t)$ for all initial modes $m$ in parallel.}}
 used to solve for $V_m^l$ on each interval $\left[ 0, \, \left( T_{l+1} - T_l \right) \right]$.
To account for the shift in time on each subinterval, we use $f^{ijk} = f(x^i, y^j, T_l + k\Delta t)$ and $\overline{K}_m^{ijk} = \sum_{n=1}^M b_n^k K_n(x^i, y^j, T_l + k \Delta t)$.
Recalling the definition of $\Theta$ in the continuous case, we write $[\Theta_m^l]^{ijk} = \sum_{n=1}^M b_n^k [V_n^l]^{ij0}$, which determines the terminal condition for all $l \neq L$.
\update{The overall computational cost is $O(J^2 N)$ without mode observations and $O\left( J^2 N \left[ M \left(1 - \frac{T_1 - T_0}{T} \right)  \, +  \frac{T_1 - T_0}{T} \right] \right)$ with mode 
observations scheduled for $T_1,  \ldots, T_L < T.$}  

\begin{algorithm}[t]
\SetKw{Compute}{Compute}
\SetKw{Solve}{Solve}
\SetKw{Initialize}{Initialize}
\SetKw{Set}{Set}

\KwIn{$[T_1, \hdots, T_L]$, a set of observation times}

\Set{$T_0 = 0 \text{ and } T_{L+1} = T$}

\For{$l = L, L-1, \hdots, 0$}{
	\Set{$N^l = (T_{l+1} - T_l)  / \Delta t$}

	\For{$m=1, \ldots,M$}{
		\uIf{l = L} {
			\Set{$[V_m^L]^{i,j,N^L} = \overline{\psi}_m^{i,j,N^L}$} 
		}
		\Else{
			\Set{$[V_m^l]^{i,j,N^l} = [\Theta_m^{l+1}]^{i,j,N^l}$}
		}
		\Compute{$[V_m^l]^{i,j,k-1}$ via \eqref{eq:fin-dis} for $k=N^l,\hdots,1; \quad i,j = 0, \hdots, J$.}
	}
}
\caption{Finite horizon  with  scheduled observations solver}
\label{alg:fin-sched-solver}
\end{algorithm}

\subsection{Infinite Horizon, Periodic Observations}
\begin{algorithm}[t]
\SetKw{Compute}{Compute}
\SetKw{Solve}{Solve}
\SetKw{Initialize}{Initialize}
\SetKw{Set}{Set}

\KwIn{$tol$, a termination threshold}

\Set{$l = 0$, $\delta = tol + 1$}

\While{$\delta > tol$}{
	\For{$m=1, \ldots,M$}{
		\uIf{l = 0}{
		\Set{$[V_m^l]^{ijN} = K_m^{ij}/\beta$}   \hspace*{3mm}(The cost of staying in place indefinitely without switches.) 
		}
		\Else{
		\Set{$[V_m^l]^{ijN} = [\Theta_m^{l-1}]^{ijN}$}
		}
		\Compute{$[V_m^l]^{i,j,k-1}$ via \eqref{eq:inf-per-dis} for $k=N,\hdots,1; \quad i,j = 0, \hdots, J$.}
	}

	\Set{$\delta = \max\limits_{m\in\mathcal{M}} \|V_m^{l} - V_m^{l-1}\|_{\infty}$}\\
	\Set{$l = l+1$}\\
}

\caption{Infinite horizon with periodic observations solver}
\label{alg:inf-iter-solver}
\end{algorithm}

The setting described in Section \ref{sec:infinite} presents another problem defined over a finite (inter-observation) time interval $[0,T].$
However, due to the nonlocal coupling in terminal conditions, it is no longer possible to compute the solution in a single sweep backwards in time.
We take an iterative approach, where we still solve backwards in time at each iteration, using the previous best approximation of the value function to compute the terminal condition.
This scheme can be interpreted as solving backwards in time for an increasing number of periods; we use superscript $l$ to indicate the current iteration and highlight the number of observation periods already considered.
Using the same upwind difference operators, we arrive at the following discretized update equations for the $l$-th iteration
\begin{align}
\left[V^l_m\right]^{i,j,k-1} &= \left[V^l_m\right]^{ijk} + \Delta t\left(e^{-\beta t^k} \overline{K}_m^{ijk} - f^{ij} \sqrt{\left(\mathcal{D}_x^{ijk} V^l_m \right)^2 + \left(\mathcal{D}_y^{ijk} V^l_m \right)^2} \right), \label{eq:inf-per-dis}\\
\left[V^l_m\right]^{ijN} &= 
\left[\Theta_m^{l-1}\right]^{ijN}.\label{eq:inf-per-dis-tc}
\end{align}
The iterative approach we take is outlined in Algorithm \ref{alg:inf-iter-solver}.
\update{The computational cost of each iteration is $O(J^2 N M).$}

\subsection{Indefinite Horizon}

Shifting to the setting of Section \ref{sec:indef}, we will assume that the boundaries of the target set $\Gamma$ are grid-aligned
and will treat all gridpoints in the interior of $\Gamma$ as if their values are infinite.
When there are no observations after $t = 0$, we are interested in approximating the solution of 
\begin{align*}
       - \frac{\partial v_m}{\partial t}(\vx,t)  &= \overline{K}_m(\vx,t) - f(\vx) \, |\grad v_m(\vx,t)|, &\vx& \in \Omega\setminus \Gamma, \, t \in [0, T);\\
        v_m(\vx,t) &= \overline{\psi}_m(\vx, t), &\vx& \in \partial \Gamma, \, t  \in [0, T);\\
        v_m(\vx,T) &= +\infty, &\vx& \not\in \Gamma,
\end{align*}
where the horizon $T$ is chosen based on Proposition \ref{prop:bounded-TG}.
For the case with no observations, the discretization is identical to equation \eqref{eq:fin-dis} (though with a stationary $f$), but 
with additional boundary and terminal conditions
\begin{align}
\label{eq:discr_target_boundary}
V_m^{i,j,k-1} &= 
\min \left\{
\overline{\psi}_m^{i,j,k-1}, \;
V_m^{ijk} + \Delta t G_m^{ijk} 
\right\},
& \text{ for } (x^i,y^j) \in \partial \Gamma, \, k \geq 0;\\
V_m^{ijN} &= 
\text{INF}, 
& \text{ for } (x^i,y^j) \in \Omega\setminus\Gamma. \label{eq:indef-dis-bc-tc}
\end{align}
The $\min$ in \eqref{eq:discr_target_boundary} represents our interpreting the boundary conditions {\em in viscosity sense}
\cite[Chapter 5]{bardi2008}.
On $\partial \Gamma,$ the planner can either enter the interior of $\Gamma$ immediately (with the expected terminal cost of  $\overline{\psi}_m$) or
 continue the process without entering, possibly to secure a better terminal cost later.   
In \eqref{eq:indef-dis-bc-tc}, INF is a large number approximating infinity.  

\begin{rem}
While the viscosity solution is not affected by this truncation of the infinite time interval in the sense outlined in Section \ref{sec:indef}, 
if the bound $T$ is tight and the grid is coarse, the numerical viscosity of our discretization might 
result in \eqref{eq:indef-dis-bc-tc} affecting the approximate solution even on a relevant part of the space-time grid.
However these artifacts disappear under the grid refinement and can be also easily avoided by using a larger $T$.
\end{rem}

\begin{algorithm}[t]
\SetKw{Compute}{Compute}
\SetKw{Solve}{Solve}
\SetKw{Initialize}{Initialize}
\SetKw{Set}{Set}

\KwIn{$T$, an upper bound on the time until the goal}

\For{$l = L, \hdots, 0$}{
	\For{$m=1, \ldots,M$}{
	\Set{$[V_m^l]^{ijN} = $ INF for $(x^i, y^j)\notin \Gamma$.}\\
	\For{$k =N, \hdots, 1$}{
		\uIf{l = L} {
			\Compute{$[V_m^l]^{i,j,k-1}$ via \eqref{eq:fin-dis} for $(x^i,y^j)\notin \Gamma$ 
			and via \eqref{eq:discr_target_boundary} for $(x^i,y^j) \in \partial \Gamma.$}}
		\Else{
			\Compute{$[V_m^l]^{i,j,k-1}$ via \eqref{eq:indef-bdd-dis} for $(x^i,y^j)\notin \Gamma$
			and via \eqref{eq:indef-bdd-dis-gamma} for $(x^i,y^j) \in \partial \Gamma.$}}
		}
	}
}

\caption{Indefinite horizon  with limited observations}
\label{alg:indef-bdd-solver}
\end{algorithm}

Turning to the case of limited on-demand observations (as in Section \ref{sec:indef-bdd}),
we can see that when the number of remaining observations is $(L-l),$
the PDE for $v_m^l$ depends only on $v_n^{l+1}$ for all $n \in \M.$
This means that we can solve the problem non-iteratively (looping through $l = L, ..., 0$).
For $0 \leq l < L$, the discrete update equations are
\begin{align}
\left[V^l_m\right]^{i,j,k-1} &= \min \left \{\left[V^l_m\right]^{ijk} + \Delta t \left[G_m^l\right]^{ijk}, \,
\left[\Theta_m^{l+1}\right]^{ijk} \right\},
& \text{on } \Omega \setminus \Gamma; \label{eq:indef-bdd-dis}\\
\left[V^l_m\right]^{i,j,k-1} &= \min \left \{\left[V^l_m\right]^{ijk} + \Delta t \left[G_m^l\right]^{ijk}, \,
\left[\Theta_m^{l+1}\right]^{ijk}, \, \overline{\psi}_m^{i,j,k-1} \right\},
& \text{on } \partial \Gamma. \label{eq:indef-bdd-dis-gamma}
\end{align}
Thus, we must approximate $M(L + 1)$ value functions, with each of them computed non-iteratively (backwards in time). 
This approach is outlined in Algorithm \ref{alg:indef-bdd-solver}, 
\update{with the overall computational cost of $O \left(J^2 N \left(M L+1 \right) \right).$} 

\begin{algorithm}[t]
\SetKw{Compute}{Compute}
\SetKw{Solve}{Solve}
\SetKw{Initialize}{Initialize}
\SetKw{Set}{Set}

\KwIn{$tol$, a termination threshold}
\KwIn{$T$, an upper bound on the time until the goal}

\Set{$l = 0$, $\delta = tol + 1$}

\While{$\delta > tol$}{
	\For{$m=1, \ldots,M$}{
	\Set{$[V_m^l]^{ijN} = \infty$ for $(x^i, y^j)\notin \Gamma$.}\\
	\For{$k =N, \hdots, 1$}{
		\uIf{l = 0} {
			\Compute{$[V_m^l]^{i,j,k-1}$ via \eqref{eq:fin-dis} for $(x^i,y^j)\notin \Gamma$ 
			and via \eqref{eq:discr_target_boundary} for $(x^i,y^j) \in \partial \Gamma.$}}
		\Else{
			\Compute{$[V_m^l]^{i,j,k-1}$ via \eqref{eq:indef-paid-dis} for $(x^i,y^j)\notin \Gamma$
			and via \eqref{eq:indef-paid-dis-gamma} for $(x^i,y^j) \in \partial \Gamma.$}}
		}
	}
	\Set{$\delta = \max\limits_{m\in\mathcal{M}} \|V_m^{l} - V_m^{l-1}\|_{\infty}$}\\
	\Set{$l = l+1$}\\
}

\caption{Indefinite horizon with paid (unlimited) observations}
\label{alg:indef-iter-solver}
\end{algorithm}

When observations 
can be instead purchased for a cost $C(\vx)$, the sequential method described above is no longer possible.
We instead take an iterative approach, outlined in Algorithm \ref{alg:indef-iter-solver}.
When observations are purchased for a cost $C,$ within the $l$-th iteration the discretized equations are instead
\begin{align}
\left[V^l_m\right]^{i,j,k-1} &= \min \left \{\left[V^l_m\right]^{ijk} + \Delta t\left[ G_m^l\right]^{ijk}, \, 
C^{ij} + \left[\Theta_m^{l-1}\right]^{ijk}\right\},
& \text{on } \Omega \setminus \Gamma;  \label{eq:indef-paid-dis}\\
\left[V^l_m\right]^{i,j,k-1} &= \min \left \{\left[V^l_m\right]^{ijk} + \Delta t\left[ G_m^l\right]^{ijk}, \, 
C^{ij} + \left[\Theta_m^{l-1}\right]^{ijk}, \,  \overline{\psi}_m^{i,j,k-1} \right\},
& \text{on } \partial \Gamma. \label{eq:indef-paid-dis-gamma}
\end{align}
\update{The computational cost of each iteration is again $O(J^2 N M).$}

\subsection{Randomly Terminated Problems}

As with the continuous case described in Section \ref{sec:rt}, we can handle randomly terminated problems 
by simply modifying the Hamiltonian for existing idefinite horizon problems.
Discretizing the continuous randomly terminated Hamiltonian we obtain
\begin{equation}
\left[R_m^l\right]^{ijk} = \overline{K}_m^{ijk} - f^{ij} \left(\sqrt{\left(\mathcal{D}_x^{ijk} V^l_m\right)^2 + \left(\mathcal{D}_y^{ijk} V^l_m\right)^2} \right) + \sum_{n=1}^M b_n^k \gamma^{}_n \left(\phi_n^{ij} - \left[V_m^l\right]^{ijk}\right).
\end{equation}
We can thus use Algorithms \ref{alg:indef-bdd-solver} and \ref{alg:indef-iter-solver} with $R_m^l$ replacing $G_m^l$ for the problems of interest.


\section{Numerical Experiments}
\label{sec:numexp}
Most of our test problems (except in subsection \ref{sec:rover}) are motivated by a security application:
an evader that seeks to minimize their cumulative exposure to surveillance while navigating through $\Omega.$
The modes are interpreted as ``surveillance patterns'', with each $K_i(\x)$ encoding the exposure at a location $\x$ when the $i$-th pattern is in effect.
(We will usually model each $K_i$ as a Gaussian centered on some source of surveillance $\hat{\x}_i$; e.g. a security camera, a watchtower, or a drone base.)
We will also generally set $\psi_i(\x) = 0$ across examples, since we assume that exposure to surveillance ends when the process stops.
We note that this is a very natural application for our framework since the evader often has no idea of the surveillance pattern currently in effect
and has no way of directly measuring the cumulative exposure.   We assume that the evader has access to data on how the surveillance patterns were shifting in the past,
allowing the estimation of pattern switching rates $\lambda_{ij}.$  When we consider versions with mode observations, we view them as updates on the current surveillance patterns obtained directly by the evader or through allies (e.g., satellite imagery). 

Unless otherwise specified, all experiments are conducted on the spatial domain $[0,1] \times [0,1]$ and discretized on a $501 \times 501$ Cartesian grid. 
We use the smallest number of time slices  $(N+1)$ that satisfies the CFL condition. 

\subsection{Rotating Surveillance}\label{sec:fin-results}

\begin{figure}
\centering
\begin{subfigure}{0.28\textwidth}
\centering
\caption{Surveillance Patterns}\label{fig:4modes-labels}
\includegraphics[height=4.5cm, trim={20 15 0 0}, clip]{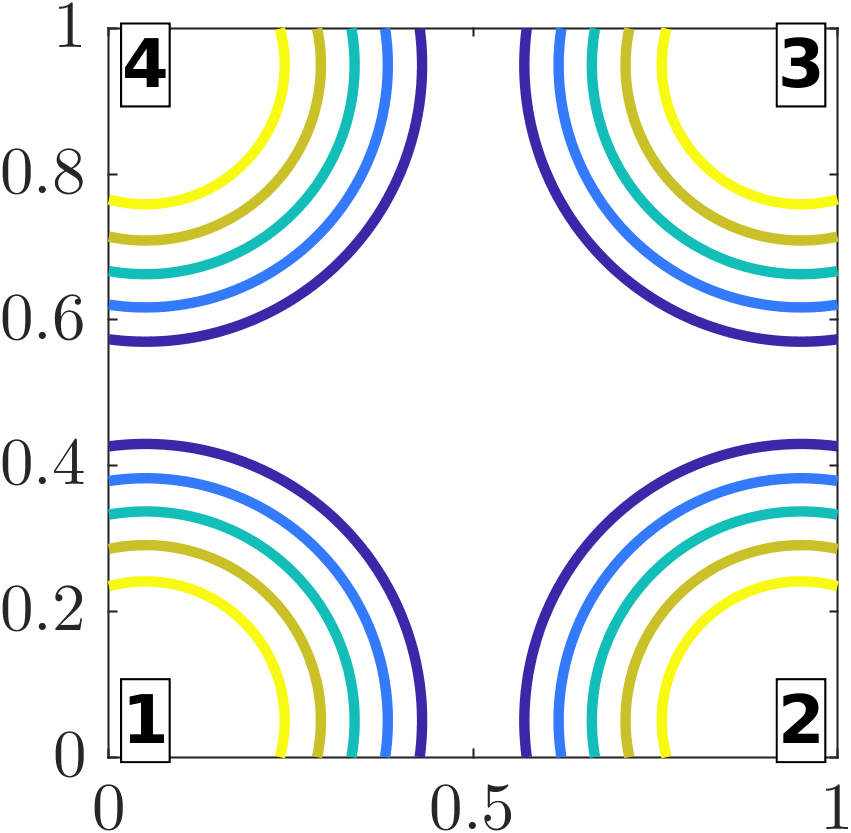}
\end{subfigure}
\hfill
\begin{subfigure}{0.35\textwidth}
\centering
\caption{$K_1(\vx)$}\label{fig:4modes-K1}
\includegraphics[height=4.5cm, trim={22 15 0 0}, clip]{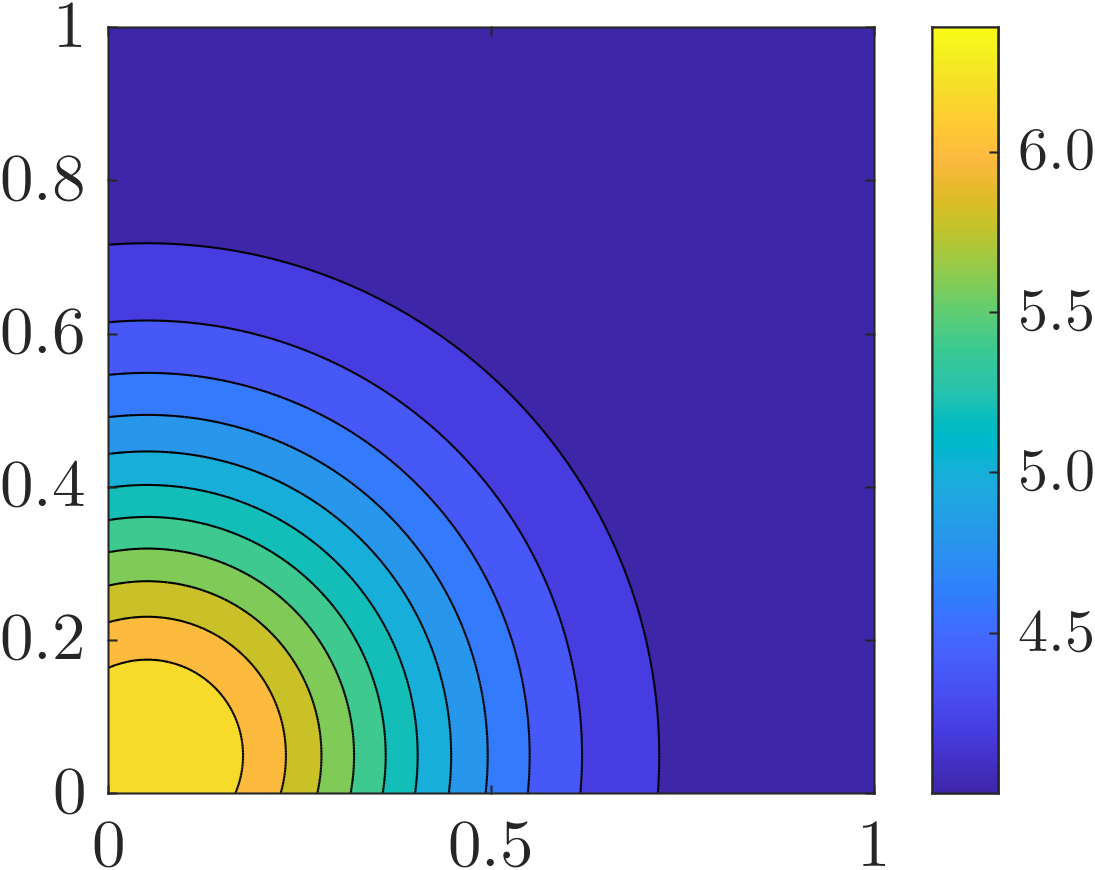}
\end{subfigure}
\hfill
\begin{subfigure}{0.35\textwidth}
\centering
\caption{$\overline{K}_s(\vx)$}\label{fig:4modes-Ks}
\includegraphics[height=4.5cm, trim={22 15 0 0}, clip]{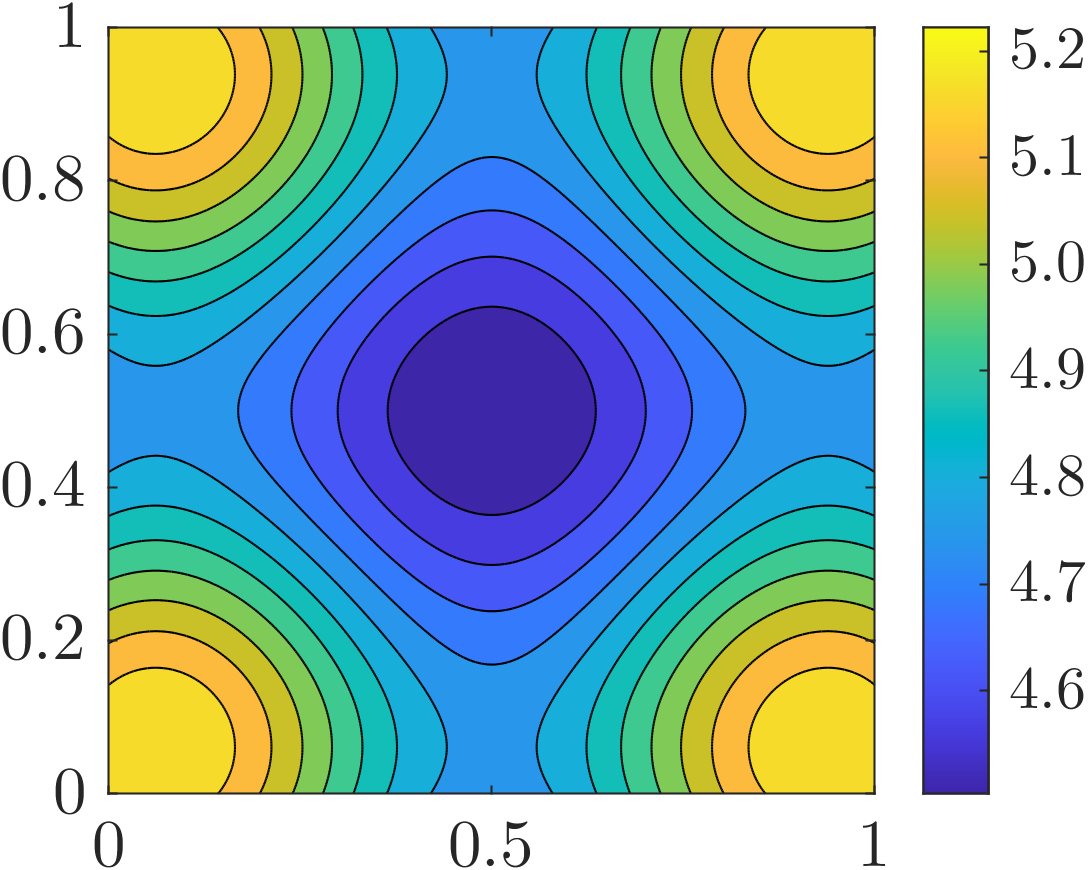}
\end{subfigure}
\caption{``Rotating Surveillance'' environment. \textbf{(a)} Labeled contour plots of mode surveillance patterns $K_i(\vx)$. Mode labels (boxed numbers) are placed at the peak of each $K_i$. Only one pattern is ``active'' at a time. \textbf{(b)} $K_1(\vx)$, the surveillance pattern in Mode 1. \textbf{(c)} $\overline{K}_s(\vx)$, the expected surveillance associated with the stationary (uniform) mode distribution.}
\label{fig:4modes-cost}
\end{figure}

The first example is an environment with four possible surveillance patterns (Figure \ref{fig:4modes-cost}), each with corresponding surveillance intensity 
(running cost) $K_i(\vx) = 4 + 9(2\pi\sigma)^{-1} \exp(-(2\sigma^{2})^{-1}(\vx - \hat{\vx}_i)(\vx - \hat{\vx}_i)^{\top})$ for $\sigma = 0.3$ and $\hat{\x}_i \in \left\{ (0.05,0.05), (0.95, 0.05), (0.95, 0.95), 0.05,0.95)\right\}$.
We suppose that the adversary is rotating counterclockwise through these patterns, as encoded by the rate matrix
\begin{equation}
\Lambda = \begin{bmatrix}
-1 &  1 &  0 &  0 \\
 0 & -1 &  1 &  0 \\
 0 &  0 & -1 &  1 \\
 1 &  0 &  0 & -1 
\end{bmatrix}.
\end{equation}
While these transitions occur with equal rates and along a fixed cycle $1 \rightarrow 2 \rightarrow 3 \rightarrow 4 \rightarrow 1 \rightarrow \hdots$, the time until the the next pattern switch is still a random variable.
The corresponding CTMC has a stationary distribution $\vb q_s = [1/4, 1/4, 1/4, 1/4],$ and we use 
$\overline{K}_s(\vx)$ to denote the expected running cost associated with $\vb q_s$.
For this simple example, we take the speed to be constant ($f(\vx) = 1$) 
and investigate the impact of anticipated observations on optimal trajectories.

\begin{figure}
\centering
\hspace{1.8cm} Subinterval: $[0,1]$  \hspace{.4cm} Subinterval: $[0,2]$ \hspace{.4cm} Subinterval: $[0,3]$  \hspace{.4cm} Subinterval: $[0,4]$\\
\begin{subfigure}{\textwidth}
\begin{minipage}[c]{2.2cm}\raggedright\mbox{}\\[-\baselineskip]No mode observations\end{minipage}
\raisebox{-0.5\height}{\includegraphics[height = 3.33 cm, trim={0 0 0 0}, clip]{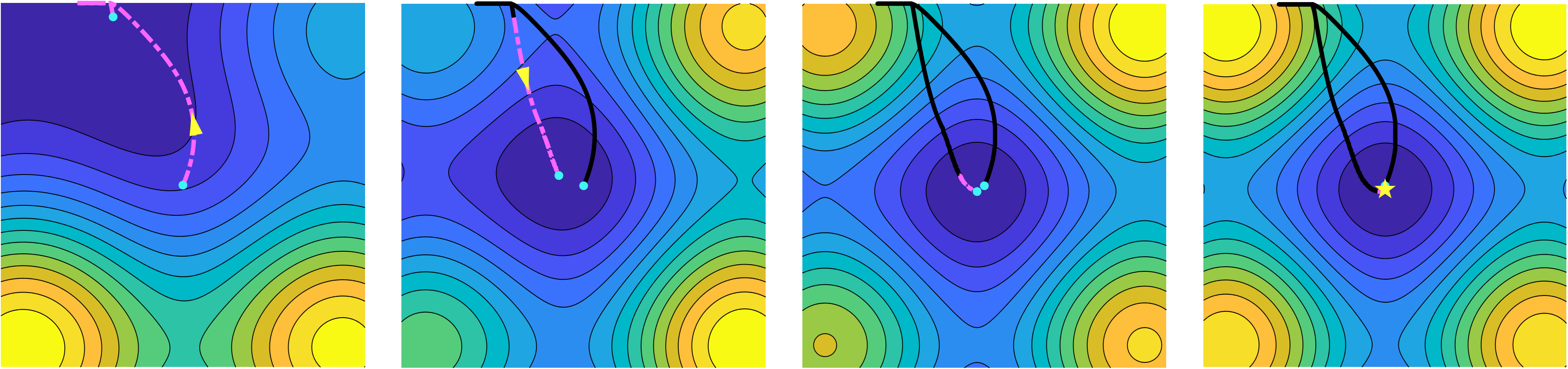}}
\end{subfigure}\\[5pt]
\begin{subfigure}{\textwidth}
\begin{minipage}[c]{2.2cm}\raggedright\mbox{}\\[-\baselineskip]Mode observations at $t \in \{1, 2, 3\}$\end{minipage}
\raisebox{-0.5\height}{\includegraphics[height = 3.31 cm, trim={0 0 0 0}, clip]{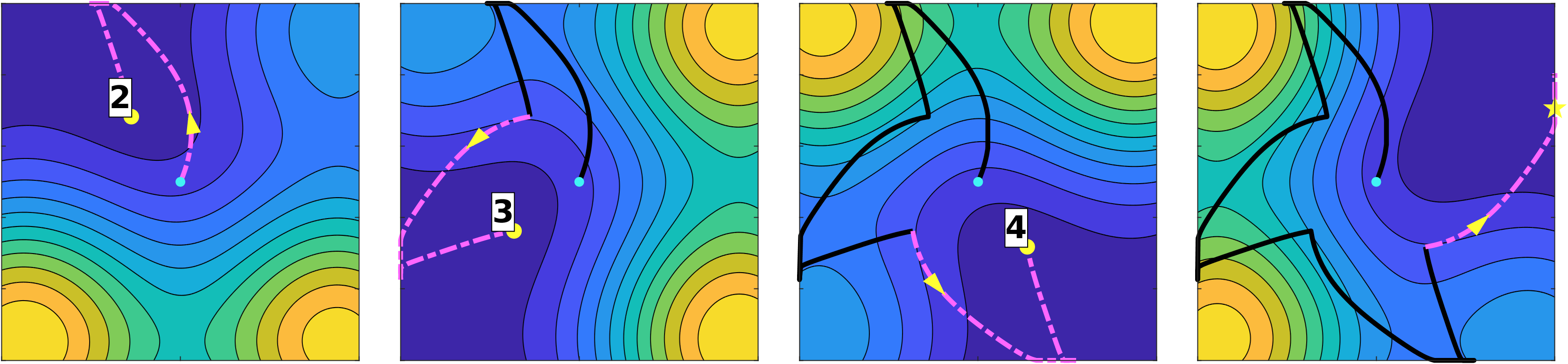}}
\end{subfigure}\\
\caption{Optimal trajectories for a finite horizon process without (top row) and with (bottom row) mode observations.
Time evolution is shown across columns, with magenta dash-dotted lines representing path components new to each column and black solid lines encoding path components shown previously. Cyan dots mark the planner's initial position and position at the end of each subinterval if an observation is not received. Yellow dots indicate the latest observations, with boxed numbers specifying the observed mode, yellow arrows indicate the direction of travel, and yellow stars indicate the planner's final position. The background is the expected surveillance at the end of each subinterval given $\mu(0) = 1$ and any other received mode observations. 
}
\label{fig:4modes-fin-traj-1}
\end{figure}

We first consider a fixed horizon process with and without mode observations (Sections \ref{sec:reduced-b-prog} and \ref{sec:finite}).
The evader remains in the domain for $t \in [0,4]$, starts at $\vy(0) = (0.5, 0.5)$, and knows the initial mode $\mu(0) = 1$.
Figure \ref{fig:4modes-fin-traj-1} compares the evolution of  optimal trajectories without mode observations (top row) and with mode observations at times $t \in \{1, 2, 3\}$ (bottom row).
In both cases, the planner initially travels to 
the upper left corner, where $K_1$ is low and $K_{\mu(t)}$ is not expected to be high for a long time.  
(Adversary will have to switch through $K_2$ and $K_3$ before getting to $K_4.$)
Without further mode observations (top row), the planner then returns to the center of the domain, where $\overline{K}_s$ has a global minimum.
On the other hand, in the bottom row the impact of mode observations can be seen before the first observation occurs.
In the first column, the planner with access to mode observations begins to travel towards the center of the domain earlier than would otherwise be optimal.
This is because a central location provides a better position from which to react to information gained from the upcoming mode observation.
The planner is then able to exploit the information gained to travel to areas 
where the exposure is expected to stay low for a long time based on the last observed mode.
In the final time interval, the trajectory remains close to the boundary, since there are no remaining observations to react to.
The shown trajectory corresponds to a specific sequence of observations, but the described properties are generic across all possible observation sequences. 

\begin{figure}[t]
\centering
\hspace{0cm} Subinterval: $[0,1]$  \hspace{.4cm} Subinterval: $[0,2]$ \hspace{.4cm} Subinterval: $[0,3]$  \hspace{.4cm} Subinterval: $[0,4]$\\
\includegraphics[height = 3.3cm, trim={0 0 0 0}, clip]{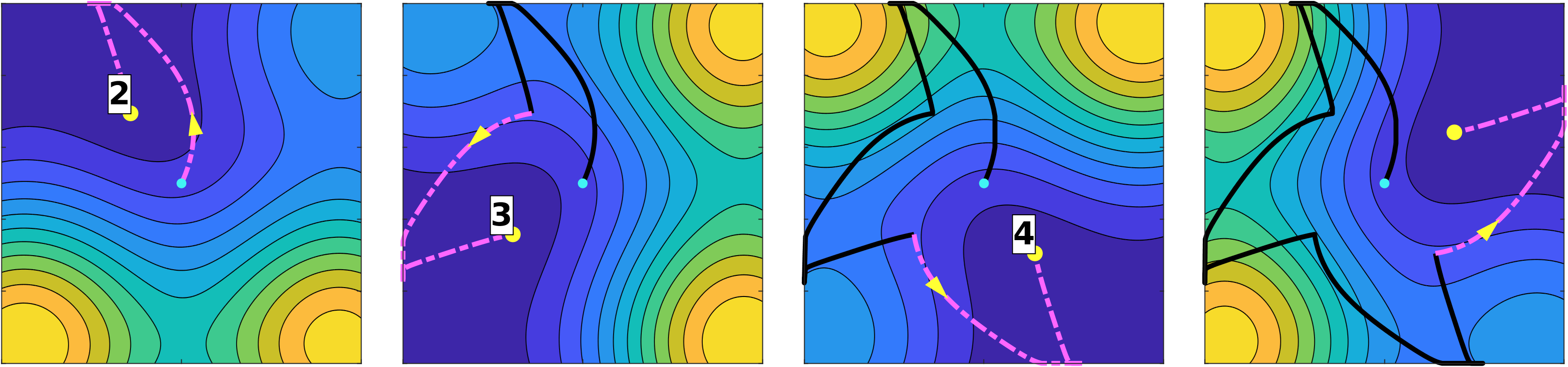}
\caption{
Optimal trajectories for an infinite horizon process with periodic mode observations.
Inter-observation period $T=1.$  Time discounting factor $\beta = 0.5.$  
Same visual format as in Figure \ref{fig:4modes-fin-traj-1}. \update{The solver required 19 iterations to converge with a tolerance of $10^{-6}.$
}
}
\label{fig:4modes-inf-traj-1}
\end{figure}

\begin{figure}[hbt]
\centering
\begin{subfigure}{0.30\textwidth}
\centering
\caption{$\beta = 0.5$}
\includegraphics[height=4.5cm, trim={22 15 40pt 0}, clip]{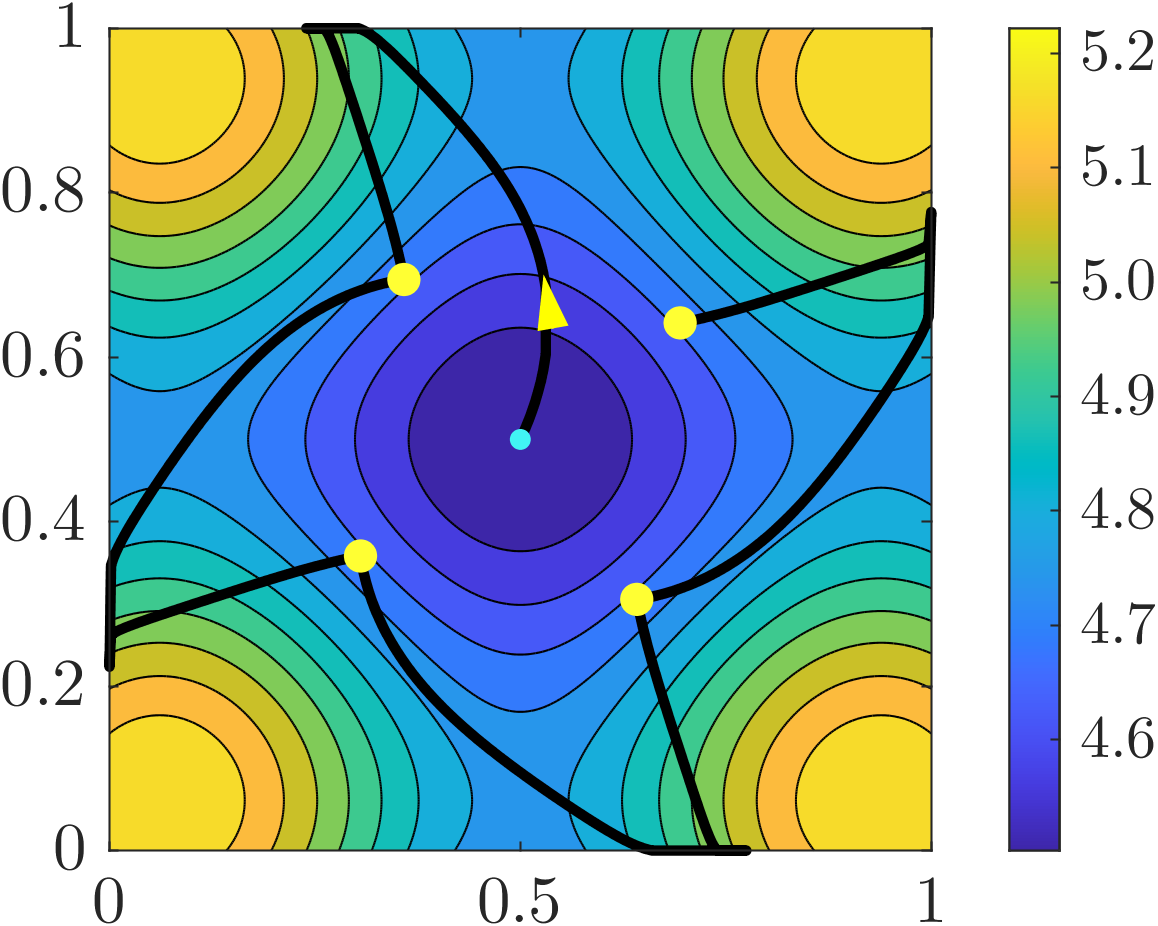}
\end{subfigure}
\begin{subfigure}{0.30\textwidth}
\centering
\caption{$\beta = 4$}
\includegraphics[height=4.5cm, trim={22 15 40pt 0}, clip]{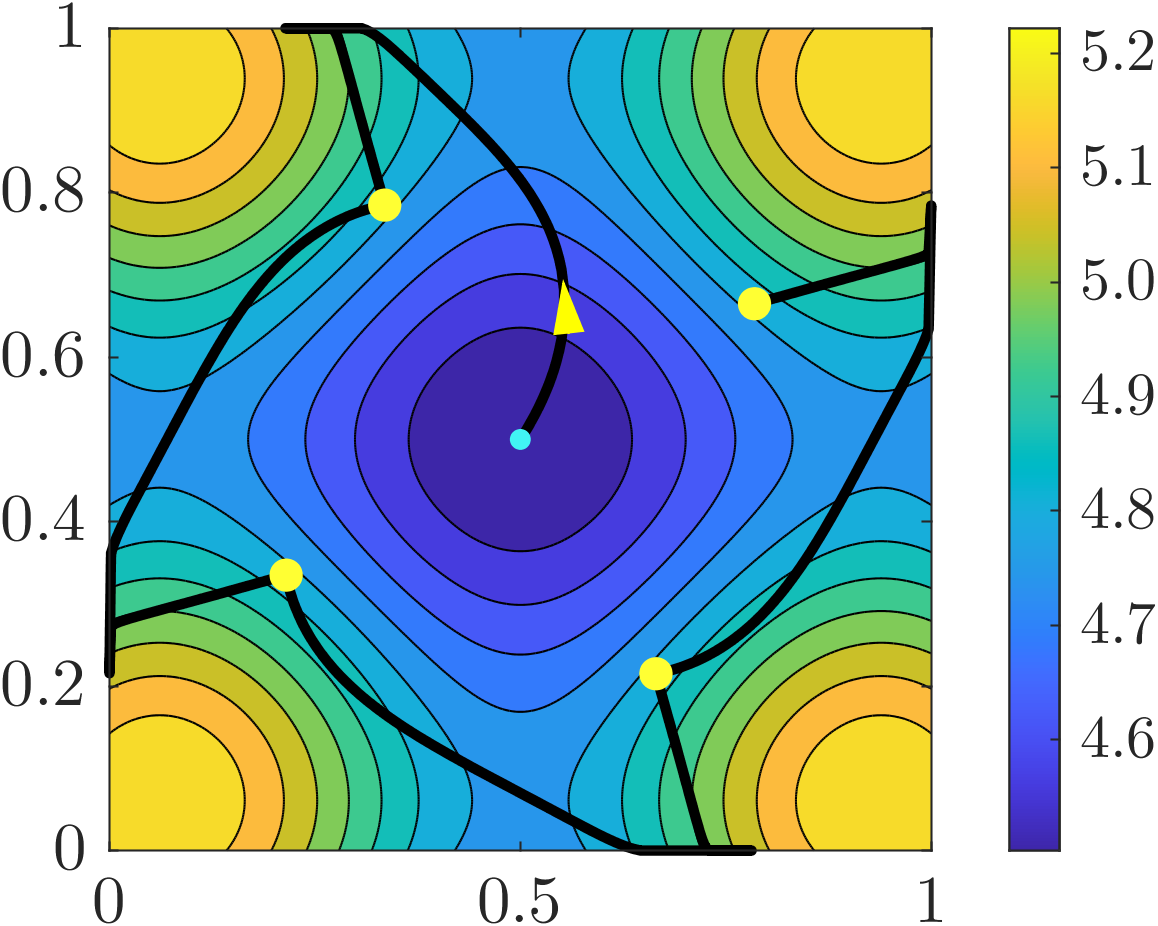}
\end{subfigure}
\begin{subfigure}{0.32\textwidth}
\centering
\caption{$\beta = 6$}
\includegraphics[height=4.5cm, trim={22 15 0 0}, clip]{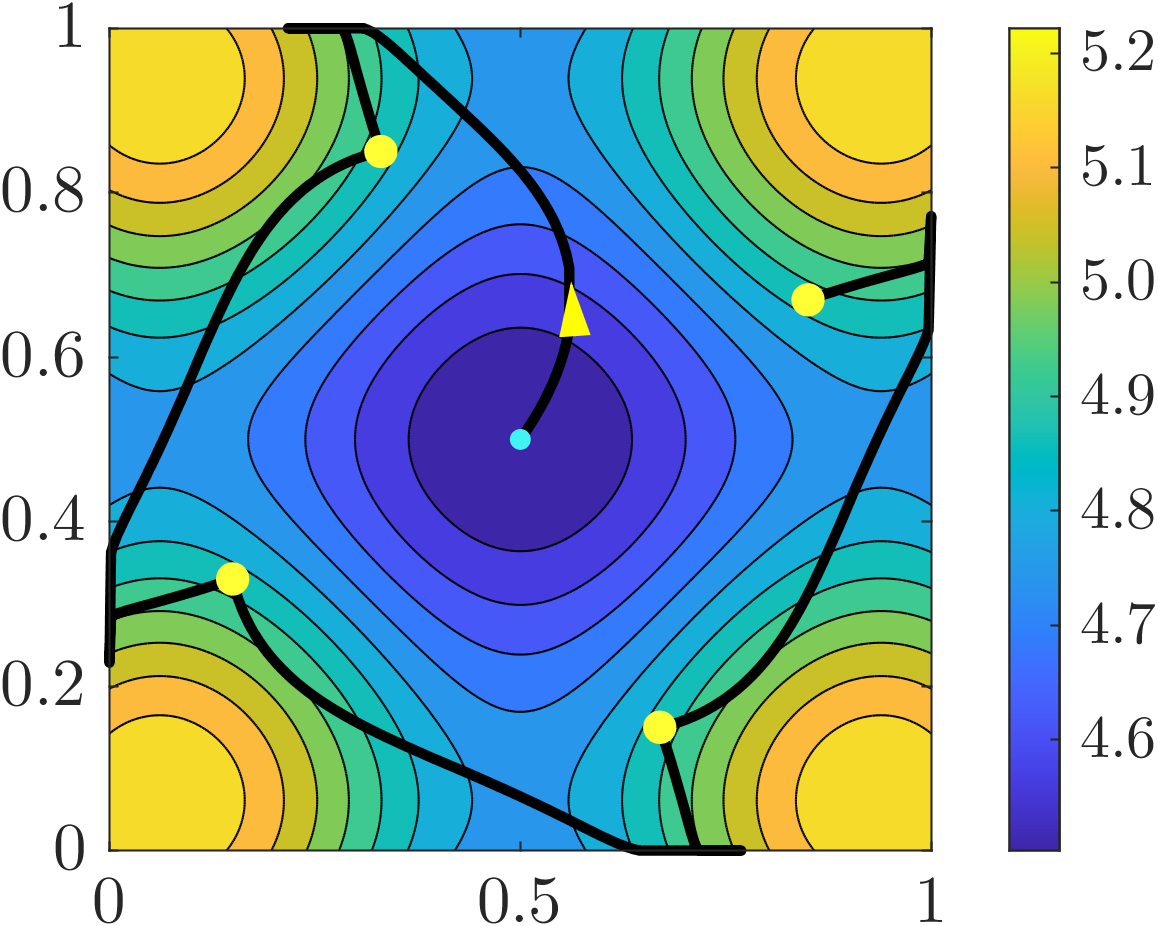}
\end{subfigure}
\caption{Optimal trajectories for infinite horizon process with periodic observations with three possible discount factors. Trajectories are shown for $t \in [0,4]$, corresponding to four periods. Cyan dots mark the planner's initial position. Yellow dots indicate observations and the observed modes are $\mu(1) = 2$, $\mu(2) = 3$, and $\mu(3) = 4$. Yellow arrows indicate the direction of travel. The background is $\overline{K}_s(\vx)$. \update{As the discount rate $\beta$ increases, the future impacts the planner to a lesser degree, and the number of iterations needed to reach convergence decreases. When $\beta = 6$, the solver requires just two iterations (compared to the 19 above) to converge to within a tolerance of $10^{-6}.$ 
}}
\label{fig:4modes-discounting}
\end{figure}

We next consider a discounted infinite horizon process with periodic observations (Section \ref{sec:infinite}) in the same environment.
We assume that the observation period is $T = 1$, so the planner receives observations at $t = 1, 2, 3, \hdots$.
Figure \ref{fig:4modes-inf-traj-1} shows a realized trajectory with the same observed modes as those in Figure \ref{fig:4modes-fin-traj-1}.
The behavior of these two planners is quite similar, 
but the infinite horizon process continues beyond the four periods that we show, and thus the planner still returns to the center of the domain in the last subinterval, in anticipation of the next observation.
Additionally, the optimal trajectory is impacted by the discount factor $\beta$.
Figure \ref{fig:4modes-discounting} shows optimal trajectories for the same sequence of observed modes but for three different values of $\beta$.
As $\beta$ increases, the planner dreads the future exposure less relative to the present exposure, and spends more time in areas with low current expected surveillance.
As a result, it takes 
smaller excursions towards the center of the domain in advance of each mode observation.

\subsection{Avoiding Barriers}
\begin{figure}[h]
\centering
\begin{subfigure}{0.25\textwidth}
\centering
\caption{$K_1(\vx)$}
\includegraphics[height=4.5cm, trim={20 15 50pt 0}, clip]{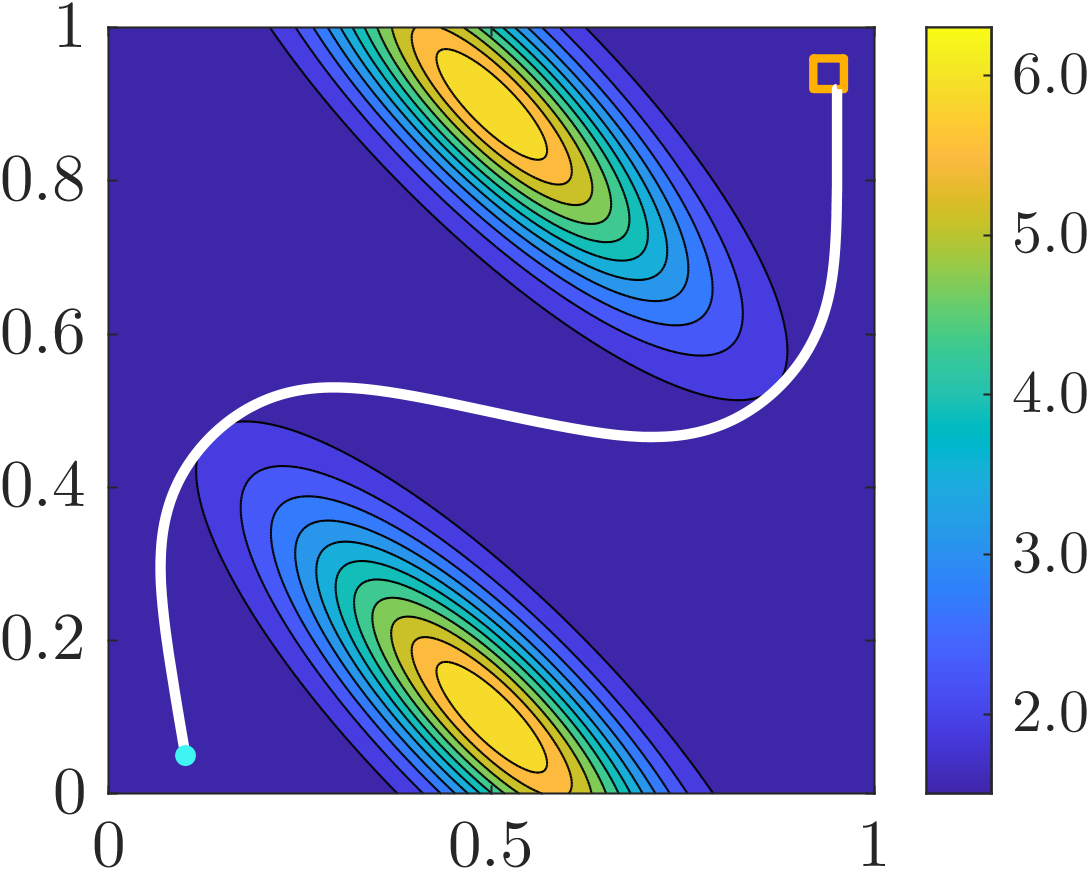}
\end{subfigure}
\begin{subfigure}{0.38\textwidth}
\centering
\caption{$K_2(\vx)$}
\includegraphics[height=4.5cm, trim={22 15 0 0}, clip]{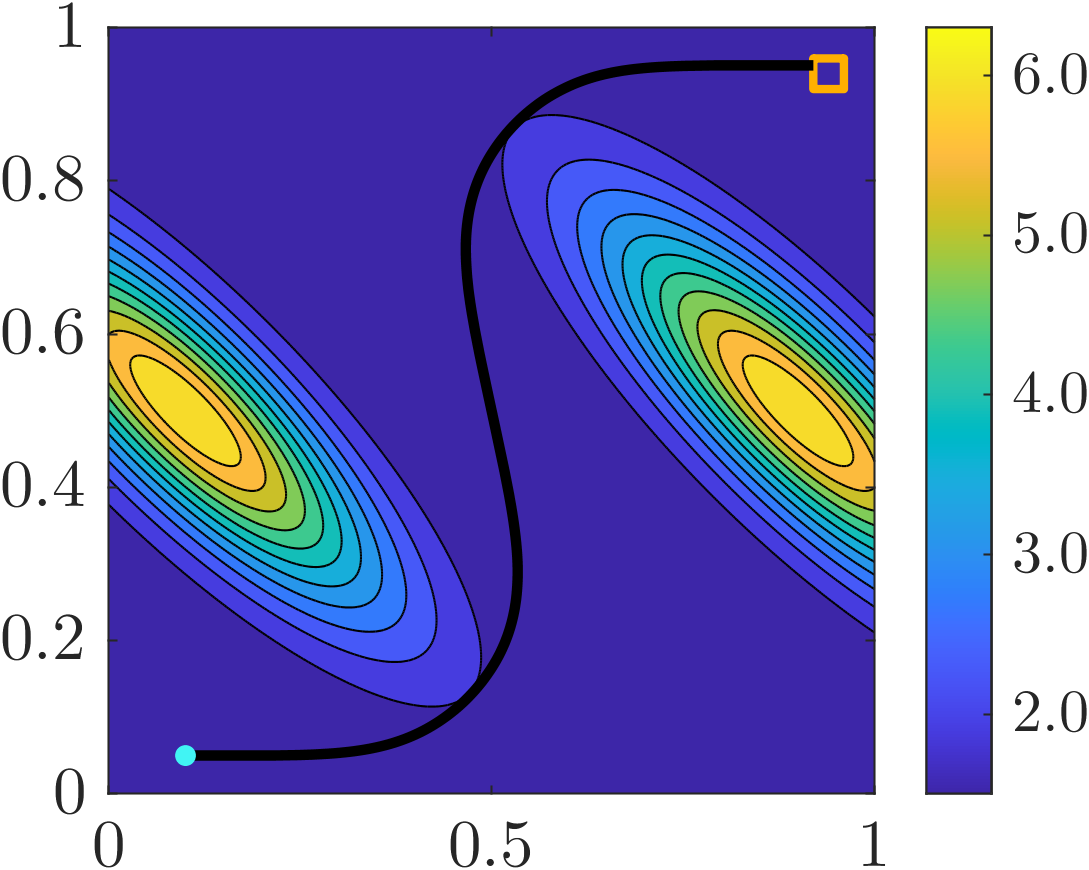}
\end{subfigure}
\begin{subfigure}{0.35\textwidth}
\centering
\caption{$\overline{K}_s(\vx)$}\label{fig:barriers-Kstat}
\includegraphics[height=4.5cm, trim={22 15 0 0}, clip]{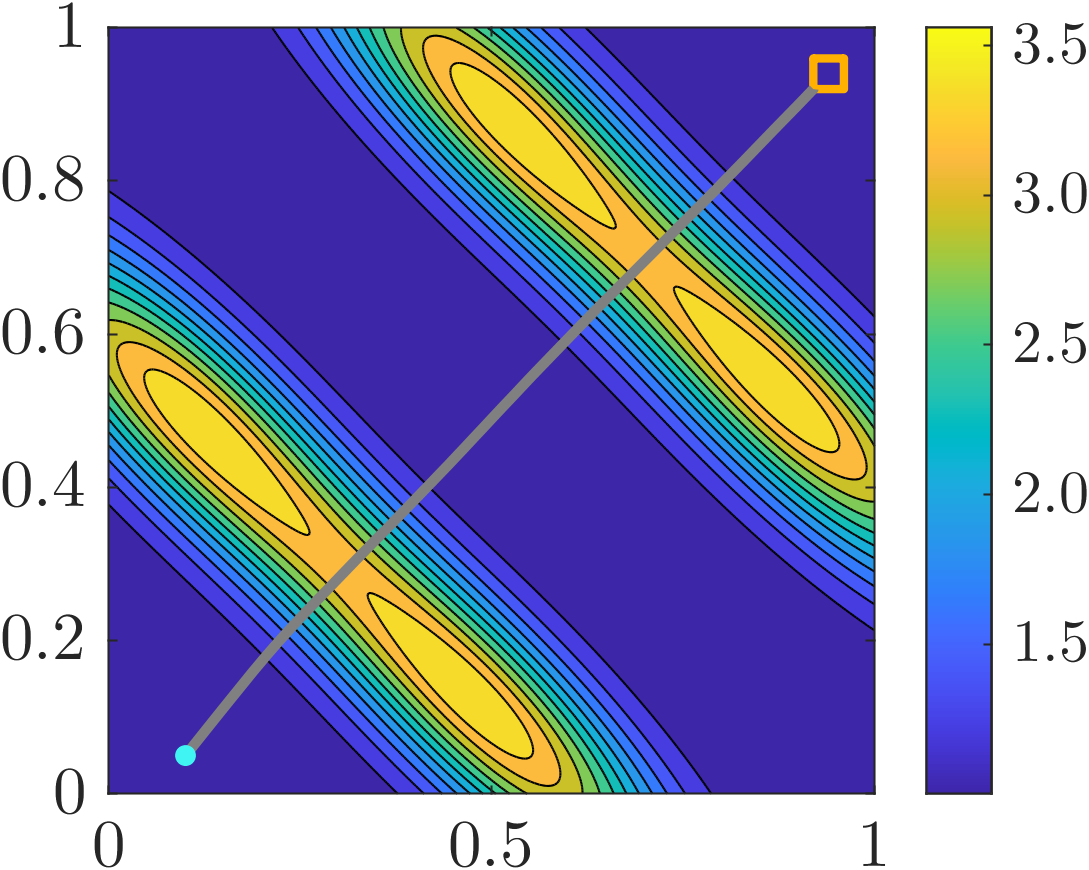}
\end{subfigure}
\caption{Surveillance patterns that form ``barriers'' along the direct path to the target (outlined in orange). 
Optimal trajectories are shown for $\lambda_{ij} = 0$ (no mode switches) and trajectory color encodes last observed mode (distribution).
White represents Mode 1 ($\vb b(0) = \vb e_1$), black Mode 2 ($\vb b(0) = \vb e_2$), and gray the stationary distribution ($\vb b(0) = [1/2, 1/2]$). 
Cyan dots represent the starting location.
 \update{Using the upper bound in Prop. \ref{prop:bounded-TG}, we solve the PDE over the time domain $[0, 14.83]$.
 } 
}
\label{fig:barriers-cost}
\end{figure}

We next investigate a simple indefinite horizon problem (Section \ref{sec:indef}), where the evader faces two possible surveillance patterns, each with two regions of elevated surveillance (``barriers'') on the way to the target (Figure \ref{fig:barriers-cost}).
Each barrier has the form  
$1 + \left(2\pi\sqrt{|\Sigma|}\right)^{-1} \exp(-2^{-1}(\vx - \hat{\vx}_i)\Sigma^{-1}(\vx - \hat{\vx}_i)^{\top})$ for $\Sigma = [3, -2.5; -2.5, 3]$.
We assume that the transition rates are symmetric ($\lambda_{12} = \lambda_{21}$) in all cases, leading to the stationary mode distribution $\vb q_s = \left[1/2, 1/2\right]$.
For a fixed initial location $\vy(0) = (0.1,0.05)$ with no mode switches and with a known initial mode $\vb b(0) = \vb e_i,$  
the optimal path to the target is a serpentine trajectory that avoids both high-surveillance areas.
This is not the case when there is a high degree of uncertainty about the current mode (e.g., $\vb b(t) = \vb q_s$).
The planner no longer has enough information to commit to either serpentine trajectory, and it instead becomes optimal to travel through the center of the domain.

\begin{figure}[h]
\centering
\begin{subfigure}{0.29\textwidth}
\centering
\caption{No observations, $\lambda_{ij} = 0.5$}
\label{fig:barriers-id-lambda05}
\includegraphics[height=4.5cm, trim={20 15 50pt 0}, clip]{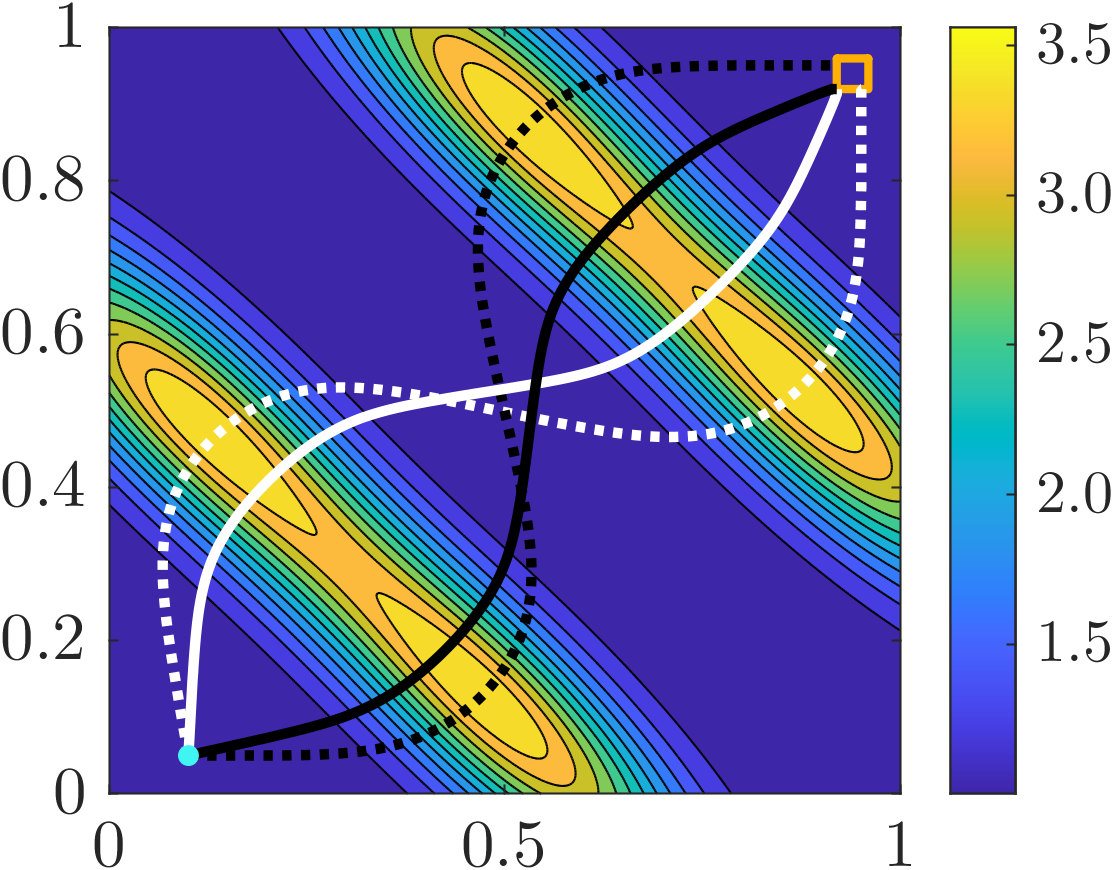}
\end{subfigure}
\begin{subfigure}{0.31\textwidth}
\centering
\caption{No observations, $\lambda_{ij} = 1.0$}
\label{fig:barriers-id}
\includegraphics[height=4.5cm, trim={20 15 50pt 0}, clip]{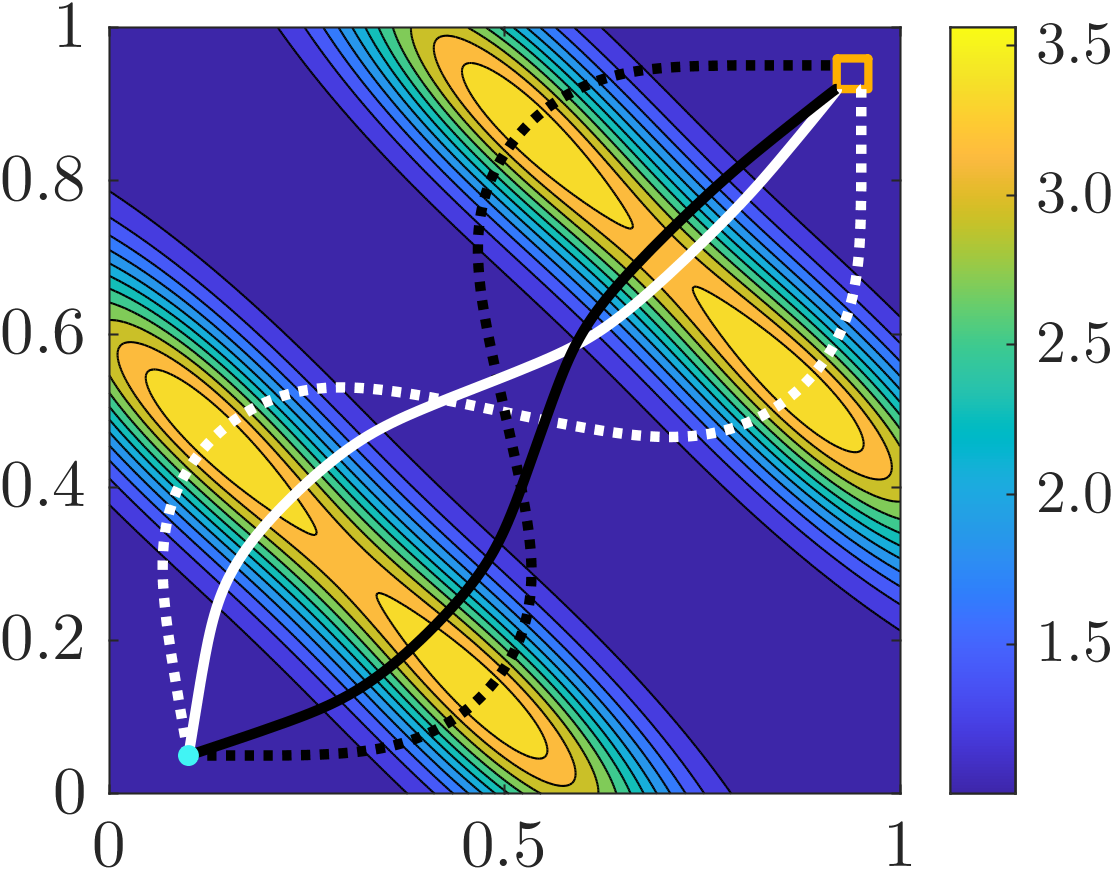}
\end{subfigure}
\begin{subfigure}{0.34\textwidth}
\centering
\caption{One observation, $\lambda_{ij} = 1.0$}
\label{fig:barriers-1obs}
\includegraphics[height=4.5cm, trim={22 15 0 0}, clip]{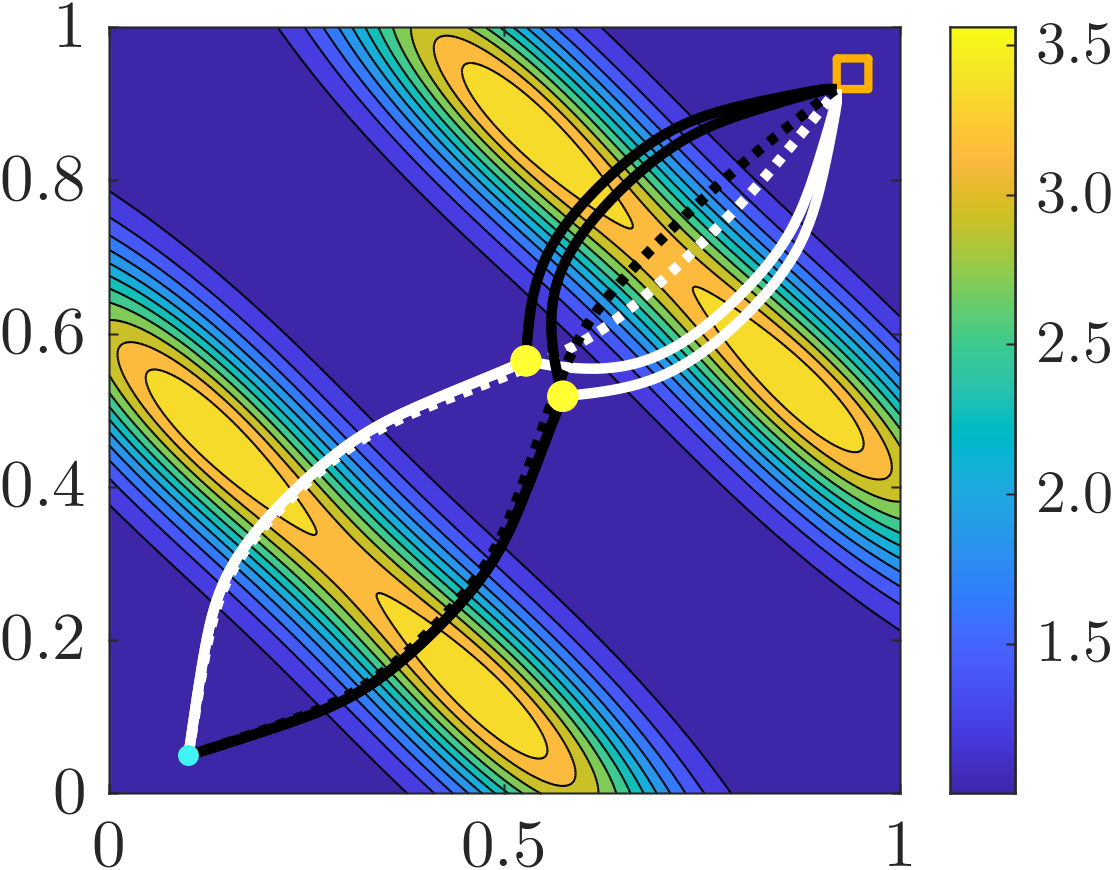}
\end{subfigure}
\caption{Impact of mode transition rates and availability of mode observations on optimal trajectories. Across all figures, optimal trajectories are shown by solid lines, 
with the color encoding the last observed mode as in Figure~\ref{fig:barriers-cost}. 
In all cases, the background is $\overline{K}_s(\vx)$.
Cyan dots represent the starting location and yellow dots represent mode observations.
In \textbf{(a)} and \textbf{(b)}, the planner only has access to the initial mode distribution and trajectories corresponding to $\lambda_{ij} = 0$ are shown (dotted lines) for reference. \textbf{(c)} The planner has access to one free on-demand observation and takes detours depending on what is observed. No-observation trajectories are also shown for reference (dotted lines). 
}
\label{fig:barriers-traj}
\end{figure}

To highlight the impact of mode transition rates on optimal trajectories, we consider two examples: ``slow'' transition rates $\lambda_{ij} = 0.5$ and ``fast'' transition rates $\lambda_{ij} = 1$.
Figure \ref{fig:barriers-traj}(a,b) provides a comparison of these two cases for both $\vb b(0) = \vb e_1$ and $\vb b(0) = \vb e_2,$ assuming no additional mode observations are available.  Even with the slow transition rates, we observe that the detours around the likely surveillance barriers are smaller than we saw in Figure \ref{fig:barriers-cost} and the second detour is much smaller than the first.  This is due to the planner's decreasing confidence about the true value of $\mu(t)$.
This feature is even more noticeable with fast transition rates, where by the time of the second barrier, the planner hardly takes any detour at all (since by then $\vb b(t)$ is close to $\bq_s$).  However, the trajectories change if that planner is allowed to request a single (free) observation, the setting described in Section \ref{sec:indef-bdd}.
Figure \ref{fig:barriers-traj}(c) shows that they take an opportunity to learn the mode just before the second barrier, with an appropriate detour used immediately after that (depending on what was observed).

\subsection{Surveillance with Obstacles}
\label{sec:num-exp-maze}
\begin{figure}
\centering
\begin{subfigure}{0.35\textwidth}
\centering
\caption{Mode Labels}
\label{fig:maze-labels}
\includegraphics[height=5cm, trim={25 15 0 0}, clip]{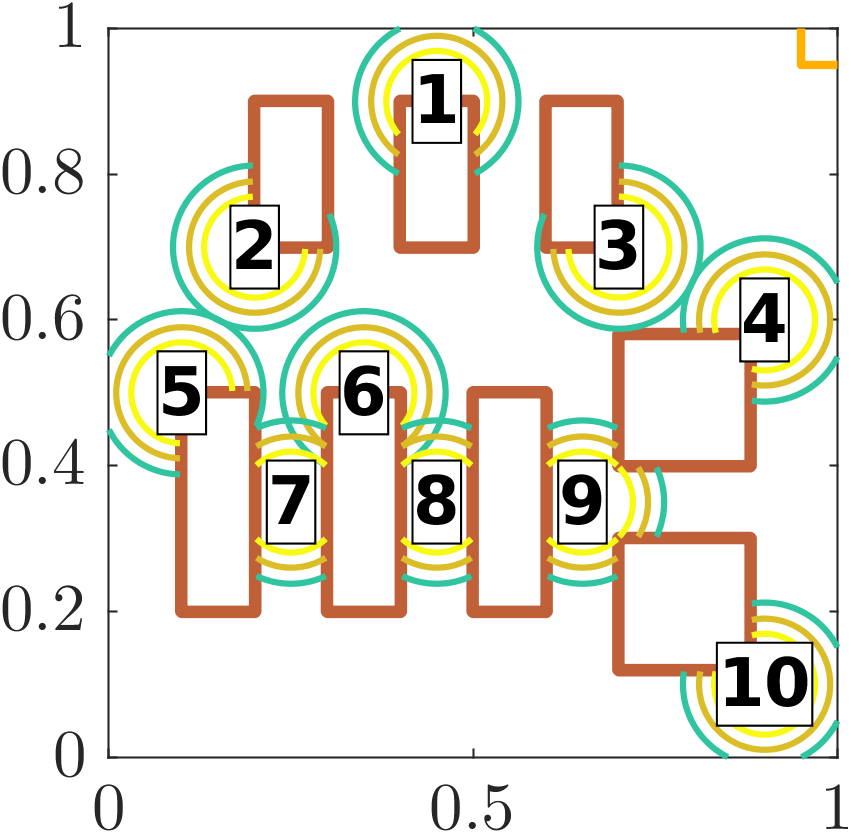}
\end{subfigure}
\hspace{0.5cm}
\begin{subfigure}{0.35\textwidth}
\centering
\caption{$\overline{K}_s(\vx)$}
\label{fig:maze-stationary-cost}
\includegraphics[height=5cm, trim={25 15 0 0}, clip]{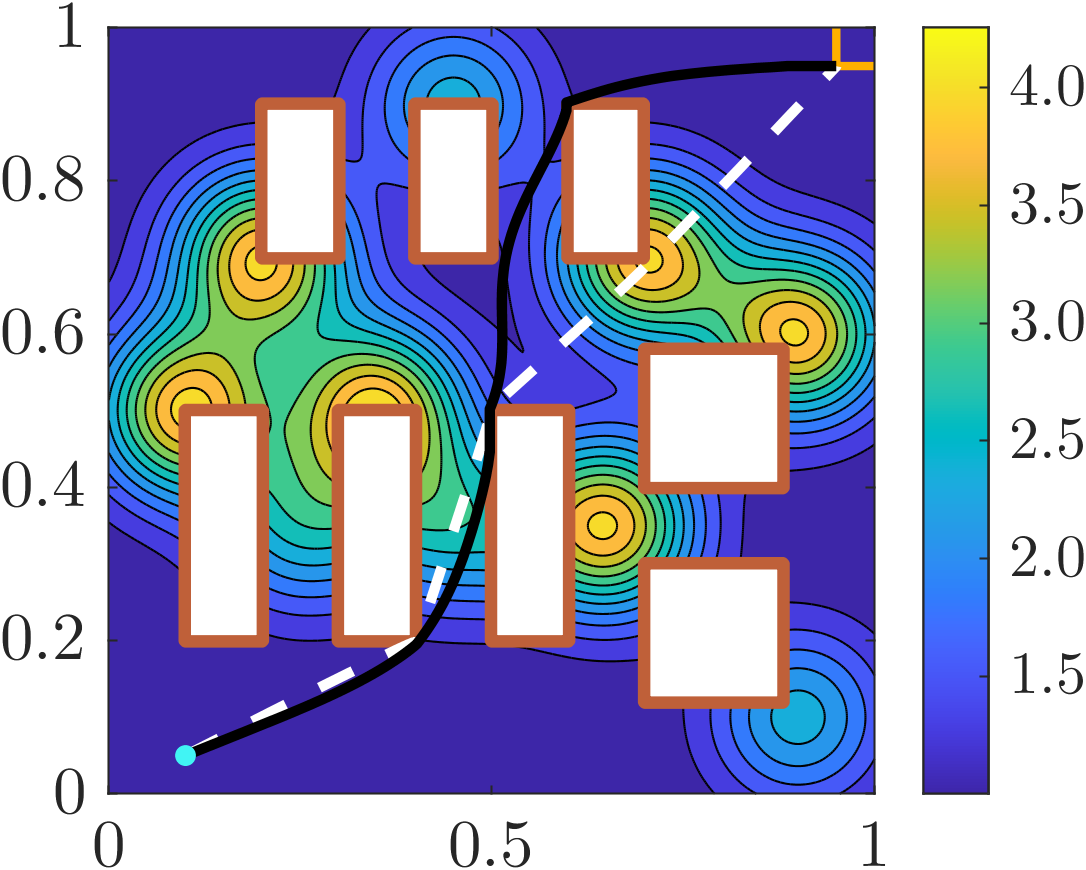}
\end{subfigure}
\caption{Surveillance patterns, obstacles, and stationary expected surveillance. Impassable obstacles are outlined in brown and the target is outlined in orange. \textbf{(a)} Labeled contour plots of mode surveillance patterns $K_i(\vx)$ (formatted as in Figure \ref{fig:4modes-labels}), which guard paths between and around obstacles. 
\textbf{(b)} The expected surveillance associated with the stationary mode distribution. Two trajectories are shown: a $\overline{K}_s$-optimal trajectory (black solid) and a time-optimal trajectory (white dashed), both with no mode observations.
\update{
If the planner begins with a belief equal to the stationary distribution, the problem becomes time-independent and both trajectories above could be computed in less than a second by solving a standard Eikonal equation. 
As an additional verification, the black curve was instead based on $v_{\vb q_s}$ computed via Alg. \ref{alg:indef-bdd-solver} with $L = 0$. The discretized equations were solved over the time interval $[0, 37.68].$
}}
\label{fig:maze-cost}
\end{figure}

This is another indefinite horizon example, in which
the planner must navigate a maze-like environment with impassable obstacles and ten possible surveillance patterns (Figure \ref{fig:maze-cost}). 
(That many modes would make  true belief-space dynamic programming computationally infeasible.)
We again take $f(\vx) = 1$ outside of obstacles 
and compute the value functions over a $401\times401$ grid.
Within Mode $i$, the surveillance intensity is given by $K_i(\vx) = 1 + 12(2\pi\sigma)^{-1} \exp(-(2\sigma^{2})^{-1}(\vx - \hat{\vx}_i)(\vx - \hat{\vx}_i)^{\top})$ for $\sigma = 0.08$ and with centers $\hat{\vx}_i$ shown in Figure \ref{fig:maze-labels}.
The adversary uses a more complicated matrix of switch rates to move through their surveillance patterns randomly:
\begin{equation}
\Lambda = \begin{bmatrix}
-2 &  1 &  1 &  0 &  0 &  0 &  0 &  0 &  0 &  0 \\
1/2& -2 & 1/2&  0 &  1 &  0 &  0 &  0 &  0 &  0 \\
1/2& 1/2& -2 &  1 &  0 &  0 &  0 &  0 &  0 &  0 \\
 0 &  0 &  1 & -2 &  0 &  0 &  0 &  0 & 1/2& 1/2\\
 0 &  1 &  0 &  0 & -2 & 1/2& 1/2&  0 &  0 &  0 \\
 0 &  0 &  0 &  0 & 1/2& -2 &  0 & 1/2&  1 &  0 \\
 0 &  0 &  0 &  0 &  1 &  0 & -2 &  1 &  0 &  0 \\
 0 &  0 &  0 &  0 &  0 &  1 &  1 & -2 &  0 &  0 \\
 0 &  0 &  0 & 1/2&  0 &  1 &  0 &  0 & -2 & 1/2\\
 0 &  0 &  0 &  1 &  0 &  0 &  0 &  0 &  1 & -2
\end{bmatrix},
\end{equation}
which has stationary distribution
$
\bv q_s = \begin{bmatrix}
\frac{1}{16}, & \frac{1}{8}, & \frac{1}{8}, & \frac{1}{8}, & \frac{1}{8}, & \frac{1}{8}, & \frac{1}{16}, & \frac{1}{16}, & \frac{1}{8}, & \frac{1}{16}
\end{bmatrix}.
$
Throughout this example, we assume that $\vb b(0) = \vb q_s$, which we interpret as the planner having no information about the initial mode beyond what can be gained from knowing $\Lambda$.

The expected stationary surveillance intensity $\overline{K}_s(\vx)$ is shown in Figure \ref{fig:maze-stationary-cost}.
High surveillance areas guard paths between and around obstacles, leading to ``decision points'' where the planner must commit to a strategy for 
dealing with
a particular obstacle (e.g., choosing to go clockwise or counterclockwise around it).
In the absence of further mode information, the optimal trajectory is similar to the time-optimal path, but in the final stretch it takes a slight detour through a gap near $\hat{\vx}_1,$ which is less likely to be surveilled.

\begin{figure}
\centering
\begin{subfigure}{0.48\textwidth}
\centering
\caption{Slow mode switches}
\includegraphics[height=3.6cm]{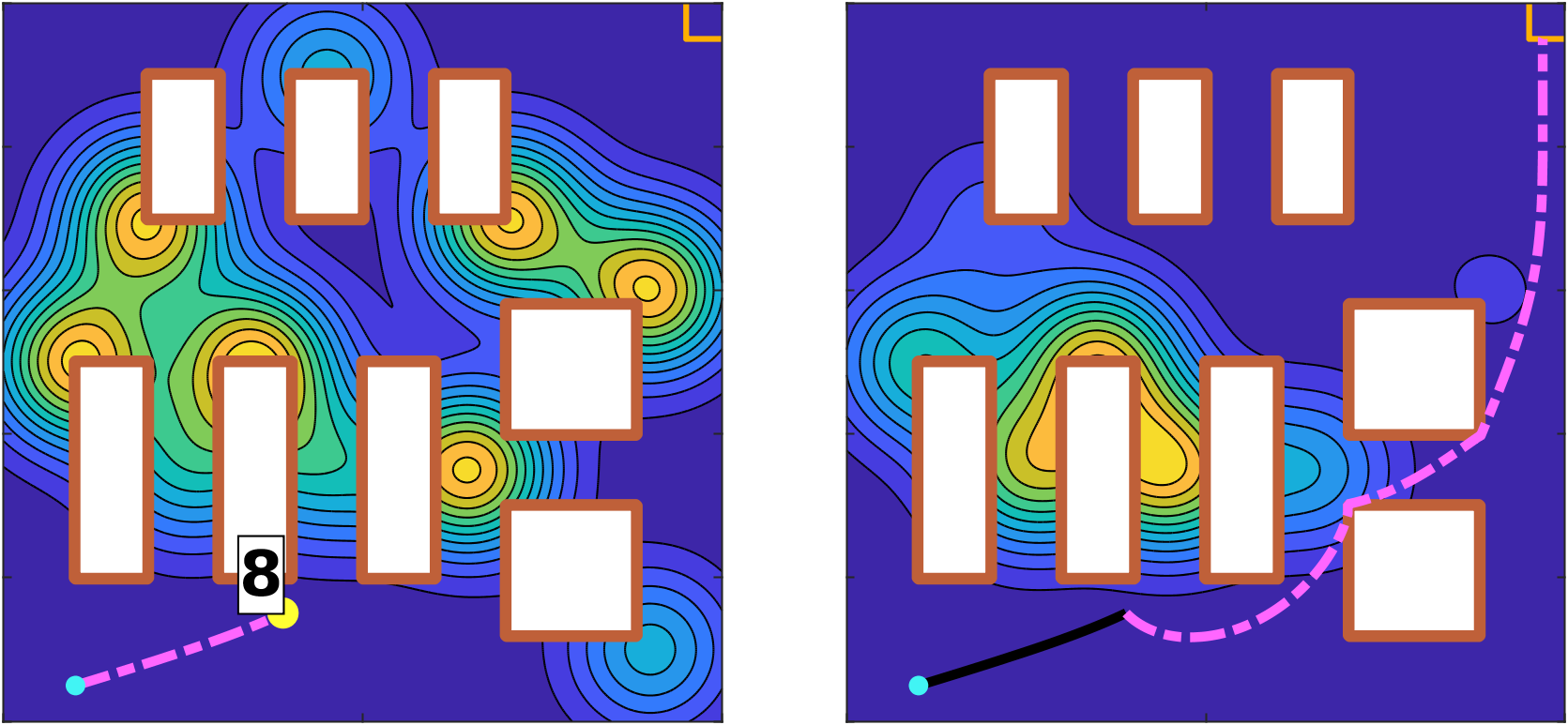}
\end{subfigure}
\hfill
\begin{subfigure}{0.48\textwidth}
\centering
\caption{Fast mode switches}
\includegraphics[height=3.6cm]{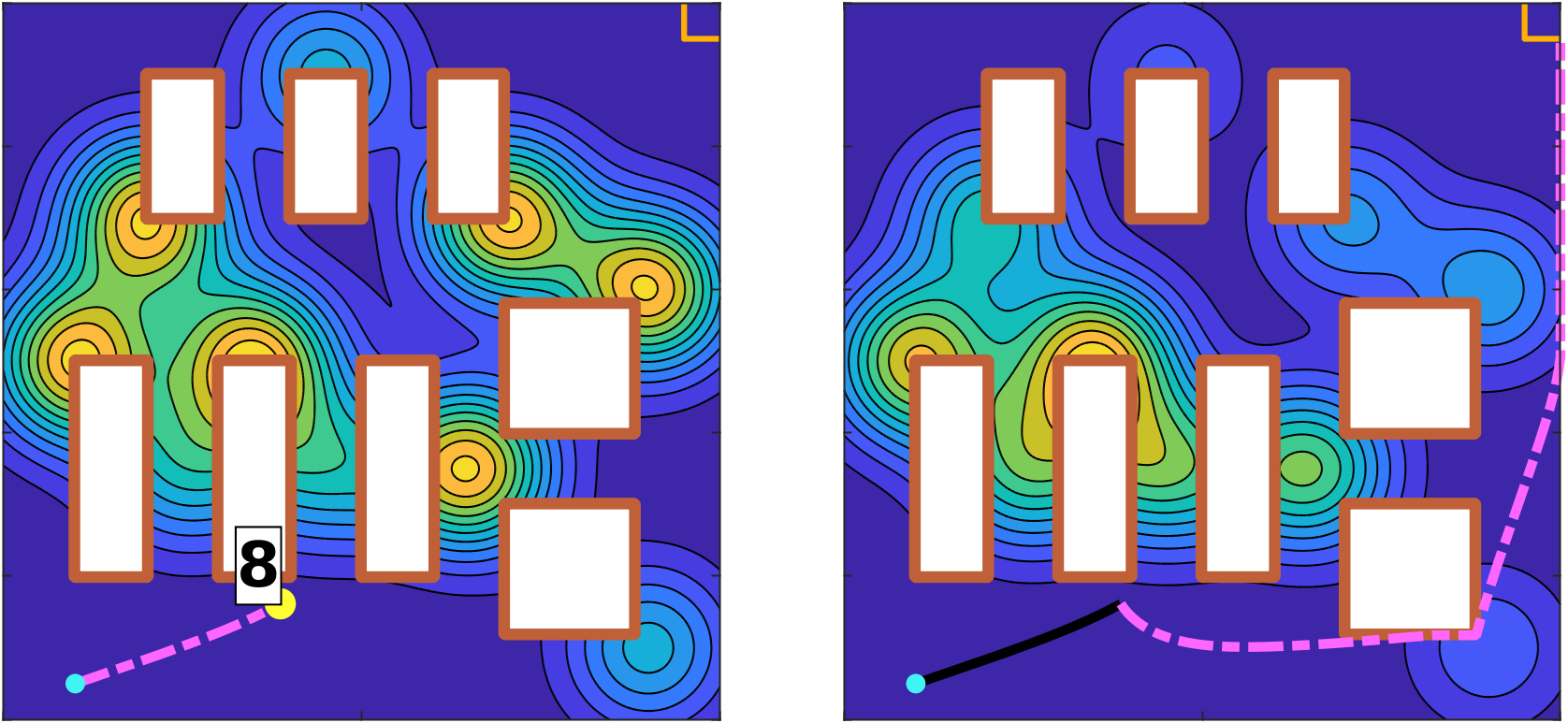}
\end{subfigure}
\caption{Impact of magnitude of $\Lambda$ on optimal trajectories. The planner begins with $\vb b(0) = \vb q_s$ and has access to one free on-demand observation. Trajectories are formatted as in Figure \ref{fig:4modes-fin-traj-1}. The background is the expected surveillance at the end of the current (magenta dot-dash) segment of the trajectory (just before an observation if applicable).
\textbf{(a)} Transition matrix: $\Lambda$.
\textbf{(b)} Transition matrix: $2\Lambda$. 
}
\label{fig:maze-traj}
\end{figure}

We will again consider ``slow'' and ``fast'' mode switches: slow corresponding to the rate matrix $\Lambda$ above, and fast corresponding to $2\Lambda$.
Figure \ref{fig:maze-traj} shows a realized optimal trajectory for both switching regimes when the planner has access to one free on-demand observation\footnote{
As pointed out in Remark~\ref{rem:stationary_q}, the initial belief in this example is rather special.  Without observations, $\vb b(0) = \vb q_s$ implies $\vb b(t) = \vb q_s$ for all $t,$ and so the optimal trajectory in Figure \ref{fig:maze-stationary-cost} can be obtained much more cheaply by solving PDE $|\nabla u | f(\x) = \overline{K}_s(\x)$ with one of the standard fast methods developed for Eikonal equations; e.g., \cite{sethian1996fast, tsai2003sweeping, ChacVlad1, potter2019ordered}.  
But the fact that $\vb b(t)$ starts to vary in time immediately after any mode observation makes it necessary to use the more expensive computational framework developed here for OOPDMP problems.
} (Section \ref{sec:indef-bdd}).
The optimal observation location is similar across both switching regimes: the planner requests an observation just before entering the corridor that would be optimal to travel through according to 
$\overline{K}_s$.
When the planner observes that this ``assumed to be optimal'' corridor is heavily guarded (Mode 8), it instead takes a detour, the length of which depends on the switching rate.
When the switching rate is low, the planner takes a shorter detour past $\hat{\vx}_9$, since it is unlikely the adversary will have switched to guarding that area by then.
However, with high switching rates,
that path becomes too risky (our mode observation \update{loses} relevance sooner), and the planner chooses to take a longer (but safer) path.

\begin{figure}
\centering
\begin{subfigure}{0.48\textwidth}
\centering
\caption{Slow mode switches}
\label{fig:F1-1obs-all}
\hspace{-0.15cm} 1 Mode observation \hspace{0.35cm} 2 Mode observations\\
\includegraphics[height=3.6cm,trim={21 17 40 0},clip]{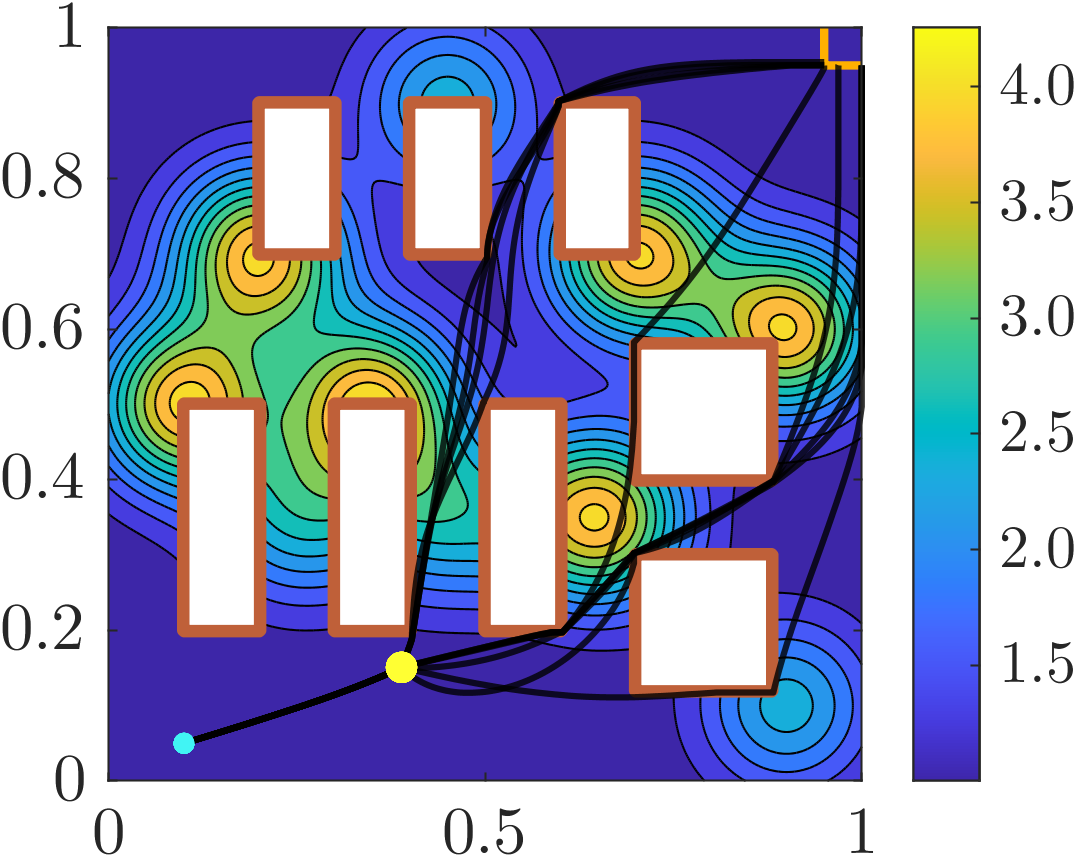}
\hspace{0.05cm}
\includegraphics[height=3.6cm,trim={21 17 40 0},clip]{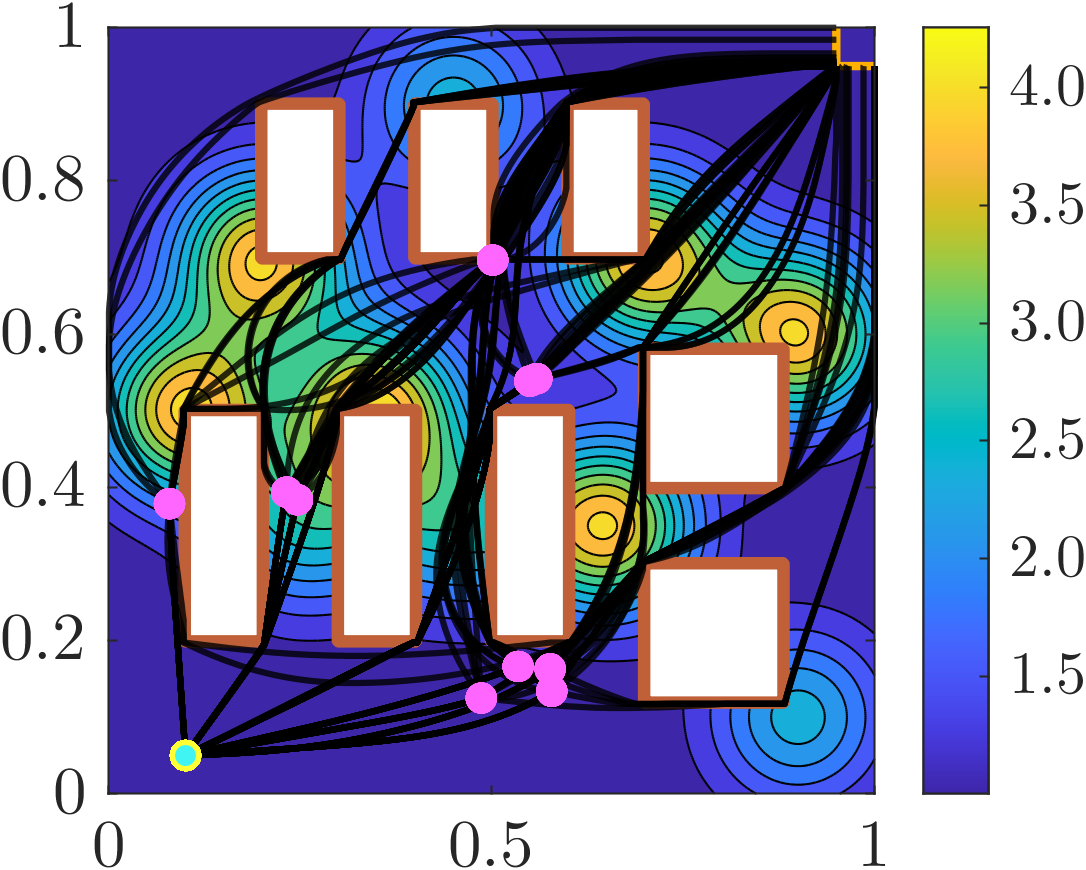}
\end{subfigure}
\hfill
\begin{subfigure}{0.48\textwidth}
\centering
\caption{Fast mode switches}
\label{fig:F2-1obs-all}
\hspace{-0.15cm} 1 Mode observation \hspace{0.35cm} 2 Mode observations\\
\includegraphics[height=3.6cm,trim={21 17 40 0},clip]{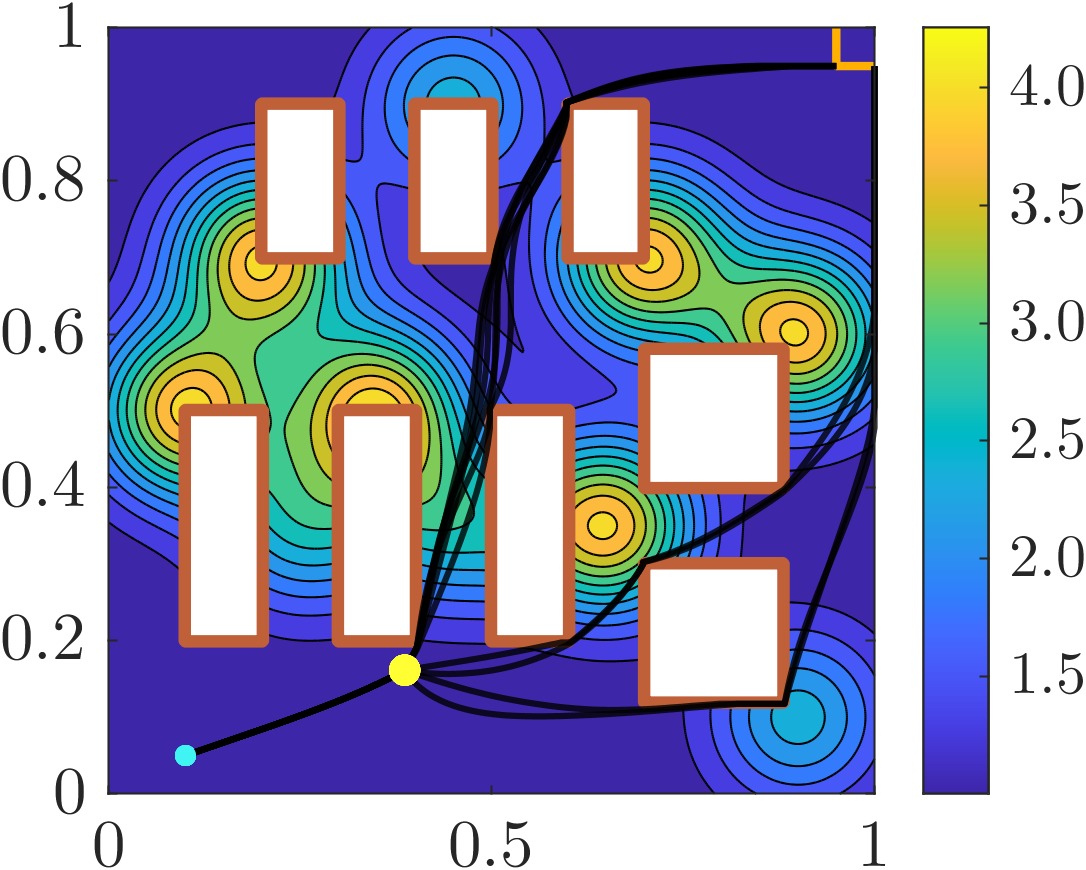}
\hspace{0.05cm}
\includegraphics[height=3.6cm,trim={21 17 40 0},clip]{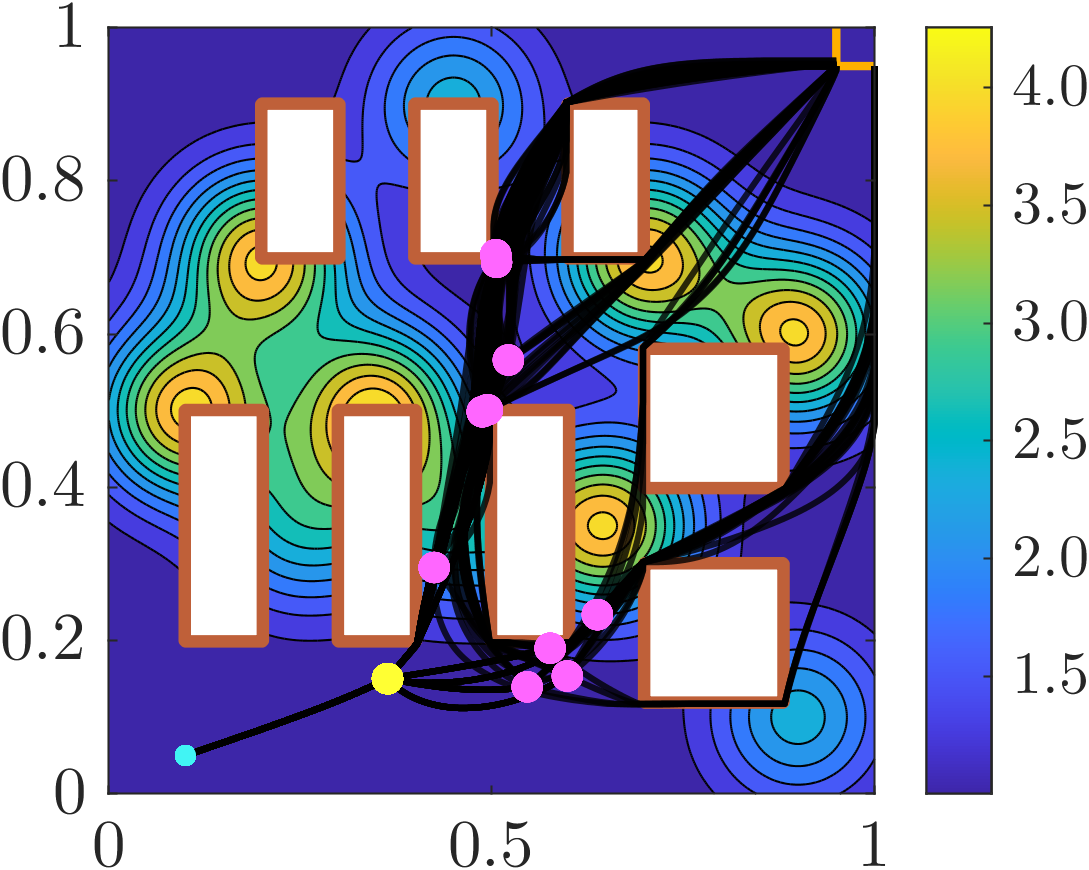}
\end{subfigure}
\caption{Impact of $\Lambda$ on optimal observation locations and resulting trajectories. The initial position (cyan dot) is fixed, observations are free and on-demand, with a limit of one (first and third images) or two (second and fourth images).
Optimal observation locations change as the rate of mode switching increases. Yellow dots encode the first observation location, and magenta dots encode the second (if applicable). Optimal trajectories (black solid lines) are shown for all possible observation sequences and $\vb b(0) = \vb q_s$. Background is $\overline{K}_s(\vx)$.
\textbf{(a)} Transition matrix: $\Lambda$
\textbf{(b)} Transition matrix: $2\Lambda$.
}
\label{fig:maze-all-paths}
\end{figure}

To investigate the impact of $\Lambda$ more thoroughly, Figure \ref{fig:maze-all-paths} shows {\em all} optimal trajectories for both switching regimes, and for one or two free on-demand observations.
(I.e., we show all trajectories that might be realized depending on the actually observed modes.)
When only one observation is available, the overall behavior in the ``slow'' and ``fast'' switching regimes is quite similar (even if examples exist where responses to the same observed mode differ significantly, as in Figure \ref{fig:maze-traj}).
But when the planner is allowed to request two observations, the magnitude of switching rates impacts optimal observation locations, in addition to the shape of the trajectories. 
When the mode switching is slow, the planner now requests the first observation almost immediately, and then uses the information to determine how to approach the maze (e.g., go around the outside vs. pass through the center).
This significantly changes the qualitative behavior of optimal trajectories, as the planner can commit to longer detours to avoid areas where $K_{\mu(t)}$ is expected to be high.
When the mode switching is fast, it instead remains optimal to ``save'' the first observation until the planner is closer to the $\overline{K}_s$-optimal decision point.
Mode switches are expected to occur rapidly enough that information gained from observations does not stay relevant long enough for the planner to commit to long detours.
Thus, the magnitude of the switching rates decreases the longevity of mode information and incentivizes making observations just before the information is needed.

\subsection{Mars Rover}
\label{sec:rover}

\begin{figure}
\centering
\begin{subfigure}{0.27\textwidth}
\centering
\caption{$f(\vx)$}\label{fig:mars-speed}
\includegraphics[height=4.3cm, trim={0 0 42pt 0}, clip]{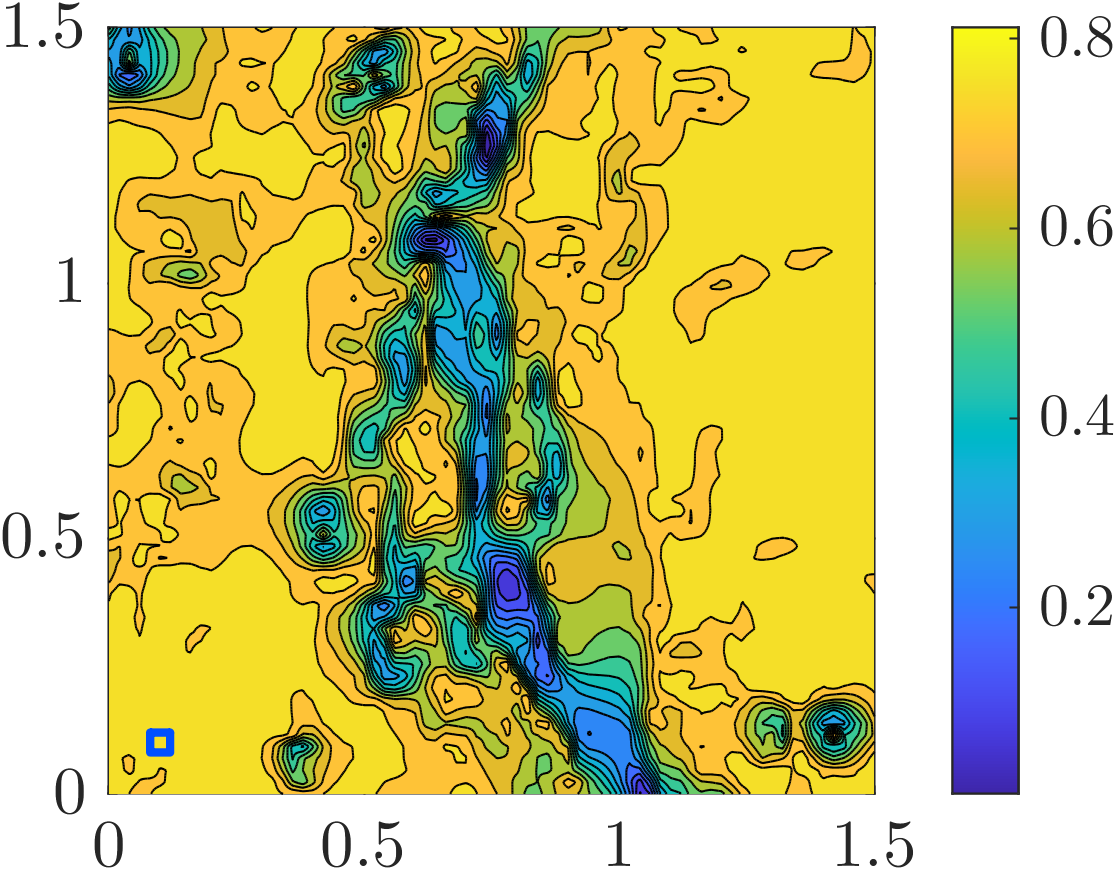}
\end{subfigure}
\begin{subfigure}{0.31\textwidth}
\centering
\caption{$f_b(\vx)$}\label{fig:mars-speed-broken}
\includegraphics[height=4.3cm]{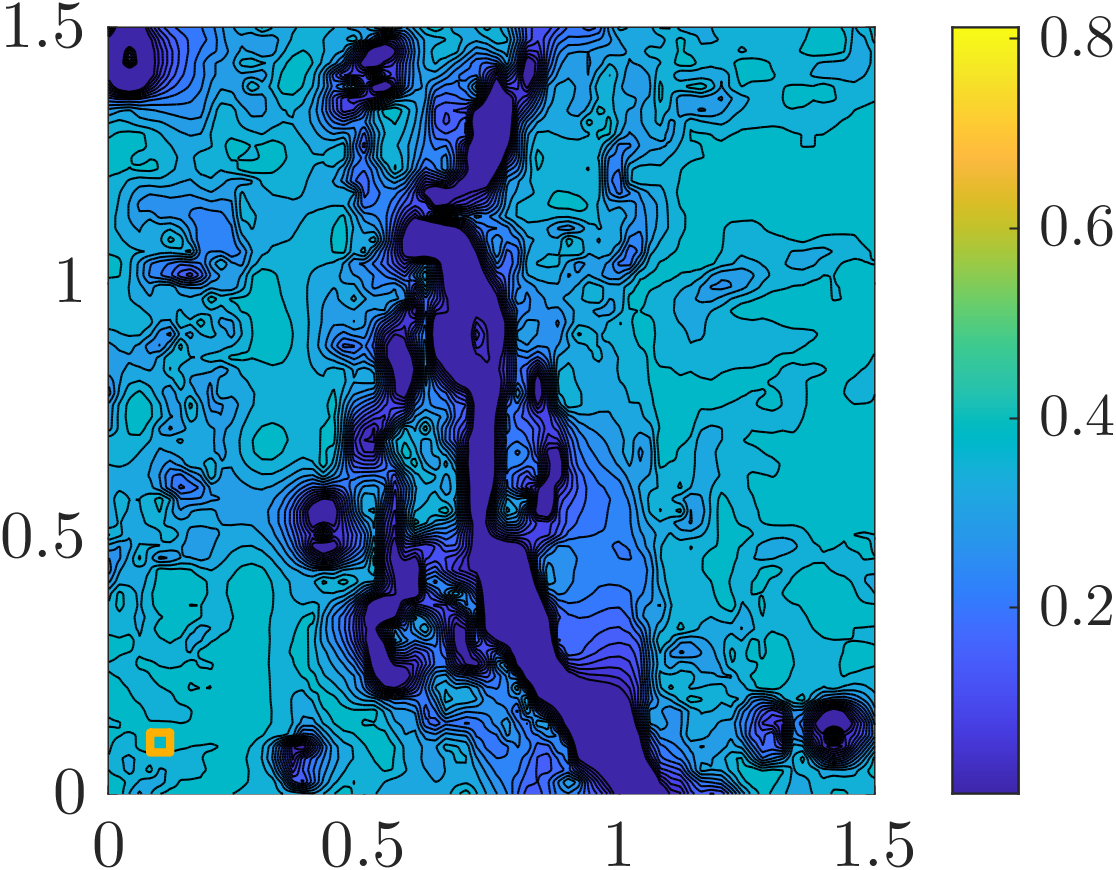}
\end{subfigure}
\hfill
\begin{subfigure}{0.40\textwidth}
\centering
\caption{$\phi(\vx)$}\label{fig:mars-terminal-cost}
\includegraphics[height=4.3cm]{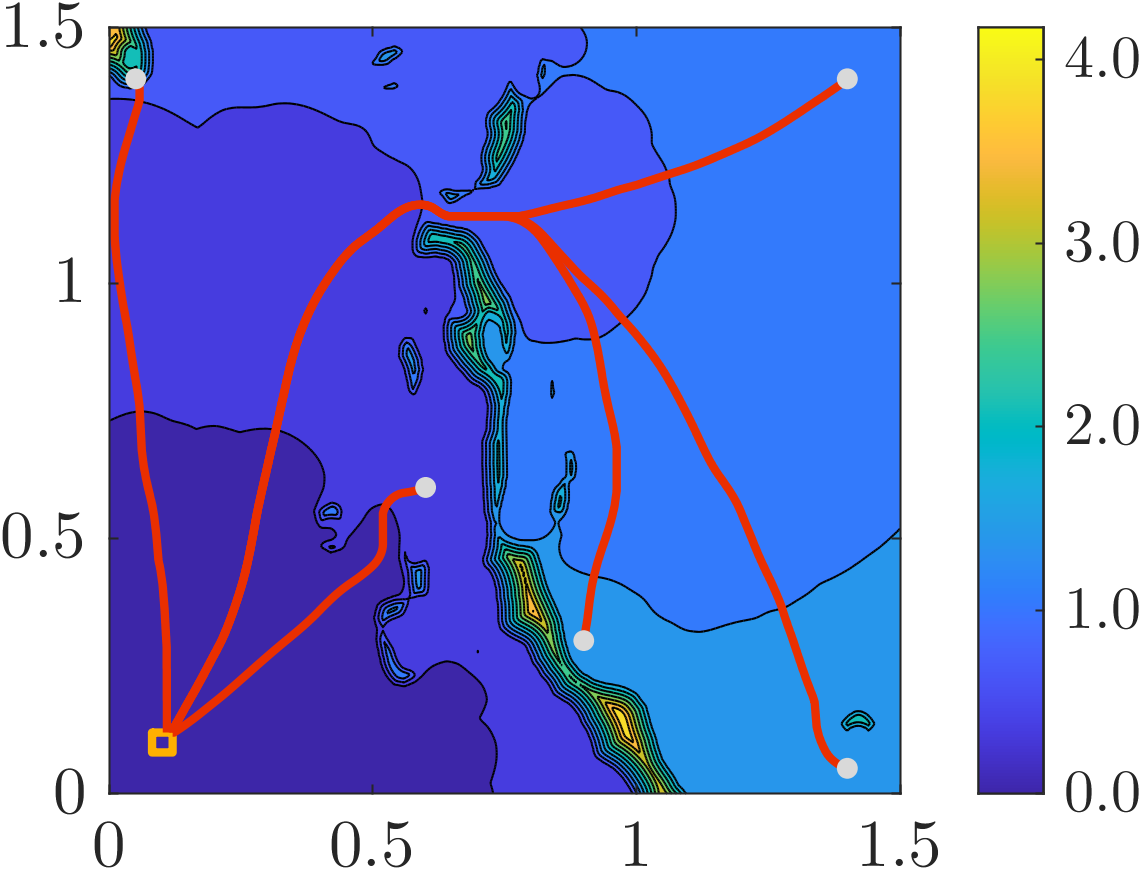}
\end{subfigure}
\caption{Speed and terminal cost for a rover that may experience observable and unobservable breakdowns. \textbf{(a)} The speed in Modes 1 and 2, before an observable breakdown has occurred. Target is shown in blue. \textbf{(b)} The speed after an observable breakdown occurs. Target is shown in orange. \textbf{(c)} The optimal time-to-target after an observable breakdown occurs. Optimal post-breakdown trajectories are shown in red for a variety of starting locations (gray dots) and the target is shown in orange. 
}
\label{fig:mars-env}
\end{figure}

The final example that we consider is motivated by optimal path planning for a Mars rover that may become damaged as it navigates hazardous terrain. 
We assume that the rover (located within the area of Jezero crater on Mars) seeks to reach a target $\Gamma$ while minimizing its expected travel time ($K_i(\vx) = 1$, $\psi_i(\vx) = 0$).
A previous PDMP model of this process assumed that ``damage'' takes the form of breakdowns that are always observable \cite{gee2022breakdowns}.
In reality, incremental damage may not immediately impact the observable dynamics of the rover.
Here, we consider a rover that may accumulate \emph{unobservable} damage that increases the chance of a future observable breakdown.
Before a breakdown has occurred, the rover can be in one of two possible modes: fully functional (Mode 1) or incrementally damaged / breakdown-prone (Mode 2), both with the same operating speed $f(\vx)$ (Figure \ref{fig:mars-speed}, units: $m/$sol, where a sol is a Martian day).
The terrain-dependent speed is computed as in \cite{gee2022breakdowns}; data was acquired using JMARS, a Mars GIS \cite{christensen2009jmars}, and speed was scaled down by a factor of $1/2$.

This process is an example of a prematurely ``terminated'' problem (Section \ref{sec:rt}), where the unobserved process ``terminates'' when an observable breakdown occurs.
\update{This is not a termination in the usual sense of the word since the rover still needs to keep moving, though now with a much lower speed  $f_b \ll f$ (Figure \ref{fig:mars-speed-broken}).}
The terminal cost $\phi(\vx)$ encodes the minimal time needed to reach the target after that breakdown from a point $\x$.
It is computed by solving $f_b|\nabla \phi| = 1$ on $\Omega \backslash \Gamma$ with the boundary condition $\phi = 0$ on $\partial \Gamma$ using the Fast Marching Method \cite{sethian1996fast}. 
Figure \ref{fig:mars-terminal-cost} shows the level sets of $\phi(\vx)$ as well as the resulting post-breakdown optimal trajectories.
Unobserved mode switches impact the likelihood of observable breakdowns via the mode-dependent termination rates $\gamma_1 = 1$ and $\gamma_2 = 12.33$ (all rates haves units sol$^{-1}$). 
Transitions between unobservable modes are assumed to occur with rates $\lambda_{12} = 5 $ and $\lambda_{21} = 0;$ i.e., that incremental unobserved damage is never fixed. But we assume that it can be discovered: the rover can ``pay'' for on-demand mode observations (Section \ref{sec:indef-paid-for}) by staying in place for a fixed amount of time \update{(assumed here to be $C(\vx) = 0.0203$ sol or, equivalently, 30 minutes)} to run a \update{full} system diagnostic.

\begin{figure}
\centering
\begin{subfigure}{0.265\textwidth}
\centering
\caption{$\vb b(0) = \vb e_1$}\label{fig:mars-traj-mode1}
\includegraphics[height=4.5cm, trim={22 12 60pt 0}, clip]{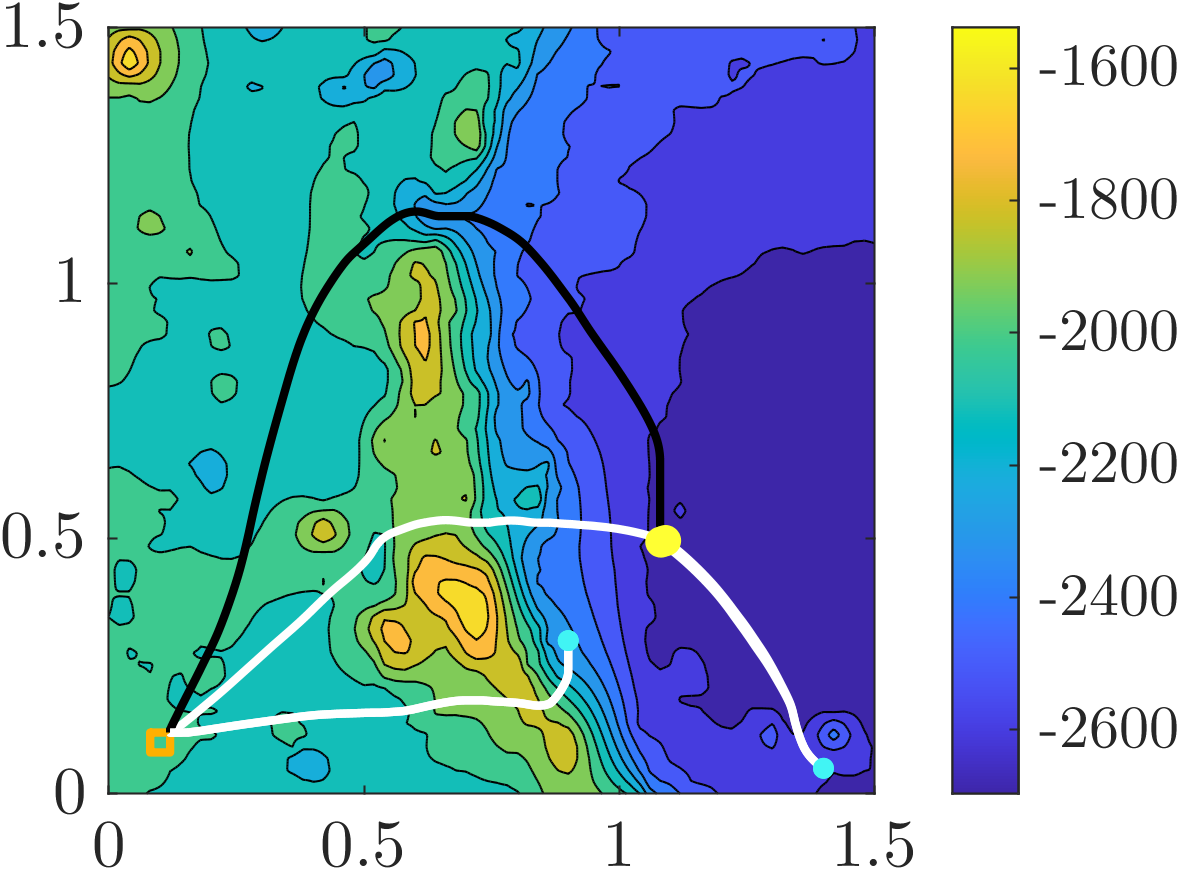}
\end{subfigure}
\begin{subfigure}{0.33\textwidth}
\centering
\caption{$\vb b(0) = \vb q_s$}\label{fig:mars-traj-qstat}
\includegraphics[height=4.5cm, trim={22 12 0 0}, clip]{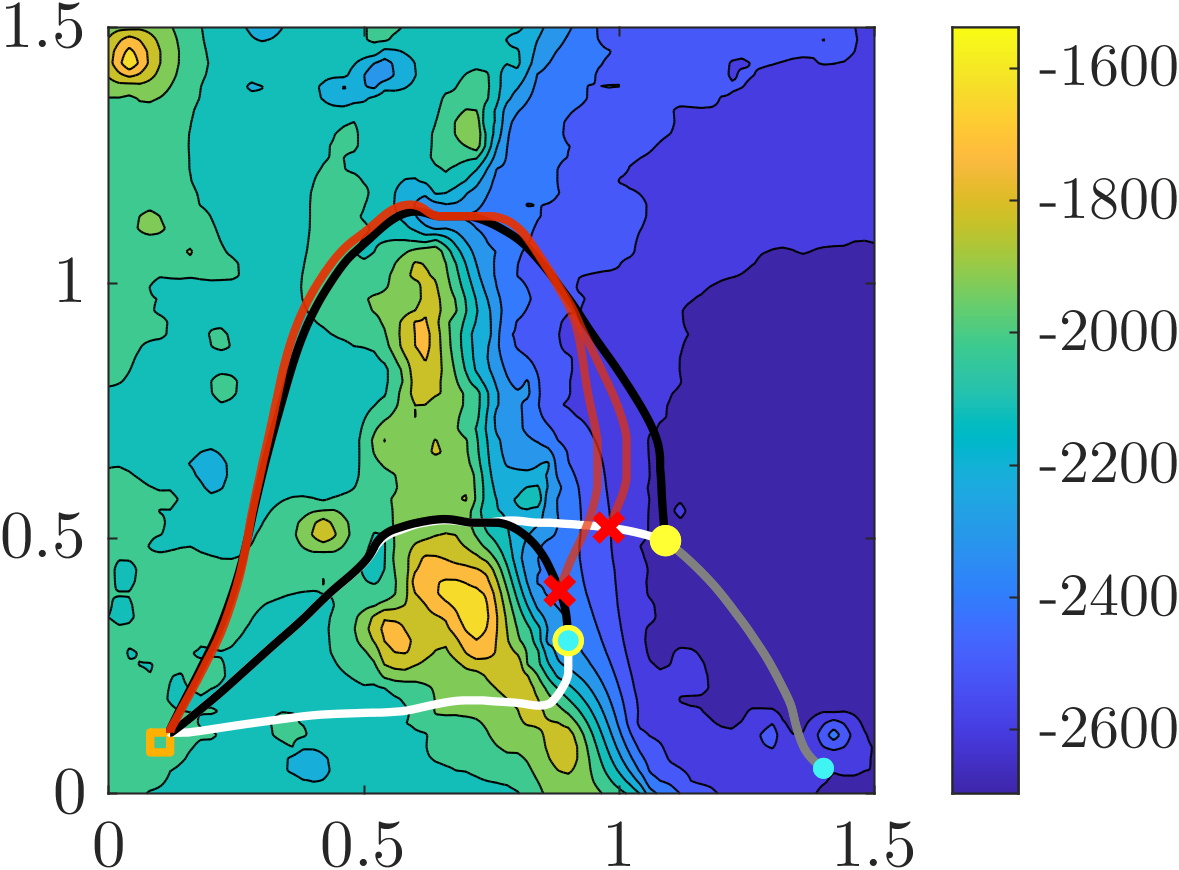}
\end{subfigure}
\hfill
\begin{subfigure}{0.35\textwidth}
\centering
\caption{Observation Region, $\vb b(t) = \vb q_s$}\label{fig:mars-obs}
\includegraphics[height=4.5cm, trim={22 12 0 0}, clip]{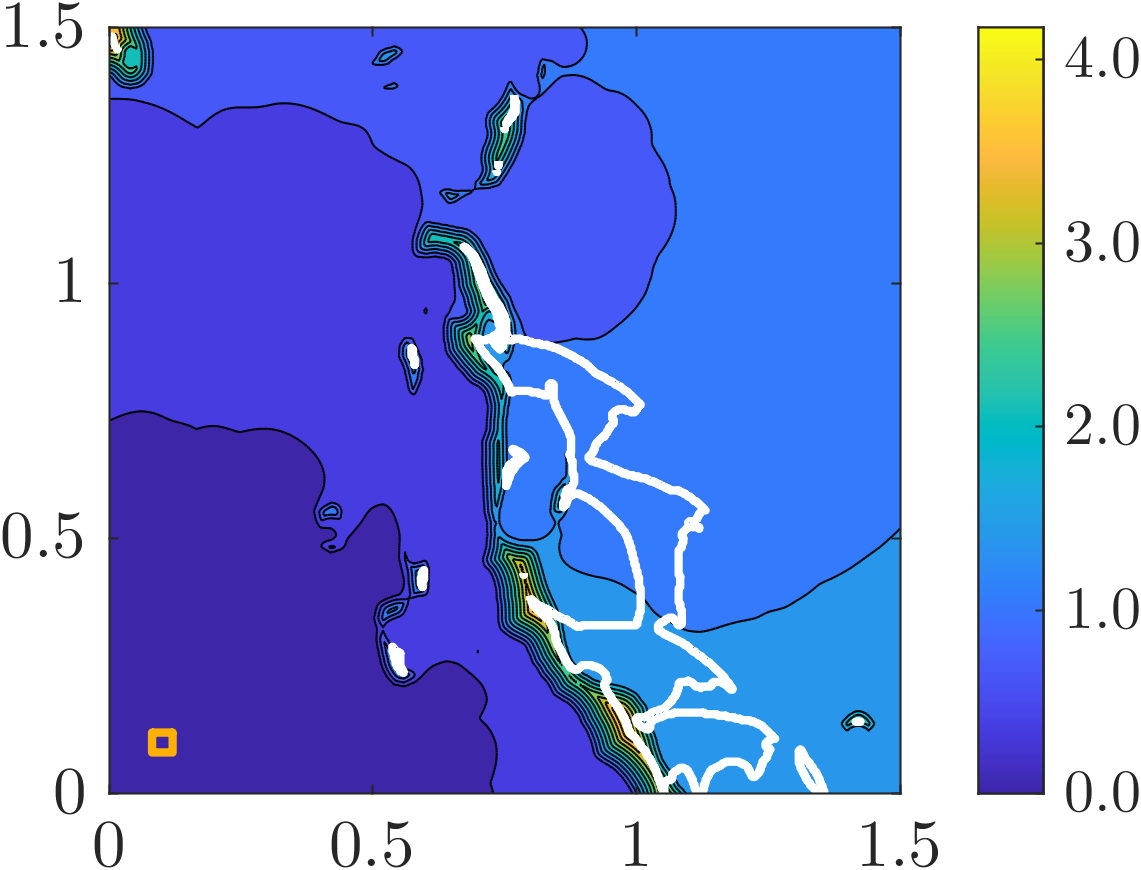}
\end{subfigure}
\caption{\textbf{(a)}-\textbf{(b)} Optimal trajectories and observation locations when the rover can observe its current mode by running diagnostics. In each panel, the background is the terrain    elevation, the target is outlined in orange, starting locations are cyan dots, observation locations are yellow dots, and trajectory color corresponds to last observed mode (distribution), as in Figure \ref{fig:barriers-traj}. Red ``x''s mark potential breakdown locations, with the resulting optimal trajectories shown in red. \textbf{(c)} Outline (white) of the region in which purchasing an observation is optimal. Background is $\phi(\vx)$. 
\update{The equation is solved numerically over the time interval $[0, 5.46]$. The computation requires 13 iterations 
to converge to within a tolerance of $10^{-6}$.}
}
\label{fig:mars-traj-obs}
\end{figure}

Figure \ref{fig:mars-traj-mode1} shows sample optimal trajectories and observation locations for two possible starting positions. 
If the rover is initially fully functional (Mode 1) and begins close to the crater rim, it is unlikely that incremental damage will occur in a short amount of time; thus, the rover does not bother with diagnostics before taking a relatively direct route over a steeply sloped region.
If the rover begins further away from the crater rim, there is a higher likelihood that it may accumulate unobserved damage by the time it reaches a direct path over the slope.
Thus, the rover will run diagnostics once it's closer to the rim to determine whether it is safe to cross directly.
If it observes Mode 1 (fully functional), a direct route is optimal.
Otherwise, the risk of a breakdown in the near future is high enough that the rover instead takes a long detour to avoid almost all steep slopes.
Here, the rover chooses to change its trajectory in anticipation of a breakdown, based on the knowledge that some damage has already occurred.

In Figure~\ref{fig:mars-traj-qstat}, we also consider a rover that begins in a ``stationary'' distribution.
In this case, the stationary distribution is not with respect to the mode-switching CTMC, but instead is the limit of equation \eqref{eq:b-rt} as $t$ becomes large.
Unless the robot starts with $\vb b(0) = \vb e_2,$ this limiting belief\footnote{
Both $\vb e_2$ and $\vb q_s$ are the equilibria of the of an ODE system satisfied by $\vb b(t)$ in the randomly terminated case. Appendix \ref{app:b-ode-bayes} presents a derivation of this ODE via Bayes theorem. It can be also derived using equations \eqref{eq:r-ode} and \eqref{eq:b-rt}.} 
will be  $\vb q_s \approx [0.5587, 0.4413]$.
As in Section \ref{sec:num-exp-maze}, setting $\vb b(0) = \vb q_s$ can be interpreted as not knowing the status of the rover beyond the fact that an observable breakdown has not yet occurred after a long operating time.
In this case, if the rover begins far from the rim of the crater, it still waits to run diagnostics until it is closer to the steeply sloped areas.
But if the rover begins close to the crater rim, it runs diagnostics immediately to determine its operating condition.
If diagnostics show that the rover is fully functional (Mode 1), it takes a path similar to the one in Figure \ref{fig:mars-traj-mode1}.
If the rover has already accumulated some damage (Mode 2), it instead takes a short detour through a more mildly sloped portion of the domain.
A much longer detour is taken in case of a real breakdown; optimal post-breakdown trajectories are shown in red for a couple of possible breakdown locations.
Figure \ref{fig:mars-obs} shows the region where it is optimal to run diagnostics\footnote{In general, the shape of this region will be time-dependent, since it is a function of the current belief. However, when $\vb b(0) = \vb q_s$, the belief (and thus the observation region) does not change in time.}
 when $\vb b(t) = \vb q_s$. In general, the diagnostics will be worthwhile when the rover is near or approaching steeply sloped regions where a possible breakdown would result in a high time penalty. 
 
Overall, this example demonstrates how unobservable incremental damage can significantly impact optimal trajectories even if it does not directly impact (yet) the process dynamics.

\section{Conclusion}
\label{sec:conclusion}
We have presented a new framework for defining and exploiting ``occasional observability'' in 
a subclass of piecewise-deterministic Markov processes (PDMPs) used to model abrupt changes in a global environment,
performance measures, or capabilities of a controlled system.  Such ``mode-switching'' PDMPs arise naturally
in many application areas and their special structure 
requires a different notion of observability.  Unlike in partially-observable general PDMPs, where each jump is immediately noted
and the full post-jump state is subject to noisy observations \cite{brandejsky2013stopping, bauerle2018popdmps}, 
in our mode-switching setting it is far more reasonable to assume that the continuous part of the system state is always fully observed 
while the mode-switches go unnoticed and the current mode is observed only occasionally (if ever).
Our Assumption \ref{a:no_mode_info_from_traj} ensures that the observed continuous state dynamics and the incurred running costs cannot be used to identify the current (unobserved) mode. 
It also allows for a simple handling of state-constraints (e.g., path planning on a domain with obstacles) since whether a control is allowable can be verified mode-independently,
and it is only the cumulative cost resulting from that control that is influenced by the random sequence of mode-switches. 
This formulation is a natural fit for many security applications (e.g., surveillance-evading path planning), models in behavioral ecology, and robotic navigation.

We have described a general method for dynamic programming on mode-belief space for such mode-switching PDMPs, but its practical usefulness is limited since the computational cost scales exponentially with the number of modes present in a problem.   We circumvent this curse of dimensionality by making an additional Assumption \ref{a:const_switching_rates}, 
which ensures that the current belief can be reconstructed from the initial belief and elapsed time.
This allows for a  lower-dimensional formulation of the control problem and much more efficient numerical methods (with computational cost scaling linearly with the number of modes).  We have presented the latter approach for a number of planning horizon/observation scheme combinations and illustrated it using several numerical experiments in surveillance avoidance and planetary rover navigation.  

An interesting future application would be to incorporate our approach  into {\em Stackelberg games}, which are frequently used in security domain.
E.g., in protecting natural resources, Authorities
might be deploying 
their limited resources to choose the best achievable pointwise-surveillance function $K$, 
and  Perpetrators (wildlife poachers or illegal forest loggers) are then planning their trajectories with that $K$ in mind \cite{CarteeVlad_Poaching}.  
The methods presented here could be used to extend such games, allowing for deployment of randomized surveillance schemes by Authorities, who would be choosing 
the best surveillance patterns $K_i$ along with the switching matrix $\Lambda.$ 

We also list several methodological extensions that will be useful to develop in the future.  First, all considered test problems were based on isotropic dynamics and running costs, but it would be easy to treat the general case by using different numerical schemes for HJB equations and quasi-variational inequalities.  
Additionally, our current use of time-explicit discretization might result in many time slices due to the CFL stability condition.  This can become a computational bottleneck, particularly for problems with ``fast layers'' (i.e., the speed $f$ might be large on a small subset of $\Omega$ only.)  But this can be avoided by using either time-implicit Eulerian schemes \cite{vladimirsky2013fast} or semi-Lagrangian discretizations \cite{falconeferretti2014}.  Second, it should not be hard to consider other mode-observation schemes.  In this paper, we showed how to use observations available at predetermined times or on-demand.  But a similar approach could be adapted to mode-observations occurring at random (driven by a non-homogeneous Poisson process) and to random terminations in infinite-horizon PDMPs.  We hope that the latter will have broad applications beyond path planning.
E.g., in modeling the economic impact of rare disasters \cite{haurie2006disaster}, some of the incremental changes (mode transitions) may not be immediately observed despite increasing the likelihood of disasters in the near future.
\update{Third, it will be useful to extend our framework to treat any combination of unobserved and immediately observed mode transitions.  The Mars rover problem considered in \S \ref{sec:rover} can be reinterpreted as  a simple example of this type since the ``termination''/breakdown is really just an immediately observed transition to the new mode with reduced capabilities.}

While our Assumption \ref{a:const_switching_rates} is very helpful for computational efficiency, it rules out all problems where the planner's choices might affect their belief about the current state.  E.g.,  it does not allow treating a planetary rover problem from Section \ref{sec:rover} if the rate of transitions into the breakdown-prone Mode 2 depends on roughness of the traversed terrain.  
Relaxing this assumption (possibly via some hybrid/reduced form of belief programming \update{or in combination with sampling-based methods}) will be an interesting and challenging direction for future work.   Other more challenging extensions include treating inexact mode observations, uncertainty in mode-switching rates, and piecewise-deterministic differential games (with a possible information asymmetry about the current mode).

\addtolength{\textheight}{-0.0cm}

%


\noindent
{\bf Acknowledgements: } The authors are grateful to Natasha Patnaik and Nagaprasad Rudrapatna, 
who performed initial computational tests on approaches in Sections \ref{sec:models}, \ref{sec:finite}, and the ``no observations'' version of \ref{sec:indef}
during a summer REU program at Cornell University. 
This research was partially funded by the National Science Foundation (awards DMS-1645643 and DMS-2111522) 
as well as the Air Force Office of Scientific Research (award FA9550-22-1-0528).

\noindent
{\bf Conflict of interest statement: } On behalf of all authors, the corresponding author states that there is no conflict of interest. 

\bibliographystyle{plain}
\bibliography{refs}

\appendix


\section{Belief-state Programming with Observations}\label{app:noisy-obs}
The framework presented in Section \ref{sec:belief-prog} can be extended to accommodate (potentially noisy) mode observations. Here we outline a belief-state programming approach for that scenario. We assume that in addition to some initial distribution, the planner receives observations (which may be noisy or exact) of the current mode at $L$ known discrete times\footnote{The general belief-state programming approach can accommodate many mechanisms for receiving information about the mode, but we choose to focus on the case that the planner receives observations at discrete times to make the comparison easier to the occasionally observed problems we considered in this paper.} $0 < T_1 \leq \hdots \leq T_L < T.$
At each observation time $T_l$ the planner receives an observation $I(\mu(T_l))$, generated according to the likelihood function $\rho(i,j)$ where $\rho(i,j) = \Prob(I(i) = j)$.\footnote{If $\rho(i,j) = \delta_{ij}$ where $\delta_{ij}$ is the Kronecker delta function, this corresponds to exact mode observations.}
In this case, the evolution of the belief itself becomes a piecewise-deterministic process.
Between observations it will be deterministic and governed by equation \eqref{eq:belief-statedep-ode} as in the no observation case.
Immediately after an observation, the belief must be updated to account for the newly acquired information.
We compute the updated belief via an update operator $U$ based on Bayes theorem. 
For previous belief $\vb q$ and observed mode $j$, the $i$-th element of $U$ is defined as
\begin{align}
U_i(\vb q,j) &= \Prob(\mu = i \mid \vb q, I(\mu) = j) \\
&= \frac{\Prob(I(\mu) = j \mid \vb q, \mu = i)  \Prob(\mu = i \mid \vb q) }{\Prob(I(\mu) = j \mid \vb q) } \\
&= \frac{\rho(i,j) \vb q_i}{\sum_{k=1}^M \rho(k,j) \vb q_k}.
\end{align}
If $\vb b(T_l^-)$ is the belief immediately prior to the observation at $T_l$ and $\vb b(T_l^+)$ is the updated belief after the observation has occurred, then we have
\begin{equation}
\vb b(T_l^+) = U \lp\vb b(T_l^-), I(\mu(T_l))\rp.
\end{equation}

The value function $w^l(\vx, \vb q, t)$ can then be defined to encode the minimal expected cost if we start at a time $t \in [T_l, T_{l+1}]$ after the observation at $T_l$ but before the observation at $T_{l+1}$.  (As before, we use the convention $T_0 = 0$ and $T_{L+1} = T.$)
Starting with the continuous state $\vy(t) =\vx$ and the  
belief $\vb b(t) = \vb q$, we have
\begin{align}
	w^l(\vx, \vb q, t) &= \inf_{\veca(\cdot)} \left\{ \int_{t}^{T_{l+1}} \overline{K}\lp\vy(s), \vb b(s), \veca(\vy(s), \vb b(s), s) \rp ds \right.\\
	&+ \left. \sum_{i=1}^M \vb b_i(T_{l+1}) \sum_{j=1}^M \rho(i,j) w^{l+1}(\vx, U(\vb b(T_{l+1}^-),j), T_{l+1})\right \}
\end{align}
Each $w^l$ can be found by solving a PDE over $\Omega \times [T_l,T_{l+1}] \times \mathbb{Q}^{M-1}.$  I.e., 
\begin{align}
         - \frac{\partial w^l}{\partial t} &= 
         \min_{\veca \in A}\left\{\overline{K}(\vx, \vb q, t, \veca) + \mathbf{f}(\vx, t, \veca) \cdot \nabla_{\vx} w^l 
         + Q\lp \vx, \vb q, t, \veca \rp \cdot \nabla_{\vb q} w^l 
         \right\},\\
        w^l(\vx,\vb q,T_{l+1}) &= \sum_{i=1}^M \vb q_i \sum_{j=1}^M \rho(i,j) w^{l+1}(\vx, U(\vb q,j), T_{l+1}), 
        \label{eq:belief-pde-noisy}
\end{align}
with the alternate terminal condition $w^L(\vx,\vb q,T) = \overline{\psi}(\vx,\vb q)$ for the final time interval.
Here we use again the functions $\overline{K}, \, \overline{\psi},$ and $Q$ defined in Section \ref{sec:belief-prog}.
This problem can be solved numerically over $\Omega \times [0,T] \times \mathbb{Q}^{M-1}$ by computing sequentially $w^L, w^{L-1}, ..., w^1, w^0.$
But unfortunately, if the PDEs are discretized on a grid or a mesh, the cost of this process scales exponentially with the number of modes.

\section{Beliefs conditioned on non-termination: a direct ODE derivation}\label{app:b-ode-bayes}
We are interested in estimating the likelihood of a mode conditioned on not terminating in a prematurely terminated problem described in Section \ref{sec:rt}.
While the interpretation presented in equation \eqref{eq:b-rt} provides a closed form expression for the belief, it is also possible to derive a nonlinear system of ODEs governing the evolution of $\vb b(t)$ via Bayes theorem.
We present that argument here.
\begin{prop}
Let $b_i(t) = \Prob(\mu(t) = i \mid \vb q, \lnot \Xi(t))$, where $\Xi(t)$ is the event that the process terminates prematurely in the interval $[0,t]$, then $\vb b(t)$ satisfies the following system of nonlinear ODEs
\begin{align}
b'_i(t) &= \sum_{j\neq i} \left( \lambda_{ji} b_j(t) - \lambda_{ij} b_i(t) \right) + \sum_{j=1}^M \left( b_j(t) \left(\gamma_j - \gamma_i\right) b_i(t)\right) \\
\vb b(0) &= \vb q
\end{align}
\end{prop}

\begin{proof}
We first compute $b_i(t + \tau)$ for some small time $\tau$, using $\tilde{\bq}$ to denote the known $\vb b(t),$
\begin{align}
b_i(t+\tau) &= \Prob(\mu(t + \tau) = i \mid \vb b(0) = \bq, \lnot \Xi(t +\tau)) \\
& = \Prob(\mu(\tau) = i \mid \vb b(0) = \tilde{\bq}, \lnot \Xi(\tau)),
\end{align}
where to obtain the second line we use the memoryless property of the exponential random variable.
An application of Bayes theorem leads to
\begin{align}
b_i(t + \tau) = \frac{\Prob\left(\lnot \Xi(\tau) \mid \mu(\tau) = i, \vb b(0) = \tilde{\bq} \right) \Prob\left(\mu(\tau) = i \mid \vb b(0) = \tilde{\bq} \right)}
{\Prob\left(\lnot \Xi(\tau) \mid \vb b(0) = \tilde{\bq}\right)}. \label{eq:b-bayes}
\end{align}
The rest of the argument involves finding small-time expansions in $\tau$ for each term in \eqref{eq:b-bayes}.

\textbf{Denominator:} The denominator is computed using a straightforward approximation of the relevant integral,
\begin{align}
\Prob\left(\lnot \Xi(\tau) \mid \vb b(0) = \tilde{\bq}\right) &= 1 - \int_0^\tau \sum_{j = 1}^M b_j(t + s) \gamma_j e^{-\gamma_j s} ds = 1 - \tau \sum_{j=1}^M \tilde{q}_j \gamma_j + o(\tau) \label{eq:b-notT-given-b}.
\end{align}

\textbf{Numerator:} The second term in the numerator is estimated via an Euler step for equation \eqref{eq:belief-statedep-ode}:
\begin{align}
\Prob\left(\mu(\tau) = i \mid \vb b(0) = \tilde{\bq}\right) 
&= \tilde{q}_i + \tau \sum_{j\neq i} \left( \tilde{q}_j \lambda_{ji} - \tilde{q}_i \lambda_{ij} \right)+ o(\tau) \label{eq:b-mu-given-b}.
\end{align}
The remaining term, $\Prob\left(\lnot \Xi(\tau) \mid \mu(\tau) = i, \vb b(0) = \tilde{\bq}\right)$, is the most involved.
Using the Law of Total Probability, we expand this into the sum
\begin{align}
&\Prob\left(\lnot \Xi(\tau) \mid \mu(\tau) = i, \vb b(0) = \tilde{\bq}\right) \nonumber \\
\qquad &=\sum_{j = 1}^M \Prob\left(\lnot \Xi(\tau) \mid \mu(\tau) = i, \mu(0) = j\right) \Prob(\mu(0) = j \mid \mu(\tau) = i, \vb b(0) = \tilde{\bq}). \label{eq:b-notT-given-mu}
\end{align}
We now consider two cases to evaluate the terms within the sum: $i = j$ and $i \neq j$. 
Since $\tau$ is small 
we assume at most one mode switch occurs in $[0,\tau]$.

\textbf{Case 1:} $i = j$. In this case, no mode switches occur in $[0, \tau]$ so the probability of not terminating over this interval is determined solely by $\gamma_i$:
\begin{align}
\Prob\left(\lnot \Xi(\tau) \mid \mu(\tau = i), \mu(0) = i\right) 
= 1 - \tau \gamma_i + o(\tau).
\end{align}
We can compute the probability that $\mu(0) = i$ using Bayes theorem:
\begin{align}
\Prob(\mu(0) = i \mid \mu(\tau) = i, \vb b(0) = \tilde{\bq}) 
&= \frac{\left(1 - \tau \sum_{j\neq i} \lambda_{ij}  + o(\tau)\right) \tilde{q}_i}{\tilde{q}_i + \tau \sum_{j=1}^M \left(\tilde{q}_j \lambda_{ji} - \tilde{q}_i \lambda_{ij}\right) + o(\tau)} 
\end{align}
Overall, the term in \eqref{eq:b-notT-given-mu} corresponding to $j = i$ is given by
\begin{align}
\frac{\tilde{q}_i - \tau \sum_{j\neq i} \tilde{q}_i\lambda_{ij} - \tau \gamma_i \tilde{q}_i + o(\tau)}{\tilde{q}_i + \tau \sum_{j\neq i}  \left( \tilde{q}_j \lambda_{ji} - \tilde{q}_i \lambda_{ij}\right) + o(\tau)}
\end{align}

\textbf{Case 2:} $j\neq i$. If $i \neq j$, then exactly one mode switch has occurred.
Let $T_s$ denote the time of this switch. 
Conditioned on the switch occurring, $T_s \sim \text{Unif}([0,\tau])$, so
\begin{align}
\Prob\left(\lnot \Xi(\tau) \mid \mu(\tau) = i, \mu(0) = j\right) &= \frac{1}{\tau} \int_0^\tau \Prob\left(\lnot \Xi(\tau) \mid \mu(\tau) = i, \mu(0) = j, T_s = s\right) \,ds\\
&= 1 - \frac{1}{\tau} \int_0^\tau s \gamma_j + (\tau - s) \gamma_i + o(\tau) \,ds.
\end{align}
Here the factor of $\frac 1 \tau$ comes from the PDF of the uniform distribution over $[0,\tau]$ and we have again applied the Law of Total Probability.
Evaluating the final integral we obtain
\begin{align}
\frac{1}{\tau} \int_0^\tau s \gamma_j + (\tau - s) \gamma_i + o(\tau) ds = \frac{\tau}{2}\left(\gamma_j + \gamma_i\right) + o(\tau)
\end{align}
and thus $\Prob\left(\lnot \Xi(\tau) \mid \mu(\tau) = i, \mu(0) = j\right) = 1 - \frac{\tau}{2}\left(\gamma_j + \gamma_i\right) + o(\tau).$
All that remains is to compute $\Prob(\mu(0) = j \mid \mu(\tau) = i, \vb b(0) = \tilde{\bq})$, for which we again apply Bayes theorem:
\begin{align}
\Prob(\mu(0) = j \mid \mu(\tau) = i, \vb b(0) = \tilde{\bq}) 
&= \frac{\left(\tau \lambda_{ji}  + o(\tau)\right) \tilde{q}_j}{\tilde{q}_i + \tau \sum_{j=1}^M \left( \tilde{q}_j \lambda_{ji} - \tilde{q}_i \lambda_{ij} \right) + o(\tau)}.
\end{align}
Thus, the terms in \eqref{eq:b-notT-given-mu} corresponding to $j\neq i$ are given by
\begin{equation}
\frac{\tau \tilde{q}_j\lambda_{ji} + o(\tau)}{\tilde{q}_i + \tau \sum_{j=1}^M \left( \tilde{q}_j \lambda_{ji} - \tilde{q}_i \lambda_{ij} \right) + o(\tau)}.
\end{equation}

\textbf{Conclusion:} We can now rewrite equation \eqref{eq:b-notT-given-mu} as
\begin{align}
\Prob\left(\lnot \Xi(\tau) \mid \mu(\tau) = i, \vb b(0) = \tilde{\bq}\right) 
&= 1 - \frac{\tau \gamma_i \tilde{q}_i + o(\tau)}{\tilde{q}_i + \tau \sum_{j\neq i}  \left( \tilde{q}_j \lambda_{ji} - \tilde{q}_i \lambda_{ij}\right) + o(\tau)} . \label{eq:b-notT-given-mu-final}
\end{align}
Combining equations \eqref{eq:b-bayes}, \eqref{eq:b-notT-given-b}, \eqref{eq:b-mu-given-b}, and \eqref{eq:b-notT-given-mu-final} and substituting $\tilde{q}_i = b_i(t)$ yields
\begin{align}
b_i(t + \tau) 
&= b_i(t) + \tau \frac{\sum_{j\neq i} \left( b_j(t) \lambda_{ji} - b_i(t) \lambda_{ij} \right) + b_i(t) \sum_{j=1}^M b_j(t) \gamma_j - \gamma_i b_i(t) + O(\tau)}{ 1 - \tau \sum_{j=1}^M b_j(t) \gamma_j + o(\tau)}.
\end{align}
which leads to the desired expression for $b'_i(t)$,
\begin{align}
\lim_{\tau \to 0} \frac{b_i(t + \tau) - b_i(t)}{\tau} &= \lim_{\tau \to 0} \frac{\sum_{j\neq i} \left(b_j(t) \lambda_{ji} - b_i(t) \lambda_{ij} \right)+ \sum_{j=1}^M b_j(t)(\gamma_j - \gamma_i) b_i(t) + O(\tau)}{ 1 - \tau \sum_{j=1}^M b_j(t) \gamma_j + o(\tau)}\\
&= \sum_{j=1}^M \lp b_j(t) \lambda_{ji} - b_i(t) \lambda_{ij} \rp + \sum_{j=1}^M b_j(t)(\gamma_j - \gamma_i) b_i(t).
\end{align}
\end{proof}

\end{document}